\documentclass[final]{siamltex}
\usepackage{amsmath,amsfonts,amssymb}
\usepackage{latexsym}
\usepackage{epsfig,overpic}
\usepackage{commath}
\usepackage{mathtools}

\usepackage{subcaption}
\usepackage[utf8]{inputenc}
\usepackage{booktabs}
\usepackage{csvsimple}
\usepackage{colortbl}

\let\pa\partial

\newcommand{\R}{\mathbb R}

\newcommand{\bI}{\mathbf I}

\newcommand{\bP}{\mathbf P}

\newcommand{\bg}{\mathbf g}
\newcommand{\bn}{\mathbf n}

\newcommand{\bu}{\mathbf u}

\newcommand{\bw}{\mathbf w}
\newcommand{\bx}{\mathbf x}

\newcommand{\T}{\mathcal T}

\newcommand{\spa}{s}
\newcommand{\ti}{q}
\newcommand{\spati}{qs}
\newcommand{\mass}{\operatorname{mass}}
\newcommand{\surf}{\operatorname{surf}}
\newcommand{\nod}{\operatorname{nod}}

\newcommand{\bnu}{\boldsymbol{\nu}}
\newcommand{\Div}{\operatorname{div}}

\newcommand{\tr}{\mathop{\rm tr}}
\newcommand{\id}{\mathop{\rm id}\nolimits}

\newcommand\restrict[1]{\raisebox{-.5ex}{$|$}_{#1}}
\newcommand{\lin}{\operatorname{lin}}
\newcommand{\dist}{\operatorname{dist}}
\newcommand{\reg}{\operatorname{reg}}
\newcommand{\triplenorm}{\ensuremath{| \! | \! |}}
\newcommand*{\ddt}[3][]{\ensuremath{\frac{\partial^{#1} #2}{\partial #3}}}

\newtheorem{remark}{Remark}

\usepackage{tikz}
\usetikzlibrary{intersections}
\usetikzlibrary{arrows}
\usetikzlibrary{patterns}
\usetikzlibrary{calc} 
\usetikzlibrary{quotes,babel,angles}
\usetikzlibrary{decorations.pathreplacing}
\usetikzlibrary{pgfplots.groupplots}
\usepackage{xintexpr}
\usepackage{pgfplots}

\usepackage{hyperref}
\usepackage{cleveref}

\crefname{equation}{}{}

\pgfplotscreateplotcyclelist{haukes_cycle_list}{
	blue,every mark/.append style={fill=blue},mark=square*\\
	red,every mark/.append style={fill=red},mark=triangle*\\
	green,every mark/.append style={fill=green},mark=diamond*\\
	orange,every mark/.append style={fill=orange},mark=pentagon*\\
	purple,every mark/.append style={fill=purple},mark=x\\
	brown,every mark/.append style={fill=brown},mark=*\\
}
\pgfplotsset{cycle list name=haukes_cycle_list}
\pgfplotsset{compat=1.15}
\pgfplotsset{	every non boxed x axis/.style={} }

\begin{document}
	
	\title{An accurate and robust Eulerian finite element method for  partial differential equations on evolving surfaces} 
	\author{Hauke Sass\thanks{Institut f\"ur
			Geometrie und Praktische  Mathematik, RWTH-Aachen University, D-52056 Aachen,
			Germany; email: {\tt sass@igpm.rwth-aachen.de}} \and Arnold Reusken\thanks{Institut f\"ur
			Geometrie und Praktische  Mathematik, RWTH-Aachen University, D-52056 Aachen,
			Germany; email: {\tt reusken@igpm.rwth-aachen.de}} }

	\maketitle
	
	
	\begin{abstract}
		In this paper we present a new Eulerian finite element method for the discretization of scalar partial differential equations on evolving surfaces. In this method we use the restriction of standard space-time finite element spaces on a fixed bulk mesh to the space-time surface. The structure of the method is such that it naturally fits to a level set representation of the evolving surface. The higher order version of the method is based on a space-time variant of a  mesh deformation that has been developed in the literature for stationary surfaces.   The discretization method that we present is of (optimal) higher order accuracy for smoothly varying surfaces with sufficiently smooth solutions. Without any modifications the method can be used for the discretization of problems with topological singularities. A numerical study demonstrates both the higher order accuracy for smooth cases and the robustness with respect to toplogical singularities.
	\end{abstract}
	

	
	\section{Introduction} \label{secIntro}
	Motivated by applications in, for example, multiphase flow and computational biology, in the past decade there has been a strongly growing interest in the development and analysis of numerical methods for partial differential equations (PDEs) on (evolving) surfaces. For an overview of recent developments in this research field we refer to \cite{DEreview,Bonito2019}. The main focus of these two review papers is on numerical methods for scalar PDEs on smooth \emph{stationary} surfaces. There are, however, also many papers in which discretization methods for scalar PDEs on \emph{evolving} surfaces are introduced or analyzed. The present paper contributes to this field of finite element methods for scalar PDEs on evolving surfaces.
	
	In this research field several interesting  methods have already been treated. The probably most  prominent method is the evolving surface finite element method (ESFEM) \cite{DEreview,DziukElliot2013a,kov2018}. This elegant method is very popular and a complete error analysis has been developed. A key feature of this method is that it is based on a Lagrangian approach in which a triangulation of the initial surface and a corresponding finite element space are transported along the flow lines of a (given) velocity field. This makes the method very attractive for problems with smoothly and slowly varying surfaces. On the other hand, it is well-known that such Lagrangian techniques have drawbacks if the geometry of the surface is strongly varying on the relevant time scales or if there are topological singularities, e.g., merging or splitting phenomena. We are not aware of any paper in which the ESFEM is applied to a surface PDE with a topological singularity.  Another class of finite element methods is based on an Eulerian approach and has been introduced in \cite{olshanskii2014eulerian,olshanskii2014error}. These so-called trace (or cut) finite element methods (TraceFEM or CutFEM) are based on restriction of standard space-time finite element spaces on a fixed bulk mesh to the space-time surface. The structure of TraceFEM is such that it naturally fits to a level set representation of the evolving surface. Such level set representations are known to be very attractive for problems with strongly varying surface geometries or with topological singularities.  Hybrid versions of this Eulerian approach have been introduced in \cite{hansbo2015characteristic,LOX_SIAM_2017}. The  method treated in \cite{hansbo2015characteristic} is based on a
	characteristic-Galerkin, i.e., Lagrangian, formulation combined with a piecewise linear CutFEM in space.  The fully Eulerian method introduced in \cite{LOX_SIAM_2017} is based on a combination of TraceFEM in space and  finite differences for the time discretization. In the numerical experiments presented in these papers the methods are applied to problems with smoothly and slowly varying surfaces.  In \cite{LOX_SIAM_2017} not only linear finite elements but also higher order ones are considered. 
	
	In the present paper we return to the Eulerian space-time technique introduced in \cite{olshanskii2014eulerian,olshanskii2014error}. We extend the method treated  in those papers in the following two directions.   Firstly, we will introduce \emph{a higher order} space-time TraceFEM. In \cite{olshanskii2014eulerian,olshanskii2014error} only piecewise linears (in space and time) are used and the issue of geometry approximation is not considered,
	The development of the higher order technique is based on the parametric TraceFEM introduced for stationary surfaces in \cite{lehrenfeld2016high} and analyzed in \cite{grande2018analysis}. The unfitted space-time finite element technique that we  propose in this paper is based on very similar ideas as used for scalar PDEs on a moving domain in \cite{HLP2022}.
	Secondly, in the method presented in this paper we include  a space-time variant of the so-called \emph{normal derivative volume stabilization} that has been developed for TraceFEM applied to PDEs on stationary surfaces, cf. \cite{burman2016cutb,grande2018analysis}.  This stabilization is important, in particular for higher order finite elements, to control the conditioning of the stiffness matrix. 
	
	Using these extensions results in a method with the following features. The method is of Eulerian type, using a fixed bulk triangulation and corresponding standard  space-time (tensor product)  finite element spaces. As input for the method one needs (only) a level set representation of the evolving surface. If the level set function approximation is sufficiently accurate (made precise further on in the paper)  and higher order bulk space-time finite element spaces are used, the resulting surface PDE discretization method is also of (optimal) higher order accuracy for smoothly varying surfaces with sufficiently smooth solutions. Without any modifications the method can be used for the discretization of problems with topological singularities.  We are not aware of any paper in which this combination of properties, i.e., higher order for smooth cases and robustness w.r.t. topological singularities, is demonstrated for a  discretization method for surface PDEs. 
	
	In this paper we introduce the method, discuss important properties of the method and present an extensive numerical study that demonstrates the features mentioned above. We do not include a stability or error analysis. An error analysis for the lowest order case (piecewise planar surface approximation and $P_1$ finite elements) is presented in \cite{Diss_Sass} and will be further treated in a forthcoming paper. Implementation of the method is not straightforward due to the use of a space-time framework and a particular space-time mesh transformation (for the higher order case). An implementation, however, is available in the open source software \texttt{ngsxfem} \cite{ngsxfem}, an add-on package to the finite element library \texttt{NGSolve}$\backslash$\texttt{Netgen}, see \cite{Schoeberl1997,ngsolve}.  
	
	The remainder of the paper is organized as follows. In \Cref{sect1} we explain the parabolic model problem that we consider and a well-posed space-time variational formulation of this problem.  The space-time TraceFEM is explained in \Cref{sectTraceFEM}. Several aspects of the method are treated. We start with the level set representation of the evolving surface, needed as input for the discretization method. We then outline the space-time mesh deformation that we use for the higher order variant of our discretization method. Several implementation aspects such as, for example, the efficient quadrature applied to the space-time integrals are addressed. In \Cref{section_numerical_experiments} we present results of numerical experiments, both for   examples with smoothly evolving surfaces and surfaces with topological singularities. These results illustrate the higher order convergence for the smooth case and the robustness of the method with respect to topological singularities.

\section{The  model problem}\label{sect1}
Let $\Omega\subset \R^3$, $T>0$  and $\mathbf{w}:\Omega\times [0,T]\rightarrow \R^3$ a given  smooth velocity field. We consider a closed, connected and orientable $C^2$-hypersurface $\Gamma(t)\subset\Omega$, $t\in [0,T]$, that is is advected $\mathbf{w}$. By $\bn$ we denote the outward pointing unit outer normal on $\Gamma(t)$. The smooth space-time manifold is denoted by
\begin{equation*}
	S\coloneqq \bigcup_{t\in [0,T]}\left(\Gamma(t)\times \{t\}\right)\subset\R^4.
\end{equation*}
We recall standard surface differential operators. For a smooth scalar function $g$ on $\Gamma(t)$, with $t$ fixed, a smooth extension to a tubular neighborhood of $\Gamma(t)$ is denoted by $g^e$. Its spatial and space-time gradients  are denoted by $\nabla$ and $\nabla_{(\bx,t)}$, respectively.
The tangential gradient is defined as $\nabla_{\Gamma} g\coloneqq (I-\bn\bn^\intercal)\nabla g^e$. For a vector valued function $\bg$ we define the surface divergence $\Div_\Gamma \bg \coloneqq \tr\big((I-\bn\bn^\intercal)\nabla \bg^e \big)$. The surface Laplace-Beltrami operator is $\Delta_{\Gamma}\coloneqq \Div_{\Gamma}\nabla_{\Gamma}$.
We consider a basic model for convection and molecular diffusion of a surface species. The conservation of mass principle combined with Fick's law for the diffusive flux leads to the parabolic surface partial differential equation
\begin{equation} \label{surfactant1}  \begin{split}
	\dot{u}+ u\Div_{\Gamma}\bw-\mu_d \Delta_{\Gamma} u&=f\quad  \text{on } \Gamma(t),\quad t\in (0,T],
	\\
	u(\cdot,0)&=0\quad \text{on } \Gamma(0), 
\end{split} \end{equation}
for the scalar unknown function $u=u(\bx,t)$.
Here $\mu_d >0$ denotes the constant diffusion coefficient and $\dot{u}$ the material derivative along the velocity field $\bw$, which can be expressed as $\dot{u} = \frac{\pa u^e}{\pa t}$ + $\bw\cdot \nabla u^e$. The right-hand side $f$ must satisfy the consistency condition $\int_{\Gamma(t)}f\dif s=0$, $t\in [0,T]$. In the remainder we assume that this condition is satisfied. 
 For details on the derivation of \cref{surfactant1} we refer to \cite{James04,DEreview,GReusken2011}. 

We outline a weak variational space-time formulation of \cref{surfactant1}, introduced in \cite{olshanskii2014eulerian}, which is the basis for the finite element method treated in \Cref{sectTraceFEM}. For this we first define suitable Hilbert spaces that can be used for a discontinuous in time Galerkin method.
For $N\in \mathbb{N}$, let the time interval $[0,T]$ be partitioned into smaller time intervals $I_n\coloneqq [t_{n-1},t_n]$, $n\in \{1,\dots,N\}$, where $0=t_0<t_1<\dots<t_N=T$. The corresponding space-time surface is denoted by $S^n\coloneqq \bigcup_{t\in I_n}( \Gamma(t)\times \{t\})$. We define the Hilbert spaces
\begin{alignat*}{2}
	H & \coloneqq \left\lbrace v\in L^2(S): \norm{\nabla_{\Gamma} v}_{L^2(S)}<\infty \right\rbrace,&&
	\norm{v}^2_H\coloneqq \int_0^T\norm{v}^2_{H^1(\Gamma(t))}\dif t
 \\
	H_n & \coloneqq\left\lbrace v\in L^2(S^n): \norm{\nabla_{\Gamma} v}_{L^2(S^n)}<\infty \right\rbrace,\,\,\,\,\,&&
	\norm{v}^2_{H_n}\coloneqq \int_{t_{n-1}}^{t_n}\norm{v}^2_{H^1(\Gamma(t))}\dif t, ~ 1 \leq n \leq N.
\end{alignat*}
The weak material derivative is based on  the linear functional 
\begin{equation}
\langle \dot{u}, \psi\rangle_n\coloneqq -\int_{t_{n-1}}^{t_n}\int_{\Gamma(t)}u\dot{\psi}+v\psi\Div_\Gamma \bw \dif s\dif t \quad \text{for all } \psi \in C_0^1(S^n).\label{def_dot_u}
\end{equation}
All functions $u\in H_n$ whose weak material derivative $\dot{u}$ define a bounded linear functional on $H_n$ form the subspace $W_n\subset H_n$, i.e.
\begin{equation*}
	W_n\coloneqq \left\lbrace u\in H_n:\dot{u}\in H'_n\right\rbrace ,\quad \norm{u}^2_{W_n}\coloneqq  \norm{v}_{H_n}^2 +\norm{\dot{u}}_{H_n'}^2. 
\end{equation*}
For $v\in H$ we define $v_n\in H_n$ as $v_n\coloneqq v\restrict{S^n}$, $ n=1,\dots,N$. The space that we use for the space-time  variational formulation is 
\begin{equation}
	W^b\coloneqq \left\lbrace v\in H: v\restrict{S^n}\in W_n, n=1,\dots,N\right\rbrace, \quad \norm{v}^2_{W^b}\coloneqq \sum_{n=1}^N\norm{v_n}_{W_n}^2.\label{brokenspace}
\end{equation}
This space is a Hilbert space. 
For further properties  of these spaces we refer to \cite{olshanskii2014eulerian}. For the discontinuous Galerkin method we need the usual jump operator to relate values between subsequent time slabs. For $v\in W^b$ we define
	\begin{alignat*}{2}
		v_+^n&\coloneqq v_+(\cdot ,t_n)=\lim_{\eta\searrow 0}v(\cdot, t_n+\eta)\in L^2(\Gamma(t_n)),\quad &&n=0,\dots,N-1,
		\\ v_-^n&\coloneqq v_-(\cdot ,t_n)=\lim_{\eta\searrow 0}v(\cdot, t_n-\eta)\in L^2(\Gamma(t_n)),  &&n=1,\dots,N, ~~v_-^0\coloneqq 0, \\
		[v]^n &\coloneqq v_+^n-v_-^n\in L^2(\Gamma(t_n)) ,&&n=0,\dots,N-1.
	\end{alignat*}
	Note that $[v]^0=v_+^0$ since we have by definition $v_-^0= 0$. On $\Gamma(t)$, $t\in[0,T]$, the $L^2$ scalar product is denoted by $(\cdot,\cdot)_{\Gamma(t)}$. Space-time integrals and the surface integral on $S$ are related by
	\begin{equation}
	\int_0^T\int_{\Gamma(t)}g\dif s\dif t=\int_S \frac{g}{\alpha}\dif \sigma, \quad \alpha\coloneqq \sqrt{1+(\bw\cdot\bn)^2}.
	\label{transformationsformelMAIN}
\end{equation}
	With these preparations we can introduce a space-time variational problem consistent with \cref{surfactant1}:
	For a given $f\in L^2(S)$ determine $u\in W^b$
	such that
\begin{align}
	\begin{aligned}
		\sum_{n=1}^N B_n(u,v) &\coloneqq \sum_{n=1}^N\Big[\langle \dot{u}_n,v_n\rangle_n+ \int_{S^n}\frac{1}{\alpha}\big( uv\Div_{\Gamma}\mathbf{w}
		+\mu_d\nabla_{\Gamma}u\cdot \nabla_{\Gamma}v\big)\dif \sigma  \\ & \quad +\left([u]^{n-1},v_+^{n-1}\right)_{\Gamma(t_{n-1})}\Big]=\int_S\frac{1}{\alpha}fv\dif s\dif t\label{prob_cont} \quad \text{for all}~v\in W^b.
	\end{aligned}
\end{align}
Well-posedness of this problem is shown in  \cite[Theorem 5.3]{olshanskii2014eulerian}.
In the formulation \cref{prob_cont} the material derivative is applied to the trial function $u_n$. By partial integration this derivative can be shifted to the test function $v_n$. As is well-known from the literature on DG methods this then  leads do a different discretization. We will also address this issue for the space-time finite element method that we propose. Hence, besides the formulation in \cref{prob_cont} we also introduce two other natural formulations, namely one with the material derivative on the test function and an antisymmetric variant. For deriving these we use the partial integration rule
\begin{align}
	\begin{aligned} 
	 & \int_{S^n}\frac{1}{\alpha}\dot{u}v\dif \sigma \\ 
	 &\qquad =-\int_{S^n}\frac{1}{\alpha}\left(u\dot{v}+uv\Div_{\Gamma}\bw\right)\dif \sigma +\left(u_-^{n},v_-^{n}\right)_{\Gamma(t_{n})}-\left(u_+^{n-1},v_+^{n-1}\right)_{\Gamma(t_{n-1})},\label{pi_cont}
	 \end{aligned}
 \end{align}
for $u,v\in C^1(S^n)$, which directly follows from the Leibniz rule and \cref{transformationsformelMAIN}.
%
Note that for $u,v\in H^1(S^n)$  the weak material derivative can be expressed as ${\langle \dot{u}, v\rangle_n = \int_{S^n}\frac{1}{\alpha}\dot{u} v \dif \sigma}$. For arbitrary $\beta \in \mathbb{R}$ and $u,v\in H^1(S^n)$ the following holds:
\begin{equation} \label{betaformel}
 \begin{split}
		 & B_n(u,v)=\int_{S^n}\frac{1}{\alpha}\left((1-\beta)\dot{u}v- \beta u \dot{v} +(1-\beta)uv\Div_{\Gamma}\mathbf{w}+\mu_d\nabla_{\Gamma}u\cdot \nabla_{\Gamma}v\right)\dif\sigma\\
		& \quad + \beta \left(u_-^n,v_-^n\right)_{\Gamma(t_{n})} +(1-\beta)\left(u_+^{n-1},v_+^{n-1}\right)_{\Gamma(t_{n-1})}-\left(u_-^{n-1},v_+^{n-1}\right)_{\Gamma(t_{n-1})}.
\end{split}
\end{equation}
This follows from \eqref{pi_cont} using the following argument.  For fixed $u,v\in H^1(S^n)$ we take the derivative of the right-hand side \eqref{betaformel} with respect to $\beta$, which yields
\begin{equation} \label{contbeta}
 -\int_{S^n}\frac{1}{\alpha}\left( \dot{u}v + u \dot{v} +uv\Div_{\Gamma}\mathbf{w}\right)\dif\sigma
		  +  \left(u_-^n,v_-^n\right)_{\Gamma(t_{n})}-\left(u_+^{n-1},v_+^{n-1}\right)_{\Gamma(t_{n-1})},
\end{equation}
which equals zero due to \eqref{pi_cont}. Hence the right-hand side in \eqref{betaformel} does not depend on $\beta$. For $\beta =0$ formula \eqref{betaformel} is correct, which follows by comparing with the definition in \eqref{prob_cont}. Particularly interesting, in view of finite element discretization methods, are the choices $\beta=0$, $\beta=\tfrac12$ and $\beta=1$. 
When one replaces the smooth surface $S_n$ by an approximate Lipschitz surface, as we will do in our finite element discretization below, the resulting three formulations are no longer equivalent.
The reason for this is that the  identities \cref{transformationsformelMAIN,pi_cont} do not hold on Lipschitz surfaces. 
\section{Space-time TraceFEM} \label{sectTraceFEM} In this section we explain the space-time finite element method that we use for discretization of \eqref{prob_cont}. Key ingredients of the method are a level set representation of the evolving surface and  a parametric mapping based on the level set function (approximation) that is used to approximate the surface and to define the (higher order) space-time finite element spaces.  In this approach it is essential that the evolving surface $\Gamma(t)$ is represented as the zero level of an evolving level set function, denoted by $\phi(\bx,t)$. This level set function is approximated by a (higher order) finite element function, denoted by $\phi_h$,  on a fixed tetrahedral mesh, cf. \Cref{sectLevelset}. Evaluating integrals on the zero level of a higher order polynomial is computationally (very) expensive. To avoid this, we use an approach similar  to the classical parametric finite element technique.  First  a spatial mesh deformation, that is based on $\phi_h(\cdot,t)$ for fixed $t$, is constructed, cf. \Cref{sectspatialtrans}. Using a tensor product approach this then results in a mapping that deforms of the space-time mesh, cf. \Cref{sectspacetimetrans}. This space-time mesh deformation mapping, that can be efficiently evaluated, has a twofold purpose. It is used to construct a parametrization of a higher order accurate space-time surface approximation that can be efficiently evaluated. This allows efficient quadrature for integrals over higher order space-time surface approximations, cf. \Cref{secimplementation}. The mapping is also used for constructing parametric finite element spaces that have higher order approximation accuracy, cf. \Cref{secmethod}. 
\subsection{Level set representation} \label{sectLevelset}
  We assume that the evolving surface  $\Gamma(t)$ is represented as the zero level of a smooth level set function ${\phi:\Omega\times [0,T]\rightarrow \R}$, i.e.
\begin{equation*}
\Gamma(t)=\left\lbrace \bx\in\Omega: \phi(\bx,t)=0 \right\rbrace, \quad t\in [0,T].
\end{equation*}
We assume  standard properties of a level set function, i.e. for all $(\bx,t)$ in a neighborhood $U$ of $S$,
	\begin{equation}
		\norm{\nabla\phi(\bx,t)}\sim 1, \quad \norm{D^2\phi(\bx,t)}\lesssim 1\label{phi_assumptions}
	\end{equation}
	hold. We define the space-time cylinders
\begin{equation*}
	Q\coloneqq \Omega\times [0,T]\subset \mathbb{R}^{4},\quad Q_n\coloneqq \Omega \times I_n, \quad n=1,\dots,N.
\end{equation*}
Let $\T$ be an element of a family $\{{\T}_h\}_{h>0}$ of shape regular tetrahedral triangulations of $\Omega$,  The triangulation $Q_{h,n}\coloneqq \T\times I_n$ divides the time slab $Q_n$ into space-time prismatic elements. The  space-time triangulation of $Q$ is denoted by ${Q_h\coloneqq \bigcup_{n=1}^N Q_{h,n}}$. 
Let $V_h^{m}$, $m\in \mathbb{N}$, be the standard  $H^1(\Omega)$-conforming finite element space on the  triangulation $\T$. For $m_{\spa},m_{\ti}\in \mathbb{N}$ a space-time finite element product space is given  by 
\begin{equation}
	V_h^{m_{\spa}, m_{\ti}}\coloneqq \big\{ v:Q \rightarrow \R : v(\bx,t)=\sum_{i=0}^{m_{\ti}}t^i v_i(\bx), v_i \in V_h^{m_{\spa}}, (\bx,t)\in Q_n, n=1,\dots,N\big\}.\label{st_fe_space}
\end{equation}
Note that functions from $V_h^{m_{\spa},m_{\ti}}$ are continuous is space but may be discontinuous in time  at the time interval boundaries $t_n$. For the geometry approximation we use a function from this space that is continuous both in space and time and has degree $k_{g,\spa}$ in $\bx$ and $k_{g,\ti}$ in $t$. More specifically, we assume
	  an approximation ${\phi_h\in V_h^{k_{g,\spa},k_{g,\ti}}}$ of the level set function $\phi$ with ${\phi_h(\bx,\cdot)\in C^0([0,T])}$, $\bx \in \Omega$ that satisfies
	\begin{equation}
		\max_{\substack{K \in \T\\ n\in \{1,\dots,N\}}}\abs{\phi_h-\phi}_{W^{m,\infty}(\left(K\times I_n\right)\cap U)}\lesssim h^{k_g+1-m},\quad 0\leq m\leq k_g+1,\label{phi_assumption_prev}
	\end{equation}
with ${k_g\coloneqq \min\{k_{g,\spa},k_{g,\ti}\}}$. Such a $\phi_h$ can be obtained by 
interpolation of $\phi$, if the latter is available. In applications, $\phi$ may be determined implicitly by a level set equation. To satisfy \eqref{phi_assumption_prev} one has to solve this equation  sufficiently accurate. The continuity in time of $\phi_h$ later ensures that the piecewise planar reference space-time surface is a connected Lipschitz manifold. 
\subsection{Spatial mesh deformation} \label{sectspatialtrans}
We briefly explain the main idea of the parametric mapping introduced in \cite{lehrenfeld2016high} to construct higher order \emph{un}fitted finite element approximations.  It is used to define a corresponding space-time mapping in the next section. 
By $\mathcal{I}_{m}$ we denote the spatial nodal $\mathcal{P}^{m}$-interpolation operator
\begin{equation}
	\mathcal{I}_{m}:C^0(\Omega)\rightarrow V_h^{m},\quad m\in \mathbb{N}.\label{interpol_space_nodal}
\end{equation}
 Let the spatially piecewise linear nodal interpolation $\hat{\phi}_h\in V_h^{1,k_{g,t}}$ be defined by 
\begin{equation}
	\hat{\phi}_h(\cdot,t)\coloneqq \mathcal{I}_1\phi_h(\cdot,t)\in V_h^1 \quad \text{for all } t\in [0,T],\label{phi_h_dach_definition}
\end{equation}
and  its corresponding piecewise planar zero level  at time $t$, which is easy to compute:
\begin{equation*}
	\Gamma_{\lin}(t)\coloneqq \left\lbrace\bx\in \Omega : \hat{\phi}_h(\bx,t)=0\right\rbrace.
\end{equation*}
By construction, $\Gamma_{\lin}(t)$ is only a second order approximation to $\Gamma(t)$. We consider all elements that are cut by the piecewise linear surface at \emph{any point in time within one time slab}, i.e.,
\begin{equation*}
	\T_n^\Gamma \coloneqq \left\lbrace K\in \mathcal{T}:  \mathrm{meas}_2((K\times I_n)\cap \Gamma_{\lin}(t))>0 \text{ for any }t\in I_n\right\rbrace,\quad n\in \{1,\dots,N\}.
\end{equation*}
The corresponding domain is denoted by $\Omega_n^\Gamma \coloneqq \{\bx\in K: K\in \T_n^{\Gamma}\}$. On this domain a mapping $\Theta_{h,t}^n \in  \big(V_h^{k_{g,\spa}})^3$ is defined that \emph{depends (only) on $\phi_h$} and deforms all elements of the triangulation $\T_n^{\Gamma}$.  The image of this mapping restricted to $\Gamma_{\lin}(t)$ defines a higher order approximation of $\Gamma(t)$ such that  $\dist_3(\Gamma(t),\Theta_{h,t}^n(\Gamma_{\lin}(t)))\lesssim h^{k_{g,\spa}+1}$ for all $n=1,\dots,N$, $t\in I_n$.  We refer to \cite{lehrenfeld2016high,grande2018analysis} for precise definitions and an analysis of approximation properties of this mapping.  In these papers this mapping only deforms elements that are cut by $\Gamma_{\lin}(t)$, which, however, is not essential for the general construction. In \cite{HLP2022,preuss2018higher} the extended deformation, i.e., the  one applied to all tetrahedra in  $\T_n^{\Gamma}$ is explained.  

\subsection{Space-time mesh deformation} \label{sectspacetimetrans}
We explain the space-time mesh deformation on an arbitrary time slab $Q_n$, $n=1,\dots,N$, also used in \cite{HLP2022,preuss2018higher}. Let $\tau_m^n \in I_n$ be discrete points and $\mathcal{X}_{\tau_m^n} \in \mathcal{P}^{k_{g,\ti}}$, $m=0,\dots, k_{g,\ti}$, the corresponding finite element (e.g. nodal) basis functions. On the time slab $Q_n$ the finite space-time finite element function 
$\phi_h \in V_h^{k_{g,\spa},k_{g,\ti}}$ can be represented as
\begin{equation*}
	\phi_h(\bx,t)=\sum_{m=0}^{k_{g,\ti}}\mathcal{X}_{\tau_m^n }(t)\phi_m(\bx), \quad \phi_m\in V_h^{k_{g,\spa}},
\end{equation*}
 We define the triangulation of space-time prisms that are intersected by $\Gamma_{\lin}(t)$ at any point in time $t\in I_n$ as $Q_{h,n}^S\coloneqq\T_n^{\Gamma}\times I_n$, $n= 1,\dots,N$, and $Q_{h}^S:= \cup_{n=1}^N Q_{h,n}^S$. The subdomains formed by these triangulations are denoted by $Q^S_n$ and $Q^S$. The space-time mesh transformation ${\Theta_{h}^{n}\in (V_h^{k_{g,\spa},k_{g,\ti}}\restrict{Q^S_{h,n}})^4}$ is defined by
	\begin{equation}
		\Theta_{h}^{n}(\bx,t) \coloneqq	 \big(\sum_{m=0}^{k_{g,\ti}}\mathcal{X}_{\tau^n_m}(t)\Theta_{h,\tau^n_m}^{n}(\bx),\, t\big)	=:\big(\Theta_{h,\spa}^{n}(\bx,\, t), t\big), \quad (\bx,t)\in Q_{n}^S.\label{Theta_Def}
	\end{equation}
	The mapping $\Theta_{h,\spa}^{n}$ denotes the spatial part of $\Theta_{h}^{n}$. For $n\in \{1,\dots,N\}$ we define the discrete space-time manifolds
	\begin{alignat}{2}
		S^n_{\lin}&\coloneqq \left\lbrace (\bx,t)\in Q_n: \hat{\phi}_h(\bx,t)=0\right\rbrace,  &S_{\lin}&\coloneqq \bigcup_{n=1}^N S_{\lin}^n,\label{slin_def}\\
		S_h^n&\coloneqq \left\lbrace (\bx,t)\in Q_n: (\bx,t) \in \Theta_h^n(S^n_{\lin})\right\rbrace, \quad &S_{h}&\coloneqq \bigcup_{n=1}^N S_h^n\label{sh_def}
	\end{alignat}
	and the corresponding time slices of $S_h^n$ are denoted by
	\begin{equation}
		\Gamma_h^n(t)\coloneqq \left\lbrace\bx\in\R^3: (\bx,t)\in S_h^n\right\rbrace,\quad t\in I_n, \quad n=1,\dots,N.\label{gamma_h_def}
	\end{equation}
\input{theta_picture}
We introduce notation for the curved space-time surface triangulation. Let the space-time surface triangulation $\mathcal{T}_{S_h^n}$ be the set of smooth three-dimensional manifolds that $S_{h}^n$ consists of, i.e. 
\[ {\T_{S_h^n}\coloneqq \{\Theta_h^n(P)\cap S_h^n: P\in Q_{h,n}\}}.
\]
 Then, we have ${S_h^n=\bigcup_{K_S\in \mathcal{T}_{S_h^n}} K_S}$. Let $\T_{S_h}\coloneqq \bigcup_{n=1}^N\T_{S_h^n}$. Note that $K_S\in \T_{S_h}$ is a three-dimensional curved polytope that can be partitioned into curved tetrahedra. This deformed triangulation $\T_{S_h}$ is not necessarily shape regular, in the sense that internal angles can be very large and the size of neighboring elements may vary strongly. This is due to the fact that $S_{\lin}$ can have arbitrary cuts  with the background bulk mesh $Q_h$. An illustration of the two-dimensional analogon of $\T_{S_h^n}$ is given in \Cref{theta_picture}.
 
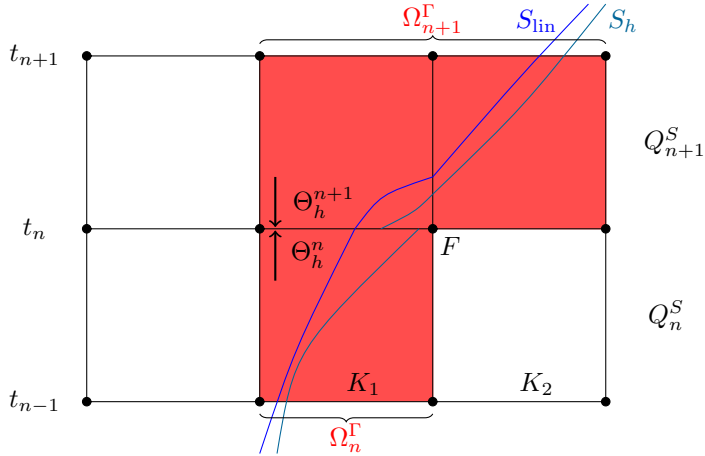
\begin{figure}[!htbp]
	\centering
	\begin{tikzpicture}[scale = 2.3]
		\newdimen\XCoord
		\newdimen\YCoord
		\newcommand{\shgrid}{dotted}
		\newcommand{\ptsize}{0.7pt}
		\newcommand{\drawmesh}[4]{
			\filldraw[#3, #4] (#1+0,#2+0) circle (\ptsize) -- (#1+1,#2+0) circle (\ptsize) -- (#1+2,#2+0) circle (\ptsize) -- (#1+3,#2+0) circle (\ptsize) ;
			\filldraw[#3, #4] (#1+0,#2+1) circle (\ptsize) -- (#1+1,#2+1) circle (\ptsize) -- (#1+2,#2+1) circle (\ptsize) -- (#1+3,#2+1) circle (\ptsize) ;
			\filldraw[#3, #4] (#1+0,#2+2) circle (\ptsize) -- (#1+1,#2+2) circle (\ptsize) -- (#1+2,#2+2) circle (\ptsize) -- (#1+3,#2+2) circle (\ptsize)  ;
			\filldraw[#3, #4] (#1+0,#2+0) -- (#1+0,#2+1) -- (#1+0,#2+2) ;
			\filldraw[#3, #4] (#1+1,#2+0) -- (#1+1,#2+1) -- (#1+1,#2+2);
			\filldraw[#3, #4] (#1+2,#2+0) -- (#1+2,#2+1) -- (#1+2,#2+2);
			\filldraw[#3, #4] (#1+3,#2+0) -- (#1+3,#2+1) -- (#1+3,#2+2) ;
		}
		
		
		\filldraw[red!70] (1,0) -- (1,1) -- (1,2) --  (3,2) --  (3,1) --  (2,1) -- (2,0);
		\drawmesh{0}{0}{}{}
		\draw[blue] (1,-0.3) .. controls (1.2, 0.3) .. (1.55,1);
		\draw[blue] (1.55,1) .. controls (1.7, 1.2) .. (2, 1.3);
		\draw[blue] (2,1.3) .. controls (2.52, 1.9) .. (2.9, 2.3);
		\node[scale=1, blue] at (2.6, 2.2) {$S_{\lin}$};
		\draw[blue!60!green] (1.1,-0.3) .. controls (1.2, 0.3) .. (1.92,1);
		\draw[blue!60!green] (1.7,1) .. controls (1.9, 1.1) .. (2, 1.2);
		\draw[blue!60!green] (2,1.2) .. controls (2.52, 1.7) .. (3, 2.3);
		\node[scale=1, blue!60!green] at (3.1, 2.2) {$S_h$};
		\node[scale=1] at (-0.3, 0) {$t_{n-1}$};
		\node[scale=1] at (-0.3, 1) {$t_{n}$};
		\node[scale=1] at (-0.3, 2) {$t_{n+1}$};
		\draw [decorate,decoration={brace,amplitude=3pt, mirror},xshift=0,yshift=-2pt]
		(1,0) -- (2,0) node [red,midway,yshift=-0.32cm] 
		{ $\Omega^\Gamma_{n}$};

		\draw [decorate,decoration={brace,amplitude=3pt},xshift=0,yshift=2pt]
		(1,2) -- (3,2) node [red,midway,yshift=0.32cm] 
		{ $\Omega^\Gamma_{n+1}$};

		\node[scale=1] at (1.28, 0.85) {$\Theta_{h}^{n}$};
		\draw[->, thick] (1.09,0.7) -- (1.09,0.99);
		\node[scale=1] at (1.37, 1.15) {$\Theta_{h}^{n+1}$};
		\draw[->, thick] (1.09,1.3) -- (1.09,1.01);
		\node[scale=1] at (3.35, 0.5) {$Q_{n}^{S}$};
		\node[scale=1] at (3.4, 1.5) {$Q_{n+1}^{S}$};
		\node[scale=1] at (1.6, 0.1) {$K_1$};
		\node[scale=1] at (2.6, 0.1) {$K_2$};
		\node[scale=1] at (2.1, 0.9) {$F$};

	\end{tikzpicture}
	\caption{Sketches of the discrete manifolds $S_{\lin}$ and $S_h$. Other than $S_{\lin}$, the surface $S_h$ is not necessarily continuous. In the calculation of ${\Theta_{h}^{n}(K_1\times I_n)}$ the degrees of freedom on the face $F$ only depend on $K_1$. In the calculation of ${\Theta_{h}^{n+1}(K_1\times I_{n+1})}$ the same degrees of freedoms depend on ${K_2\in\T^\Gamma_{n+1}}$ as well.}
	\label{deform_disc}
\end{figure}

\begin{remark} \label{remconsist}\rm 
	In general we have $\Omega_n^{\Gamma} \neq \Omega_{n+1}^\Gamma$, and therefore  $\Theta_{h}^{n} \neq \Theta_{h}^{n+1}$.  
	Take a fixed $t=t_n$ and let $K \in \T$ be a tetrahedron that is cut by $\Gamma_{\rm lin}(t_n)$. By construction we then have $K \in \T_n^\Gamma \cap \T_{n+1}^\Gamma$. The sets of neighboring tetrahedra of $K$ in $ \T_n^\Gamma$ and in $\T^\Gamma_{n+1}$, however, are not necessarily the same. Due to an averaging at the finite element nodes used in the construction of $\Theta_{h,t}^n$, cf.  \cite{lehrenfeld2016high,grande2018analysis},  the mappings $\Theta_{h,t_n}^{n}$ and $\Theta_{h,t_n}^{n+1}$ are in general not the same on $K$. Thus in general we have
	that for certain $\bx\in \Omega_n^{\Gamma}\cap \Omega^\Gamma_{n+1}$
	\begin{equation*}
		\Theta_h^{n+1}(\bx,t_n)=\Theta_{h,t_n}^{n+1}(\bx)\neq\Theta_{h,t_n}^{n}(\bx)=\Theta_h^{n}(\bx,t_n), 
	\end{equation*}
	which in particular implies $\Gamma_h^{n+1}(t_n) \neq \Gamma^n(t_n)$, cf.  \Cref{deform_disc} for an illustration.
	Due to the time stepping procedure in the finite element method, we weakly pass the discrete solution at the end point of a time interval $I_n$, i.e. on $\Gamma^n(t_n)$, to the next time interval, i.e. to $\Gamma_h^{n+1}(t_n)$. As these surfaces are not equal we need a suitable projection when defining the corresponding discrete bilinear form below.
\end{remark}

Below we will define and compute surface differential operators using discrete normal vectors on $S_{\lin}$ and $S_h$, both in the spatial and in the space-time sense. For a given $n$ and almost every $(\bx,t) \in Q_{n}^S$ we define (with $D\Theta$ the spatial Jacobian of $\Theta$)
	\begin{equation}
		\bn_{\lin}(\bx,t) \coloneqq  \frac{\nabla\hat{\phi}_h(\bx,t)}{\norm{\nabla\hat{\phi}_h(\bx,t)}}, \quad \bn_h(\Theta_h^n(\bx,t))\coloneqq \frac{D^{-\intercal}\Theta_{h,\spa}^n(\bx,t)\bn_{\lin}(\bx,t)}{\norm{D^{-\intercal}\Theta_{h,\spa}^n(\bx,t)\bn_{\lin}(\bx,t)}},  
		\label{nh_defs}
	\end{equation}
	with $\hat{\phi}_h$ as defined in \cref{phi_h_dach_definition} and $\Theta_{h,\spa}^n$ as in \Cref{Theta_Def}.
	Restricted to $(\bx,t)$ with $\bx \in \Gamma_{\rm lin}(t)$ these vectors are indeed the unit (spatial) normal vectors to these surface approximations.	We also introduce space-time variants, for almost every $(\bx,t) \in Q_{n}^S$:
	\begin{align}
		\bn_{S_{\lin}}(\bx,t) &\coloneqq  \frac{\nabla_{(\bx,t)} \hat{\phi}_h(\bx,t)}{\norm{\nabla_{(\bx,t)} \hat{\phi}_h(\bx,t)}} \label{nslin_def} \\
		\bn_{S_{h}}(\Theta_h^n(\bx,t)) & \coloneqq  \frac{D_{(\bx,t)}^{-\intercal}\Theta_h^n(\bx,t)\bn_{S_{\lin}}(\bx,t)}{\norm{D^{-\intercal}_{(\bx,t)}\Theta_h^n(\bx,t)\bn_{S_{\lin}}(\bx,t)}} .\label{defintion_nsh_higherorder}
	\end{align}
For $(\bx,t)$ restricted to $S_{\rm lin}^n$  these normals are the unit space-time normals to the respective space-time surfaces, cf.  \cite[Subsection 4.2.4]{Diss_Sass} for more discussion.
\subsection{Discrete surface differential operators}
We consider the prisms that are cut by $S_{\lin}$, i.e. prisms in $Q_{h}^S$, and apply the mesh deformation to obtain the sets:
\begin{equation*}
	Q_{\Theta,n}^{S}\coloneqq \Theta_h^n(Q^S_{n}),~n=1,\dots,N, \quad Q_{\Theta}^{S}\coloneqq\bigcup_{n=1}^N Q_{\Theta,n}^{S}.
\end{equation*}
 Almost everywhere on the domain $Q_{\Theta}^{S}$ we define  the projections $\bP_{h}\coloneqq\bI-\bn_{h}\bn_{h}^\intercal$ and $\bP_{S_h}\coloneqq\bI-\bn_{S_h}\bn_{S_h}^\intercal$. For a fixed $t$ and smooth extension $u^e$  ($\bu^e$) of $u$ ($\bu$) defined on $\Gamma_h(t)$ we define
 \begin{equation*}
{\nabla_{\Gamma_h}u\coloneqq\bP_h\nabla u^e}, \quad	\Div_{\Gamma_h}\mathbf{u}\coloneqq \tr\left(\bP_h \nabla\mathbf{u}^e\right). 
\end{equation*}
Similarly, for sufficiently smooth $u$ ($\bu$) defined on $K_S\in \T_{S_h}$ we define
\begin{equation*}
{\nabla_{S_h}u\coloneqq\bP_{S_h}\nabla_{(\bx,t)}u^e}, \quad	\Div_{S_h}\mathbf{u}\coloneqq \tr\left(\bP_{S_h} \nabla_{(\bx,t)}\mathbf{u}^e\right). 
\end{equation*}
It is convenient to introduce broken spaces of ``sufficiently smooth'' functions on $S_h$ (similar to \cref{brokenspace}) and $Q_{\Theta}^{S}$:
\begin{alignat*}{3}
	W^{h,b} & \coloneqq \bigg\lbrace v\in L^2(S_h): &&v\restrict{S_h^n}\in H^1(S_h^n), &&n=1,\dots,N\bigg\rbrace, \\
	V_{\reg,h} & \coloneqq \bigg\lbrace v\in L^2(Q_{\Theta}^{S}): \,&&v\restrict{Q_{\Theta,n}^{S}}\in H^1(Q_{\Theta,n}^{S}),\,&&n=1,\dots,N, v\restrict{S_h}\in W^{h,b}\bigg\rbrace.
\end{alignat*}
For more information concerning Sobolev spaces on Lipschitz manifolds we refer to the literature, e.g.  \cite{elliott2017unified}
\subsection{Integration over Lipschitz space-time surfaces} \label{subsection_different_st_integrals}
The approximate surface $S_h^n$ defined in \eqref{sh_def} in general has only Lipschitz smoothness. This loss of smoothness, compared to the smooth space-time surface $S$, has consequences for (partial) integration formulas, which play a key role in the derivation and analysis of the discretization method. In this section, we introduce discrete variants of the transformation formula \cref{transformationsformelMAIN} and the partial integration identity \cref{pi_cont}. These will be used in the derivation of the discretization method in \Cref{secmethod} and in the discussion of implementation aspects in \Cref{secimplementation}. Almost everywhere on $S_h^n$, $n\in\{1,\dots,N\}$, we define 
\begin{equation}
	V_h\coloneqq \frac{-\partial \left(\hat{\phi}_h\circ(\Theta_h^n)^{-1}\right)}{\partial t}\norm{\nabla\left(\hat{\phi}_h\circ(\Theta_h^n)^{-1}\right)}^{-1}.\label{Vh_def_ho}
\end{equation}
\begin{remark} \rm
Note that the discrete space-time surface $S_h^n$ is defined via the surface $S_{\rm lin}^n$, which is the zero level of $\hat{\phi}_h$, and the parametric mapping $\Theta_h^n$, cf. \eqref{sh_def}. The function $V_h$ is a discrete analogue to the normal velocity $\bw\cdot\bn$ on $S$ in the following sense. In the continuous setting we use the level set equation $\frac{\partial \phi}{\partial t} + \bw \cdot \nabla \phi =0$  to see that on $S$  we have the relations
\begin{equation} \label{wn_lset}
	\bw\cdot\bn =\bw\cdot \frac{\nabla \phi}{\norm{\nabla \phi}}=\frac{-1}{\norm{\nabla \phi}}\frac{\partial \phi}{\partial t}.
\end{equation}
The function $V_h$ is a discrete approximation of the function on the right-hand side in \eqref{wn_lset}.
\end{remark}

We introduce a discrete version of the integral transformation \cref{transformationsformelMAIN}.
\begin{theorem}\label{theo_discrete_trafo}
For $g_h\in L^2(S_h^n)$ we have
	\begin{equation}
		\int_{t_{n-1}}^{t_n}\int_{\Gamma^n_h(t)}g_h\dif s_h\dif t=\int_{S_h^n}\frac{g_h}{\sqrt{1+V_h^2}}\dif \sigma_h, \quad n\in \{1,\dots, N\},\label{ho_integral_trafo}
	\end{equation}
	where $V_h$, $S_h^n$ and $\Gamma_h^n(t)$ are defined in \cref{Vh_def_ho}, \cref{sh_def} and \cref{gamma_h_def}, respectively.
\end{theorem}

\begin{figure}[!htbp]
	\centering
	\begin{tikzpicture}[scale = 0.64]
		\newcommand{\z}{0.05}
		\newcommand{\Hoehe}{8.0}
		\newcommand{\Tiefe}{2.0}
		\newcommand{\Shift}{3.0}
		
		\newcommand{\lambdaA}{0.5}
\newcommand{\PunktA}{(\lambdaA*0+\Shift-\lambdaA*\Shift+0.6,-0.2+\lambdaA*0+\Tiefe-\lambdaA*\Tiefe)}
\newcommand{\lambdaB}{0.7}
\newcommand{\PunktB}{(\lambdaB*5,0-0.1)}
		\newcommand{\lambdaC}{0.15}
	\newcommand{\PunktC}{(\Shift,\Tiefe+\lambdaC*\Hoehe)}
	\newcommand{\lambdaD}{0.45}
	\newcommand{\PunktD}{(5,\lambdaD*\Hoehe)}
	\newcommand{\lambdaE}{0.6}
	\newcommand{\PunktE}{(5+\Shift,\Tiefe+\lambdaE*\Hoehe)}
	\newcommand{\lambdaF}{0.8}
	\newcommand{\PunktF}{(10,\lambdaF*\Hoehe)}
	\newcommand{\lambdaG}{0.5}
	\newcommand{\PunktG}{(\lambdaG*5+\lambdaG*\Shift+10+\Shift-10*\lambdaG-\lambdaG*\Shift,\Tiefe+\Hoehe-0.35+0.5)}
	\newcommand{\lambdaH}{0.5}
	\newcommand{\PunktH}{(\lambdaH*10+10+\Shift-\lambdaH*10-\lambdaH*\Shift+0.6,\Tiefe-\lambdaH*\Tiefe+\Hoehe+0.5)}
	\newcommand{\AtildeShifth}{-0.5}\newcommand{\AtildeShiftv}{1.0}
	\newcommand{\PunktAtilde}{(\lambdaA*0+\Shift-\lambdaA*\Shift+\AtildeShifth-0.5,\lambdaA*0+\Tiefe-\lambdaA*\Tiefe+\AtildeShiftv)}
	\newcommand{\BtildeShifth}{1.0}\newcommand{\BtildeShiftv}{-0.4}
	\newcommand{\PunktBtilde}{(\lambdaB*5+\BtildeShifth,-0.2+\BtildeShiftv)}
	\newcommand{\GtildeShifth}{-1.0}\newcommand{\GtildeShiftv}{+1.0}
	\newcommand{\PunktGtilde}{(\lambdaG*5+\lambdaG*\Shift+10+\Shift-10*\lambdaG-\lambdaG*\Shift+\GtildeShifth,\Tiefe+\Hoehe+\GtildeShiftv)}
	\newcommand{\HtildeShifth}{1.0}\newcommand{\HtildeShiftv}{-0.6}
	\newcommand{\PunktHtilde}{(\lambdaH*10+10+\Shift-\lambdaH*10-\lambdaH*\Shift+\HtildeShifth,\Tiefe-\lambdaH*\Tiefe+\Hoehe+\HtildeShiftv)}
		
		\newcommand{\ABh}{0.0}\newcommand{\ABv}{-0.2}
		\newcommand{\BAh}{0.0}\newcommand{\BAv}{-0.4}
		\newcommand{\BDh}{0.0}\newcommand{\BDv}{0.2}
		\newcommand{\DBh}{0.0}\newcommand{\DBv}{0.5}
		\newcommand{\ACh}{0.0}\newcommand{\ACv}{0.2}
		\newcommand{\CAh}{0.0}\newcommand{\CAv}{-0.4}
		\newcommand{\CEh}{0.0}\newcommand{\CEv}{0.2}
		\newcommand{\ECh}{0.0}\newcommand{\ECv}{0.2}
		\newcommand{\DFh}{0.0}\newcommand{\DFv}{-0.5}
		\newcommand{\FDh}{0.0}\newcommand{\FDv}{0.0}
		\newcommand{\CDh}{0.0}\newcommand{\CDv}{0.1}
		\newcommand{\DCh}{0.0}\newcommand{\DCv}{0.1}
		\newcommand{\DEh}{0.0}\newcommand{\DEv}{0.2}
		\newcommand{\EDh}{0.0}\newcommand{\EDv}{0.2}
		\newcommand{\EFh}{0.0}\newcommand{\EFv}{-0.3}
		\newcommand{\FEh}{0.0}\newcommand{\FEv}{0.0}
		\newcommand{\EGh}{0.0}\newcommand{\EGv}{-0.3}
		\newcommand{\GEh}{0.0}\newcommand{\GEv}{-0.8}
		\newcommand{\FHh}{0.0}\newcommand{\FHv}{-0.5}
		\newcommand{\HFh}{0.0}\newcommand{\HFv}{-0.3}
		\newcommand{\GHh}{-0.1}\newcommand{\GHv}{-0.05}
		\newcommand{\HGh}{-0.1}\newcommand{\HGv}{-0.05}
		\newcommand{\AAtildeh}{0.0}\newcommand{\AAtildev}{-0.3}
		\newcommand{\AtildeAh}{0.0}\newcommand{\AtildeAv}{-0.2}
		\newcommand{\BBtildeh}{0.0}\newcommand{\BBtildev}{0.0}
		\newcommand{\BtildeBh}{0.0}\newcommand{\BtildeBv}{0.0}
		\newcommand{\HHtildeh}{0.01}\newcommand{\HHtildev}{0.03}
		\newcommand{\HtildeHh}{0.1}\newcommand{\HtildeHv}{0.1}
		\newcommand{\GGtildeh}{0.1}\newcommand{\GGtildev}{0.1}
		\newcommand{\GtildeGh}{0.1}\newcommand{\GtildeGv}{0.1}
		
		\newcommand{\PunktAB}{(0.333*\lambdaA*0+0.333*\Shift-0.333*\lambdaA*\Shift+0.666*\lambdaB*5+\ABh,0.333*\lambdaA*0+0.333*\Tiefe-0.333*\lambdaA*\Tiefe+0.666*0+\ABv)}
		\newcommand{\PunktBA}{(0.666*\lambdaA*0+0.666*\Shift-0.666*\lambdaA*\Shift+0.333*\lambdaB*5+\BAh,0.666*\lambdaA*0+0.666*\Tiefe-0.666*\lambdaA*\Tiefe+0.333*0+\BAv)}
		\newcommand{\PunktBD}{(0.333*\lambdaB*5+0.666*5+\BDh,0.333*0+0.666*\lambdaD*\Hoehe+\BDv)}
		\newcommand{\PunktDB}{(0.666*\lambdaB*5+0.333*5+\DBh,0.666*0+0.333*\lambdaD*\Hoehe+\DBv)}
		\newcommand{\PunktAC}{(0.333*\lambdaA*0+0.333*\Shift-0.333*\lambdaA*\Shift+0.666*\Shift+\ACh,0.333*\lambdaA*0+0.333*\Tiefe-0.333*\lambdaA*\Tiefe+0.666*\Tiefe+0.666*\lambdaC*\Hoehe+\ACv)}
		\newcommand{\PunktCA}{(0.666*\lambdaA*0+0.666*\Shift-0.666*\lambdaA*\Shift+0.333*\Shift+\CAh,0.666*\lambdaA*0+0.666*\Tiefe-0.666*\lambdaA*\Tiefe+0.333*\Tiefe+0.333*\lambdaC*\Hoehe+\CAv)}
		\newcommand{\PunktCE}{(0.333*\Shift+0.666*5+0.666*\Shift+\CEh,0.333*\Tiefe+0.333*\lambdaC*\Hoehe+0.666*\Tiefe+0.666*\lambdaE*\Hoehe+\CEv)}
		\newcommand{\PunktEC}{(0.666*\Shift+0.333*5+0.333*\Shift+\ECh,0.666*\Tiefe+0.666*\lambdaC*\Hoehe+0.333*\Tiefe+0.333*\lambdaE*\Hoehe+\ECv)}
		\newcommand{\PunktDF}{(0.333*5+0.666*10+\DFh,0.333*\lambdaD*\Hoehe+0.666*\lambdaF*\Hoehe+\DFv)}
		\newcommand{\PunktFD}{(0.666*5+0.333*10+\FDh,0.666*\lambdaD*\Hoehe+0.333*\lambdaF*\Hoehe+\FDv)}
		\newcommand{\PunktCD}{(0.333*\Shift+0.666*5+\CDh,0.333*\Tiefe+0.333*\lambdaC*\Hoehe+0.666*\lambdaD*\Hoehe+\CDv)}
		\newcommand{\PunktDC}{(0.666*\Shift+0.333*5+\DCh,0.666*\Tiefe+0.666*\lambdaC*\Hoehe+0.333*\lambdaD*\Hoehe+\DCv)}
		\newcommand{\PunktDE}{(0.666*5+0.666*\Shift+0.333*5+\DEh,0.666*\Tiefe+0.666*\lambdaE*\Hoehe+0.333*\lambdaD*\Hoehe+\DEv)}
		\newcommand{\PunktED}{(0.333*5+0.333*\Shift+0.666*5+\EDh,0.333*\Tiefe+0.333*\lambdaE*\Hoehe+0.666*\lambdaD*\Hoehe+\EDv)}
		\newcommand{\PunktEF}{(0.333*5+0.333*\Shift+0.666*10+\EFh,0.333*\Tiefe+0.333*\lambdaE*\Hoehe+0.666*\lambdaF*\Hoehe+\EFv)}
		\newcommand{\PunktFE}{(0.666*5+0.666*\Shift+0.333*10+\FEh,0.666*\Tiefe+0.666*\lambdaE*\Hoehe+0.333*\lambdaF*\Hoehe+\FEv)}
		\newcommand{\PunktEG}{(0.333*5+0.333*\Shift+0.666*\lambdaG*5+0.666*\lambdaG*\Shift+0.666*10+0.666*\Shift-0.666*10*\lambdaG-0.666*\lambdaG*\Shift+\EGh,0.333*\Tiefe+0.333*\lambdaE*\Hoehe+0.666*\Tiefe+0.666*\Hoehe+\EGv)}
		\newcommand{\PunktGE}{(0.666*5+0.666*\Shift+0.333*\lambdaG*5+0.333*\lambdaG*\Shift+0.333*10+0.333*\Shift-0.333*10*\lambdaG-0.333*\lambdaG*\Shift+\GEh,0.666*\Tiefe+0.666*\lambdaE*\Hoehe+0.333*\Tiefe+0.333*\Hoehe+\GEv)}
		\newcommand{\PunktFH}{(0.333*10+0.666*\lambdaH*10+0.666*10+0.666*\Shift-0.666*\lambdaH*10-0.666*\lambdaH*\Shift+\FHh,0.333*\lambdaF*\Hoehe+0.666*\Tiefe-0.666*\lambdaH*\Tiefe+0.666*\Hoehe+\FHv)}
		\newcommand{\PunktHF}{(0.666*10+0.333*\lambdaH*10+0.333*10+0.333*\Shift-0.333*\lambdaH*10-0.333*\lambdaH*\Shift+\HFh,0.666*\lambdaF*\Hoehe+0.333*\Tiefe-0.333*\lambdaH*\Tiefe+0.333*\Hoehe+\HFv)}
		\newcommand{\PunktGH}{(0.333*\lambdaG*5+0.333*\lambdaG*\Shift+0.333*10+0.333*\Shift-0.333*10*\lambdaG-0.333*\lambdaG*\Shift+0.666*\lambdaH*10+0.666*10+0.666*\Shift-0.666*\lambdaH*10-0.666*\lambdaH*\Shift+\GHh,0.333*\Tiefe+0.333*\Hoehe+0.666*\Tiefe-0.666*\lambdaH*\Tiefe+0.666*\Hoehe+\GHv)}
		\newcommand{\PunktHG}{(0.666*\lambdaG*5+0.666*\lambdaG*\Shift+0.666*10+0.666*\Shift-0.666*10*\lambdaG-0.666*\lambdaG*\Shift+0.333*\lambdaH*10+0.333*10+0.333*\Shift-0.333*\lambdaH*10-0.333*\lambdaH*\Shift+\HGh,0.666*\Tiefe+0.666*\Hoehe+0.333*\Tiefe-0.333*\lambdaH*\Tiefe+0.333*\Hoehe+\HGv)}
		\newcommand{\PunktAAtilde}{(\lambdaA*0+\Shift-\lambdaA*\Shift+0.666*\AtildeShifth+\AAtildeh,\lambdaA*0+\Tiefe-\lambdaA*\Tiefe+0.666*\AtildeShiftv+\AAtildev)}
		\newcommand{\PunktAtildeA}{(\lambdaA*0+\Shift-\lambdaA*\Shift+0.333*\AtildeShifth+\AtildeAh,\lambdaA*0+\Tiefe-\lambdaA*\Tiefe+0.333*\AtildeShiftv+\AtildeAv)}
		\newcommand{\PunktBBtilde}{(\lambdaB*5+0.333*\BtildeShifth+\BBtildeh,0.333*\BtildeShiftv+\BBtildev)}
		\newcommand{\PunktBtildeB}{(\lambdaB*5+0.666*\BtildeShifth+\BtildeBh,0.666*\BtildeShiftv+\BtildeBv)}
		\newcommand{\PunktHHtilde}{(\lambdaH*10+10+\Shift-\lambdaH*10-\lambdaH*\Shift+0.666*\HtildeShifth+\HHtildeh,\Tiefe-\lambdaH*\Tiefe+\Hoehe+0.666*\HtildeShiftv+\HHtildev)}
		\newcommand{\PunktHtildeH}{(\lambdaH*10+10+\Shift-\lambdaH*10-\lambdaH*\Shift+0.333*\HtildeShifth+\HtildeHh,\Tiefe-\lambdaH*\Tiefe+\Hoehe+0.333*\HtildeShiftv+\HtildeHv)}
		\newcommand{\PunktGGtilde}{(\lambdaG*5+\lambdaG*\Shift+10+\Shift-10*\lambdaG-\lambdaG*\Shift+0.666*\GtildeShifth+\GGtildeh,\Tiefe+\Hoehe+0.666*\GtildeShiftv+\GGtildev)}
		\newcommand{\PunktGtildeG}{(\lambdaG*5+\lambdaG*\Shift+10+\Shift-10*\lambdaG-\lambdaG*\Shift+0.333*\GtildeShifth+\GtildeGh,\Tiefe+\Hoehe+0.333*\GtildeShiftv+\GtildeGv)}
		
		\node[left] at (-1.0,0.0) {$t_{n-1}$};
		\node[left] at (-1.0,\Hoehe) {$t_n$};
		
		\newcommand{\FdM}{blue}
		
		
		\node[below, color=\FdM] at \PunktBtilde {$\Gamma_h^n(t_{n-1})$};
		\node[left, color=\FdM] at \PunktGtilde {$\Gamma_h^n(t_n)$};
		
			\newcommand{\FF}{30}
			\fill[\FdM!\FF] \PunktA .. controls \PunktBA and \PunktAB .. \PunktB -- \PunktB .. controls \PunktDB and \PunktBD .. \PunktD -- %
			\PunktD .. controls \PunktCD and \PunktDC .. \PunktC -- \PunktC .. controls \PunktAC and \PunktCA .. \PunktA;
			\fill[\FdM!\FF] \PunktC .. controls \PunktDC and \PunktCD .. \PunktD -- \PunktD .. controls \PunktFD and \PunktDF .. \PunktF -- %
			\PunktF .. controls \PunktEF and \PunktFE .. \PunktE -- \PunktE .. controls \PunktCE and \PunktEC .. \PunktC;
			\fill[\FdM!\FF] \PunktE .. controls \PunktFE and \PunktEF .. \PunktF -- \PunktF .. controls \PunktHF and \PunktFH .. \PunktH -- %
			\PunktH .. controls \PunktGH and \PunktHG .. \PunktG -- \PunktG .. controls \PunktEG and \PunktGE .. \PunktE;

		\newcommand{\Radius}{2pt}
		\filldraw[fill=black, ultra thick] \PunktA circle (\Radius);
		\filldraw[fill=black, ultra thick] \PunktB circle (\Radius);
		\filldraw[fill=black, ultra thick] \PunktC circle (\Radius);
		\filldraw[fill=black, ultra thick] \PunktD circle (\Radius);
		\filldraw[fill=black, ultra thick] \PunktE circle (\Radius);
		\filldraw[fill=black, ultra thick] \PunktF circle (\Radius);
		\filldraw[fill=black, ultra thick] \PunktG circle (\Radius);
		\filldraw[fill=black, ultra thick] \PunktH circle (\Radius);
		
		\draw[ultra thick, color=\FdM] \PunktA .. controls \PunktBA and \PunktAB .. \PunktB;
		\draw[thick, color=\FdM]       \PunktA .. controls \PunktCA and \PunktAC .. \PunktC;
		\draw[thick, color=\FdM]       \PunktC .. controls \PunktEC and \PunktCE .. \PunktE;
		\draw[thick, color=\FdM]       \PunktE .. controls \PunktGE and \PunktEG .. \PunktG;
		\draw[thick, color=\FdM]       \PunktB .. controls \PunktDB and \PunktBD .. \PunktD;
		\draw[thick, color=\FdM]       \PunktD .. controls \PunktFD and \PunktDF .. \PunktF;
		\draw[thick, color=\FdM]       \PunktF .. controls \PunktHF and \PunktFH .. \PunktH;
		\draw[ultra thick, color=\FdM] \PunktH .. controls \PunktGH and \PunktHG .. \PunktG;
		\draw[ultra thick, color=\FdM] \PunktA .. controls \PunktAtildeA and \PunktAAtilde .. \PunktAtilde;
		\draw[ultra thick, color=\FdM] \PunktB .. controls \PunktBtildeB and \PunktBBtilde .. \PunktBtilde;
		\draw[ultra thick, color=\FdM] \PunktG .. controls \PunktGtildeG and \PunktGGtilde .. \PunktGtilde;
		\draw[ultra thick, color=\FdM] \PunktH .. controls \PunktHtildeH and \PunktHHtilde .. \PunktHtilde;
		
		\draw[thick, dashed, color=\FdM]       \PunktC .. controls \PunktDC and \PunktCD .. \PunktD;
		\draw[thick, dashed, color=\FdM]       \PunktD .. controls \PunktED and \PunktDE .. \PunktE;
		\draw[thick, dashed, color=\FdM]       \PunktE .. controls \PunktFE and \PunktEF .. \PunktF;
		
	\draw[ultra thick] (0.0,0.0) .. controls(2.5, -0.2) .. (5.0,0.0);
	\draw[ultra thick] (5.0,0.0) .. controls (8, -0.4) .. (10.0,0.0);
	\draw[dashed] (\Shift,\Tiefe) .. controls (3+\Shift, \Tiefe-0.4) .. (5.0+\Shift,\Tiefe);
	\draw[dashed] (5+\Shift,\Tiefe) .. controls (7+\Shift, \Tiefe+0.3) .. (10.0+\Shift,\Tiefe);
	\draw[dashed] (0.0,0.0) .. controls(\Shift/1.5, \Tiefe/4) .. (\Shift,\Tiefe);
	\draw[dashed] (\Shift,\Tiefe) .. controls(\Shift+0.2, \Tiefe/3) .. (5.0,0.0);
	\draw[dashed] (5.0,0.0) .. controls(6, \Tiefe/1.5) .. (5.0+\Shift,\Tiefe);
	\draw[dashed] (5.0+\Shift,\Tiefe) .. controls (7.7+\Shift/2, \Tiefe/1.7) .. (10.0,0.0);
	\draw[ultra thick] (10.0,0.0) .. controls (10 + \Shift/3, \Tiefe/1.5) .. (10.0+\Shift,\Tiefe);
	
	\draw[ultra thick] (0.0,0.0+\Hoehe) .. controls (3, 0.3+\Hoehe) .. (5.0,0.0+\Hoehe);
	\draw[ultra thick] (5.0,0.0+\Hoehe) .. controls (8, -0.3+\Hoehe) .. (10.0,0.0+\Hoehe);
	\draw[ultra thick] (\Shift,\Tiefe+\Hoehe) .. controls (3+\Shift, \Tiefe+\Hoehe+0.1) .. (5.0+\Shift,\Tiefe+\Hoehe);
	\draw[ultra thick] (\Shift +5 ,\Tiefe+\Hoehe) .. controls (7+\Shift, \Tiefe+\Hoehe+0.2) .. (10.0+\Shift,\Tiefe+\Hoehe);
	\draw[ultra thick] (0.0,0.0+\Hoehe) .. controls (\Shift/1.5, \Tiefe/2+\Hoehe) .. (\Shift,\Tiefe+\Hoehe);
	\draw[ultra thick] (\Shift,\Tiefe+\Hoehe) .. controls (\Shift/2+1.5, \Hoehe+\Tiefe/1.7) .. (5.0,0.0+\Hoehe);
	\draw[ultra thick] (5.0,0.0+\Hoehe) .. controls (4+ \Shift/2, \Tiefe/2+\Hoehe) .. (5.0+\Shift,\Tiefe+\Hoehe);
	\draw[ultra thick] (5.0+\Shift,\Tiefe+\Hoehe) .. controls(7.5+\Shift/3, \Tiefe/4 + \Hoehe) .. (10.0,0.0+\Hoehe);
	\draw[ultra thick] (10.0,0.0+\Hoehe) .. controls(10+\Shift/2.4, \Tiefe/2+\Hoehe) .. (10.0+\Shift,\Tiefe+\Hoehe);
		
		\draw[ultra thick] (0.0,0.0) -- (0.0,0.0+\Hoehe);
		\draw[ultra thick] (5.0,0.0) -- (5.0,0.0+\Hoehe);
		\draw[ultra thick] (10.0,0.0) -- (10.0,0.0+\Hoehe);
		\draw[dashed]      (\Shift,\Tiefe) -- (\Shift,\Tiefe+\Hoehe);
		\draw[dashed]      (5.0+\Shift,\Tiefe) -- (5.0+\Shift,\Tiefe+\Hoehe);
		\draw[ultra thick] (10.0+\Shift,\Tiefe) -- (10.0+\Shift,\Tiefe+\Hoehe);
		\node[above, color=\FdM] at (3*\Shift,\Hoehe*0.57) {$S_{h}^n$};
		\node[above, color=\FdM] at (0.925*\Shift,\Hoehe*0.04) {$K_S^1$};
		\node[above, color=\FdM] at (1.81*\Shift,\Hoehe*0.7) {$K_S^2$};
		\draw[color=blue, dotted ] (1.81*\Shift,\Hoehe*0.7) -- (1.81*\Shift,\Hoehe*0.575);
		\node[above, color=\FdM] at (4.3,2.75) {$F$};
		\draw[->, color=red, ultra thick] (2.6, 0.34) -- ++(250:1.1);
		\node[above, color=red] at (1.5,-1) {$-\bnu_{\partial}$};
		\draw[->, color=red, ultra thick] (3.83, 3.5) -- ++(70:1.1);
		\node[above, color=red] at (3.55,4.3) {$\bnu_h$};
		\draw[->, color=red!70, ultra thick] (3.83, 3.5) -- ++(220:1.1);
		\node[above, color=red!70] at (3.77,2.15) {$\bnu_h$};
		\draw[->, color=red, ultra thick] (11, 9.4) -- ++(70:1.1);
		\node[above, color=red] at (11,10.35) {$\bnu_{\partial}$};
		\node[above, color=black] at (-0.4*\Shift,\Hoehe*0.5) {$Q_{\Theta,n}^{S}$};
	\end{tikzpicture}
	\caption{Illustration of the conormals. At the time slab boundaries, the conormals coincide with $\pm \bnu_{\partial}$. Due to the non-smoothness of $S_h$, $(\bnu_{h}\restrict{K_S^1})\restrict{F}\neq -(\bnu_{h}\restrict{K_S^2})\restrict{F}$ at a common face $F$ is possible.}
	\label{conormal}
\end{figure}
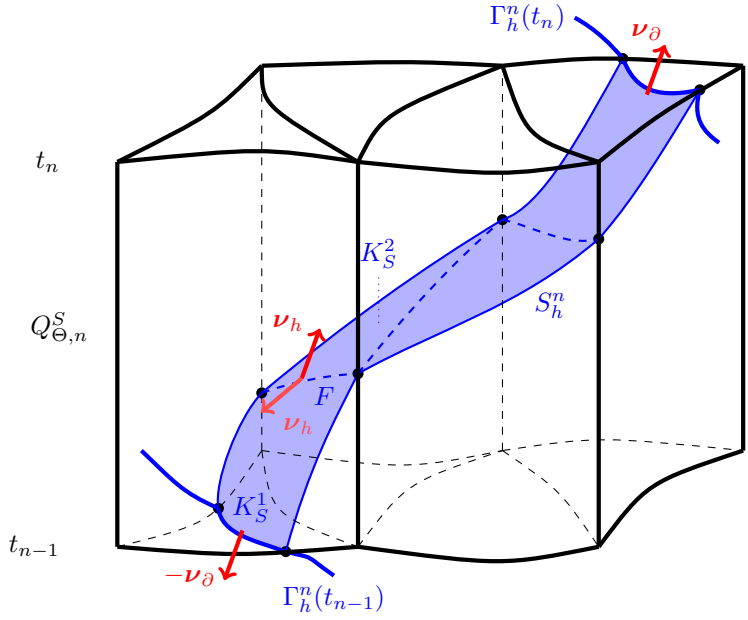

A proof is given in \cite[Theorem 4.27]{Diss_Sass}. Below, in \Cref{PI_Theorem}, we derive  an integration by parts identity on the discrete space-time manifold $S_h^n$. For this we need further definitions. First we define a discrete version of the material derivative  in \Cref{def_dot_u}. Let $\bw_S^\intercal\coloneqq(\bw^\intercal,1)$. For $v\in H^1(S_h^n)$, $n\in \{1,\dots,N\}$, we define the discrete material derivative as
\begin{equation}
	\mathring{v}\coloneqq \bw_S\cdot \nabla_{S_h}v=(\mathbf{P}_{S_h}\bw_S)\cdot \nabla_{S_h}v.\label{weakmatderivative}
\end{equation}
This derivative is tangential to the corresponding discrete space-time surface.  For $n= 1,\dots,N,$ let $\mathcal{F}_I^n\coloneqq \{ \partial K_S^1\cap\partial K_S^2:K_S^1,K_S^2\in \mathcal{T}_{S_h^n},K_S^1\neq K_S^2\}$ be the set of interior boundary faces of the elements in $\mathcal{T}_{S_h^n}$ and	${\mathcal{F}_T^n\coloneqq \{\partial K_S\cap \Gamma_h^n(t_n):K_S\in\mathcal{T}_{S_h^n}\}}$  the set of element boundaries of $\mathcal{T}_{S_h^n}$ that lie on the top time slab boundary. Similarly,  the bottom time slab boundary ${\mathcal{F}_B^n\coloneqq \{\partial K_S\cap \Gamma_h^n(t_{n-1}):K_S\in\mathcal{T}_{S_h^n}\}}$. Hence, ${\partial S_h^n= \mathcal{F}_T^n \cup \mathcal{F}_B^n}$. Almost everywhere on $\partial S_h^n$ the vector  
\begin{equation}
	\bnu_{\partial}\coloneqq \frac{1}{\sqrt{1+V_h^2}}\begin{pmatrix}V_h\bn_h\\1\end{pmatrix}\label{ndelta_defs}
\end{equation}
is the unit normal of $\partial S_h^n$ in positive $t$-direction. Hence, $\bnu_{\partial}$ is the unit outer normal of the top boundary of $S_h^n$ and $-\bnu_{\partial }$ is the unit outer normal of the bottom boundary of $S_h^n$. Note that $\bnu_{\partial}$ jumps between time slabs.   For $K_S\in \T_{S_h}$ the conormals are denoted by ${\bnu_h}\restrict{K_S}$, i.e.,  the unit outer normal on the faces of $K_S$ that is tangential to $K_S$. For  $K_S\in\T_{S_h^n}$  with $K_S\cap  \Gamma_h^n(t_n)\neq \emptyset$ the unit normal $\bnu_{\partial}$ is the conormal vector of $K_S$ at the top boundary of $S_h^n$. Analogously, if $K_S\cap \Gamma_h^n(t_{n-1})\neq \emptyset$, the vector $-\bnu_{\partial}$ is the conormal vector of $K_S$ at the bottom boundary of $S_h^n$, cf. \Cref{conormal} for an illustration. Consider an interior boundary face $F \in \mathcal{F}_I^n$ with $F=K_S^1\cap K_S^2$, $K_S^1,K_S^2\in \mathcal{T}_{S_h^n}$. For a function $v$ defined on $K_S^1 \cup K_S^2$  we define the (vector valued) conormal jump on $F$:
\begin{equation} \label{def_conormal_jump}
	{\left[v\right]_{\bnu}}\restrict{F} := 
	\big(v\restrict{K_S^1}{\bnu_h}\restrict{K_{S}^1} + v\restrict{K_S^2}{\bnu_h}\restrict{K_{S}^2}\big)\restrict{F}.
\end{equation}
Note that in general for the Lipschitz surface $S_h^n$ we do not have $C^1$ smoothness across $F$ and thus ${\bnu_h}\restrict{K_{S}^1} \neq - {\bnu_h}\restrict{K_{S}^2}$ on $F$. If the function $v$ is continuous across $F$ we have ${[v]_{\bnu}}\restrict{F} = v\restrict{F} {[1]_{\bnu}}\restrict{F}$.  
For one-sided values we use the standard notation
\begin{align*}
			u_+^n \coloneqq u_+(\cdot ,t_n)=\lim_{\eta\searrow 0}u(\cdot, t_n+\eta), \quad
			u_-^n \coloneqq u_-(\cdot ,t_n)=\lim_{\eta\searrow 0}u(\cdot, t_n-\eta).
\end{align*}
Below we  use a generic approximation of (the extension of) $\alpha= \sqrt{1+(\bw\cdot\bn)}$ on the discrete surface space-time surface $S_h$, denoted by  $\alpha_h$.
 One specific possibility is $\alpha_h=\sqrt{1+V_h^2}$. However, as we will see below, other choices may be better.
\begin{theorem}\label{PI_Theorem}
	On $S_h$, let $R\coloneqq \frac{1}{\alpha_h}\bw_S\cdot \bnu_{\partial }$. For $u,v\in H^1(S_h^n)$, $n=1,\dots,N$, the following identity holds:
	\begin{align}
	\begin{aligned}
	&\int_{S_h^n}\frac{1}{\alpha_h}\mathring{u}v\dif \sigma_h=-\int_{S_h^n}\frac{1}{\alpha_h}u\mathring{v}\dif \sigma_h +\int_{\Gamma_h^n(t_n)}u_{-}^{n}v_{-}^{n}R_-^n\dif s_h\\
	&\qquad-\int_{\Gamma_h^n(t_{n-1})}u_{+}^{n-1}v_{+}^{n-1}R_+^{n-1}\dif s_h  + \sum_{F\in \mathcal{F}_I^n}\int_{F} u v \bw_S\cdot \left[\frac{1}{\alpha_h}\right]_{\bnu} \dif F\\
	&\qquad  -\sum_{K_S\in  \mathcal{T}_{S_h^n}}\int_{K_S} uv\Div_{S_h}\left( \frac{1}{\alpha_h} \bP_{S_h}\bw_S\right)\dif \sigma_h.\label{PI-unschoen}
	\end{aligned}
	\end{align}
\end{theorem}

A proof is given in \Cref{proof_PI}. The partial integration \cref{PI-unschoen} identity contains several geometric perturbation terms that do not arise for the smooth surface $S$. These terms play a role in the definition and error analysis of the space time finite element method, cf. the next section. The proof of \Cref{PI_Theorem} relies on the surface triangulation of $S_h^n$. 
Using \cref{ho_integral_trafo} the surface integrals $\int_{S_h^n} \cdot~\dif \sigma_h $  can be reformulated as product space-time integrals $\int_{t_{n-1}}^{t_n} \int_{\Gamma_h^n(t)} \cdot~ \dif s_h \dif t$.  However, a reformulation  of the integrals $\int_{F} \cdot~\dif F$ and $\int_{K_S} \cdot~\dif \sigma_h$ in \Cref{PI_Theorem} as product space-time integrals is not natural. 

\subsection{Space-time trace finite element discretization} \label{secmethod}
In this section we introduce a fully discrete higher order space-time discretization of \cref{prob_cont}. This  method essentially uses the same approach as in \cite{olshanskii2014eulerian,olshanskii2014error}, but there are important differences, as explained in the Introduction, cf. \Cref{secIntro}. 
The method is a space-time Eulerian method that uses standard tensor product finite element spaces, combined with a parametric mapping to obtain a feasible method with higher order accuracy.
More precisely, for $k_\spa,k_\ti\in \mathbb{N}$ we consider the same product space $V_h^{k_\spa,k_\ti}$, cf. \eqref{st_fe_space}, as used for the level set approximation and combine it with the space-time parametric mapping \eqref{Theta_Def}, which depends on the (higher order accurate) level set approximation $\phi_h$. We define the parametric finite element space
\begin{equation}
	V_{h,\Theta}^{k_\spa,k_\ti}\coloneqq \left\lbrace v_h:Q_{\Theta}^{S}\rightarrow \R:  (\bx,t)\mapsto v_h(\Theta_h^n(\bx,t))\in V_h^{k_\spa,k_\ti}\restrict{Q^S_{n}}, n=1,\dots,N \right\rbrace. \label{Vh_Theta_Def}
\end{equation}
Note that $V_{h,\Theta}^{k_\spa,k_\ti}\subset V_{\reg,h}$. For $\beta \in [0,1]$ and $u,v\in V_{\reg,h}$ the following discrete bilinear form on one time slab is defined (with $u_-^0\coloneqq 0$)
\begin{align}
	\begin{aligned}\label{discreteformendetail}
		B_{h,n}^\beta(u,v)\coloneqq &\int_{S_h^n}\frac{1}{\alpha_h}\left((1-\beta)\mathring{u}v-\beta u\mathring{v} +(1-\beta)uv\Div_{\Gamma_h}\mathbf{w}+ \mu_d \nabla_{\Gamma_h}u\cdot \nabla_{\Gamma_h}v\right)\dif \sigma_h\\
		&\qquad+\beta\left(R_-^n u_-^n,v_-^n\right)_{\Gamma_h^n(t_n)}+(1-\beta)\left(R_+^{n-1} u_+^{n-1},v_+^{n-1}\right)_{\Gamma_h^n(t_{n-1})}\\
		&\qquad-\left( R_+^{n-1} u^{n-1}_-\circ \Theta_h^{n-1}\circ(\Theta_h^{n})^{-1},v_+^{n-1}\right)_{\Gamma_h^n(t_{n-1})},\\
	\end{aligned}
\end{align} 
which is the discrete analogon of \eqref{betaformel} , and its sum
\begin{equation}
	B_h^\beta(u,v)\coloneqq \sum_{n=1}^N B_{h,n}^{\beta}(u,v).\label{bh_beta}
\end{equation}
\begin{remark} \label{rembeta} \rm Opposite to the bilinear form \eqref{betaformel} the one in \eqref{discreteformendetail} is in general \emph{not} independent of $\beta$. If for fixed arguments $(u,v)$ we take the derivative with respect to $\beta$ and use \eqref{PI-unschoen} we obtain
\begin{align} 
  &\frac{\partial}{\partial \beta}B_{h,n}^\beta(u,v)  =- \int_{S_h^n}\frac{1}{\alpha_h}\left(\mathring{u}v+ u\mathring{v} +uv\Div_{\Gamma_h}\mathbf{w}\right)\dif \sigma_h \nonumber \\ & \qquad \qquad \qquad
		\quad +\left(R_-^n u_-^n,v_-^n\right)_{\Gamma_h^n(t_n)}-\left(R_+^{n-1} u_+^{n-1},v_+^{n-1}\right)_{\Gamma_h^n(t_{n-1})} \label{betaformel1} \\
		&\qquad  =  \int_{S_h^n}\frac{1}{\alpha_h} uv\Div_{\Gamma_h}\mathbf{w} \dif \sigma_h
	  -\sum_{K_S\in  \mathcal{T}_{S_h^n}}\int_{K_S} u_hv_h\Div_{S_h}\left( \frac{1}{\alpha_h} \bP_{S_h}\bw_S\right)\dif \sigma_h \label{betaformel2} \\ & \qquad \quad + \sum_{F\in \mathcal{F}_I^n}\int_{F} u_h v_h \bw_S\cdot \left[\frac{1}{\alpha_h}\right]_{\bnu} \dif F. \label{betaformel3}
\end{align}
Opposite to \eqref{contbeta} this expression is in general not zero. Note that in \eqref{betaformel2}-\eqref{betaformel3} there are no integral terms $(\cdot,\cdot)_{\Gamma_h^n(t_{n})}$. This is due to the fact that in the definition of the method \eqref{discreteformendetail} we use an appropriate weighting with $R$ such that in  the partial integration step, cf. \eqref{PI-unschoen}, these terms vanish. The perturbation terms \eqref{betaformel2}-\eqref{betaformel3} are due to geometric errors and can be controlled by the accuracy of the discrete space-time surface approximation $S_h^n$.
\end{remark}

 In our space-time discretization we use a variant of \textit{volume normal derivative stabilization} \cite{grande2018analysis,burmanembedded}, namely: 
 \begin{equation}
	s(u,v)\coloneqq \xi \int_{Q_{\Theta}^{S}}(\bn_h\cdot \nabla u) (\bn_h\cdot \nabla v)\dif \left(\bx,t\right).\label{NormalStabi}
\end{equation}
 Based on  the literature  we take the parameter range
 \begin{equation}
	h\lesssim\xi \lesssim h^{-1}.\label{stabi_parameter_restrictions}
\end{equation}
We add this stabilization to the bilinear form:
\begin{equation}
	B_h^{\beta,{\rm stab}}(u,v)\coloneqq B_h^\beta(u,v)+s(u,v). \label{stabbilform}
\end{equation}
Before defining the fully discrete problem we need a suitable approximation of the right hand side $f$. 
 For $t\in I_n$ we denote by $\mu_h=\mu_h(t)$ the quotient of the surface measures of $\Gamma(t)$ and $\Gamma_h^n(t)$, i.e. $\mu_h\dif s_h=\dif s$. We assume $f_h\in L^2(S_h)$ to be an approximation of the exact data $f$ satisfying
\begin{equation}
	\norm{f_h-\mu_h f^e}_{L^2(S_h)}\lesssim h^{k+1},\label{f_approx}
\end{equation}
where $k\coloneqq \min\{k_\spa,k_\ti,k_{g,\spa}, k_{g,\ti}\}$. From the literature it is known how such a data approximation can be determined, see e.g. \cite[Remark 4.43]{Diss_Sass} or \cite[Remark 5]{grande2018analysis} in similar settings.
We now define the \emph{space-time discrete variational problem}:
 Given a right-hand side $f_h\in L^2(S_h)$ that satisfies \cref{f_approx}, determine $u_h\in V_{h,\Theta}^{k_\spa,k_\ti}$ such that
 \begin{equation}
 	B_h^{\beta, {\rm stab}}(u_h,v_h)=\int_{S_h} \frac{f_h v_h}{\sqrt{1+V_h^2}}\dif \sigma_h 
 	\quad \text{for all } v_h\in V_{h,\Theta}^{k_\spa,k_\ti}.\label{discreteproblem}
 \end{equation}
 We are particularly interested in the case $\beta=0$ (material derivative on $u$), $\beta =1 $ (material derivative on $v$)  and $\beta=\frac12$ (antisymmetric form).
Implementation aspects of this method are discussed in \Cref{secimplementation} below. In this paper we do not address rigorous discretization error analysis. This topic will be treated in a forthcoming paper. Rigorous error bounds for the case of bi-linear finite element approximations are presented in \cite{Diss_Sass}. 
In the remarks below we discuss a few aspects of this discretization method that we consider to be relevant for a better understanding of the method. 
\begin{remark} \label{RemA} \rm 
We comment on the role of the \textit{volume normal derivative stabilization} $s(\cdot,\cdot)$.
In the setting of a space-time Eulerian method that is based on a trace technique it is very natural to use the standard nodal basis of the outer space $V_{h,\Theta}^{k_\spa,k_\ti}$. It is easy to see that in general (in particular for higher order discretizations) the traces of these outer nodal basis functions on the (space-time) surface become linearly dependent and thus do not form a basis of the trace space.   
Without the stabilization term $s(\cdot,\cdot)$ the bilinear form \cref{stabbilform} only contains surface integrals on $S_h$. In the stabilization term we have an integral over the \emph{volumetric} domain $Q_\Theta^{S}$, consisting of all (deformed) space-time prisms that are cut by the piecewise linear in space surface approximation $S_{\lin}$.  By adding   the volumetric normal derivative we get \emph{control of the variation of the discrete solution in normal direction}.   
This has several advantages, the most important one being that the stiffness matrix (w.r.t. the outer nodal basis) is nonsingular and there is no blow up of the condition number of this matrix. It is well-known that without stabilization, very small cuts of $S_{\lin}$ with the outer triangulation can lead to a very large condition number of the discretization matrix. 
A further relevant property is that the stabilization is \emph{consistent} in the following sense. If one neglects geometric errors then the stabilization term vanishes if we insert  $u=u^e$, i.e., the continuous solution constantly extended in spatial direction. This property is important for  optimal  order discretization accuracy of the method with  stabilization. A detailed analysis of these conditioning and discretization accuracy propertes for the case of a stationary surface $\Gamma$ is given in \cite{LOX_SIAM_2017}. In the present paper we do not address the conditioning of the discretization matrix. In experiments (not presented here) we observe that there is no blow up of condition numbers and that these have a similar behaviour as for the space-time finite element method applied to a parabolic PDE in the volumetric domain $\Omega \times I$. Results that illustrate the discretization   accuracy of the method are given in \Cref{section_numerical_experiments}.
\end{remark}

\begin{remark} \label{RemB} \rm
	In the discrete bilinear form \cref{discreteformendetail} we use a generic function $\alpha_h$, which is a discrete analogon of the  function $\alpha$ appearing in the continuous bilinear form in \cref{prob_cont}. The formula \cref{ho_integral_trafo} suggests that $\alpha_h=\sqrt{1+V_h^2}$ may be a good choice.
	In numerical experiments in \Cref{section_numerical_experiments} we observe that this choice
leads to suboptimal results in certain cases. In \cite[Subsection 5.2.2]{Diss_Sass} a theoretical analysis is given that explains these suboptimal results for the case $k_{g,\spa}=k_{g,\ti}=1$. It is better to use an $\alpha_h$ that is a more accurate approximation of $\alpha$ than the one suggested above. One feasible possibility is the following. Let $\tilde{\phi}_{h} \in V_h^{k_{g,\spa}+1,k_{g,\ti}+1}$ be a one order higher order approximation of the level set function $\phi$ that satisfies
\begin{equation}
	\big|\tilde{\phi}_{h}-\phi\big|_{W^{m,\infty}(Q^S)}\lesssim h^{k_g+2-m},\quad 0\leq m\leq k_g + 1,\label{phi_h_phi_alpha}
\end{equation}
with ${k_g\coloneqq \min\{k_{g,\spa},k_{g,\ti}\}}$, cf. \Cref{phi_assumption_prev}. As an (improved) alternative to  \eqref{Vh_def_ho} we define
\begin{equation}
	\tilde{V}_h\coloneqq - \frac{\partial \tilde{\phi}_h}{\partial t} \norm{\nabla \tilde{\phi}_h}^{-1}.\label{vh_tilde_def}
\end{equation}
The function $\tilde{V}_h$ is   an approximation of order $k_g+1$ of the normal velocity $\bw \cdot \bn$ that does not use a mesh deformation. 
Thus we get the candidate $\alpha_h = \sqrt{1+\tilde{V}_h^2}$, which is smooth on the space-time cylinders in $Q_h^S$ but in general discontinuous across the faces of these cylinders. We therefore apply a straightforward Oswald type projection as follows.  
Let $C(\T_n^\Gamma)$ and $C(Q_{h,n}^S)$ be the space of piecewise continuous functions on the triangulation $\T_n^\Gamma$ and $Q_{h,n}^S$, respectively. By $P_{h,n}:C(\T_n^\Gamma)\rightarrow V_h^{k_{g,\spa}}\restrict{\Omega_n^\Gamma}$ we denote an Oswald-type averaging operator  \cite{oswald1993bpx,lehrenfeld2016high}. We introduce the corresponding space-time projection $\hat P_{h,n}:C(Q_{h,n}^S)\rightarrow V_h^{k_{g,\spa},k_{g,\ti}}$, $n=1, \ldots,N$:
\begin{equation*}
	\hat P_{h,n} v(\bx,t)\coloneqq \sum_{m=0}^{k_{g,\ti}}\mathcal{X}^{\nod}_{\tau_m^n}(t)P_{h,n} v (\bx,\tau_m^n),
\end{equation*}
where $\mathcal{X}^{\nod}_{\tau_m^n}$ denote the nodal basis functions of $\mathcal{P}^{k_{g,\ti}}$ with respect to the discrete points $\tau_m^n\in I_n$.
The choice $\alpha_h=\hat P_{h,n}(\sqrt{1+\tilde{V}^2_h})$, $n=1,\ldots,N$, is a finite element function that is a better approximation of $\alpha^e$ than $\sqrt{1+{V}^2_h}$. This choice is used  in the numerical experiments in \Cref{section_alpha_exps}.  
\end{remark}
\begin{remark} \rm 
We comment  on the function $R$ used in the boundary terms in \cref{discreteformendetail}, which serves as a weighting of the time slab boundary integrals. These weights can be included without significant additional computational costs. The difference $|R-1|$ is controlled by geometric errors. In case of an exact surface we have $R=1$. 
Adding the $R$-weighting leads  to a more consistent formulation in a sense as discussed in \Cref{rembeta}. If we include these weights the boundary terms \eqref{betaformel1} vanish in the partial integration on the discrete surface. 
In the error analysis in \cite[Chapter 5]{Diss_Sass} this cancellation plays a key role in the derivation of optimal discretization error bounds. The analysis does not yield optimal bounds if we replace $R$ by 1. On the other hand, in numerical experiments (not presented here) we obtain optimal order convergence also if we replace $R$ by 1. Below in \Cref{section_numerical_experiments} we use the method with  the $R$-weighting in the boundary terms.
\end{remark}
\begin{remark} \rm 
As mentioned above, cf. \Cref{rembeta}, the discretization method depends on $\beta$. 
Natural choices are $\beta=0$, $\beta=1$, $\beta=\frac12$. The form with $\beta=0$ is ``closest'' to the original strong formulation \eqref{surfactant1}.  In the error analysis, cf.  \cite[Chapter 5]{Diss_Sass},  this leads to a relatively simple consistency error analysis. The method for $\beta=\frac12$ is antisymmetric with respect to the material derivative, which almost immediately leads to a stable method, without using partial integration. The method for $\beta=1$ and with $R$ replaced by $1$ is exactly mass conserving for $k_{g,\spa}=1$ on the discrete level cf. \Cref{remark_exact_mass_conservation} below. We will consider different  $\beta$ values  in the numerical experiments  in \Cref{section_alpha_exps} below.
\end{remark}

\begin{remark} \label{Remtransfer}
\rm
We explain the factor $\Theta_h^{n-1}\circ (\Theta_h^n)^{-1}$ that appears in the last term in \eqref{discreteformendetail}. This term is included to (weakly) transfer the solution from  the previous time step to an initial condition in the current one. As explained in \Cref{remconsist}, the space-time manifold $S_h$ is not necessarily continuous between time slabs. To evaluate $u_-^{n-1}$ on the triangulation deformed by $\Theta_h^n$, the points on $\Gamma_h^n(t_{n-1})$ are mapped to $\Gamma_{\lin}(t_{n-1})$ using $(\Theta_h^n)^{-1}$. Then, using $\Theta_h^{n-1}$, these points are mapped to $\Gamma_h^{n-1}(t_{n-1})$ where $u_-^{n-1}$ is defined on.
\end{remark}
\begin{remark} \label{remstab2} \rm In \cite{olshanskii2014error} another stabilization, namely of the form
 \begin{equation}  
	s_1(u,v)\coloneqq \xi_1\sum_{n=1}^N\int_{t_{n-1}}^{t_n}\int_{\Gamma^n_h(t)}u\dif s_h\int_{\Gamma^n_h(t)}v\dif s_h\dif t,\quad \xi_1\geq 0, 
\end{equation}
is used. This stabilization is needed in the stability analysis presented in that paper.  Based on numerical experiments it seems that this term is not essential for stability of the method.  In this paper, in particular also in the numerical experiments below,  we do not include this stabilization term.
\end{remark}

\subsection{Implementation aspects} \label{secimplementation} 
\begin{figure}[!htbp]
	\centering
	\begin{tikzpicture}[scale = 0.75]
		\newcommand{\Hoehe}{8.0}
		\newcommand{\Tiefe}{2.0}
		\newcommand{\Shift}{3.0}
		\newcommand{\lambdaA}{0.5}
		\newcommand{\PunktA}{(\lambdaA*0+\Shift-\lambdaA*\Shift,\lambdaA*0+\Tiefe-\lambdaA*\Tiefe)}
		\newcommand{\lambdaB}{0.8}
		\newcommand{\PunktB}{(\lambdaB*5,0)}
		\newcommand{\lambdaC}{0.2}
		\newcommand{\PunktC}{(\Shift,\Tiefe+\lambdaC*\Hoehe)}
		\newcommand{\lambdaD}{0.65}
		\newcommand{\PunktD}{(5,\lambdaD*\Hoehe)}
		\newcommand{\lambdaE}{0.5}
		\newcommand{\PunktE}{(5+\Shift,\Tiefe+\lambdaE*\Hoehe)}
		\newcommand{\lambdaF}{0.7}
		\newcommand{\PunktF}{(10,\lambdaF*\Hoehe)}
		\newcommand{\lambdaG}{0.5}
		\newcommand{\PunktG}{(\lambdaG*5+\lambdaG*\Shift+10+\Shift-10*\lambdaG-\lambdaG*\Shift,\Tiefe+\Hoehe)}
		\newcommand{\lambdaH}{0.5}
		\newcommand{\PunktH}{(\lambdaH*10+10+\Shift-\lambdaH*10-\lambdaH*\Shift,\Tiefe-\lambdaH*\Tiefe+\Hoehe)}
		\newcommand{\AtildeShifth}{-0.5}\newcommand{\AtildeShiftv}{1.0}
		\newcommand{\PunktAtilde}{(\lambdaA*0.5+\Shift-\lambdaA*\Shift+\AtildeShifth,\lambdaA*0+\Tiefe-\lambdaA*\Tiefe+\AtildeShiftv)}
		\newcommand{\BtildeShifth}{1.0}\newcommand{\BtildeShiftv}{-0.2}
		\newcommand{\PunktBtilde}{(\lambdaB*5+\BtildeShifth,0+\BtildeShiftv)}
		\newcommand{\GtildeShifth}{-1.0}\newcommand{\GtildeShiftv}{+1.0}
		\newcommand{\PunktGtilde}{(\lambdaG*5+\lambdaG*\Shift+10+\Shift-10*\lambdaG-\lambdaG*\Shift+\GtildeShifth,\Tiefe+\Hoehe+\GtildeShiftv)}
		\newcommand{\HtildeShifth}{1.0}\newcommand{\HtildeShiftv}{-0.5}
		\newcommand{\PunktHtilde}{(\lambdaH*10+10+\Shift-\lambdaH*10-\lambdaH*\Shift+\HtildeShifth,\Tiefe-\lambdaH*\Tiefe+\Hoehe+\HtildeShiftv)}
		
		\newcommand{\ABh}{0.0}\newcommand{\ABv}{-0.2}
		\newcommand{\BAh}{0.0}\newcommand{\BAv}{-0.4}
		\newcommand{\BDh}{0.0}\newcommand{\BDv}{0.2}
		\newcommand{\DBh}{0.0}\newcommand{\DBv}{0.5}
		\newcommand{\ACh}{0.0}\newcommand{\ACv}{0.2}
		\newcommand{\CAh}{0.0}\newcommand{\CAv}{-0.4}
		\newcommand{\CEh}{0.0}\newcommand{\CEv}{0.2}
		\newcommand{\ECh}{0.0}\newcommand{\ECv}{0.2}
		\newcommand{\DFh}{0.0}\newcommand{\DFv}{-0.5}
		\newcommand{\FDh}{0.0}\newcommand{\FDv}{0.0}
		\newcommand{\CDh}{0.0}\newcommand{\CDv}{0.1}
		\newcommand{\DCh}{0.0}\newcommand{\DCv}{0.1}
		\newcommand{\DEh}{0.0}\newcommand{\DEv}{0.2}
		\newcommand{\EDh}{0.0}\newcommand{\EDv}{0.2}
		\newcommand{\EFh}{0.0}\newcommand{\EFv}{-0.3}
		\newcommand{\FEh}{0.0}\newcommand{\FEv}{0.0}
		\newcommand{\EGh}{0.0}\newcommand{\EGv}{-0.3}
		\newcommand{\GEh}{0.0}\newcommand{\GEv}{-0.8}
		\newcommand{\FHh}{0.0}\newcommand{\FHv}{-0.5}
		\newcommand{\HFh}{0.0}\newcommand{\HFv}{-0.3}
		\newcommand{\GHh}{-0.1}\newcommand{\GHv}{-0.05}
		\newcommand{\HGh}{-0.1}\newcommand{\HGv}{-0.05}
		\newcommand{\AAtildeh}{0.0}\newcommand{\AAtildev}{-0.3}
		\newcommand{\AtildeAh}{0.0}\newcommand{\AtildeAv}{-0.2}
		\newcommand{\BBtildeh}{0.0}\newcommand{\BBtildev}{0.0}
		\newcommand{\BtildeBh}{0.0}\newcommand{\BtildeBv}{0.0}
		\newcommand{\HHtildeh}{0.01}\newcommand{\HHtildev}{0.03}
		\newcommand{\HtildeHh}{0.1}\newcommand{\HtildeHv}{0.1}
		\newcommand{\GGtildeh}{0.1}\newcommand{\GGtildev}{0.1}
		\newcommand{\GtildeGh}{0.1}\newcommand{\GtildeGv}{0.1}
		
		\newcommand{\PunktAB}{(0.333*\lambdaA*0+0.333*\Shift-0.333*\lambdaA*\Shift+0.666*\lambdaB*5+\ABh,0.333*\lambdaA*0+0.333*\Tiefe-0.333*\lambdaA*\Tiefe+0.666*0+\ABv)}
		\newcommand{\PunktBA}{(0.666*\lambdaA*0+0.666*\Shift-0.666*\lambdaA*\Shift+0.333*\lambdaB*5+\BAh,0.666*\lambdaA*0+0.666*\Tiefe-0.666*\lambdaA*\Tiefe+0.333*0+\BAv)}
		\newcommand{\PunktBD}{(0.333*\lambdaB*5+0.666*5+\BDh,0.333*0+0.666*\lambdaD*\Hoehe+\BDv)}
		\newcommand{\PunktDB}{(0.666*\lambdaB*5+0.333*5+\DBh,0.666*0+0.333*\lambdaD*\Hoehe+\DBv)}
		\newcommand{\PunktAC}{(0.333*\lambdaA*0+0.333*\Shift-0.333*\lambdaA*\Shift+0.666*\Shift+\ACh,0.333*\lambdaA*0+0.333*\Tiefe-0.333*\lambdaA*\Tiefe+0.666*\Tiefe+0.666*\lambdaC*\Hoehe+\ACv)}
		\newcommand{\PunktCA}{(0.666*\lambdaA*0+0.666*\Shift-0.666*\lambdaA*\Shift+0.333*\Shift+\CAh,0.666*\lambdaA*0+0.666*\Tiefe-0.666*\lambdaA*\Tiefe+0.333*\Tiefe+0.333*\lambdaC*\Hoehe+\CAv)}
		\newcommand{\PunktCE}{(0.333*\Shift+0.666*5+0.666*\Shift+\CEh,0.333*\Tiefe+0.333*\lambdaC*\Hoehe+0.666*\Tiefe+0.666*\lambdaE*\Hoehe+\CEv)}
		\newcommand{\PunktEC}{(0.666*\Shift+0.333*5+0.333*\Shift+\ECh,0.666*\Tiefe+0.666*\lambdaC*\Hoehe+0.333*\Tiefe+0.333*\lambdaE*\Hoehe+\ECv)}
		\newcommand{\PunktDF}{(0.333*5+0.666*10+\DFh,0.333*\lambdaD*\Hoehe+0.666*\lambdaF*\Hoehe+\DFv)}
		\newcommand{\PunktFD}{(0.666*5+0.333*10+\FDh,0.666*\lambdaD*\Hoehe+0.333*\lambdaF*\Hoehe+\FDv)}
		\newcommand{\PunktCD}{(0.333*\Shift+0.666*5+\CDh,0.333*\Tiefe+0.333*\lambdaC*\Hoehe+0.666*\lambdaD*\Hoehe+\CDv)}
		\newcommand{\PunktDC}{(0.666*\Shift+0.333*5+\DCh,0.666*\Tiefe+0.666*\lambdaC*\Hoehe+0.333*\lambdaD*\Hoehe+\DCv)}
		\newcommand{\PunktDE}{(0.666*5+0.666*\Shift+0.333*5+\DEh,0.666*\Tiefe+0.666*\lambdaE*\Hoehe+0.333*\lambdaD*\Hoehe+\DEv)}
		\newcommand{\PunktED}{(0.333*5+0.333*\Shift+0.666*5+\EDh,0.333*\Tiefe+0.333*\lambdaE*\Hoehe+0.666*\lambdaD*\Hoehe+\EDv)}
		\newcommand{\PunktEF}{(0.333*5+0.333*\Shift+0.666*10+\EFh,0.333*\Tiefe+0.333*\lambdaE*\Hoehe+0.666*\lambdaF*\Hoehe+\EFv)}
		\newcommand{\PunktFE}{(0.666*5+0.666*\Shift+0.333*10+\FEh,0.666*\Tiefe+0.666*\lambdaE*\Hoehe+0.333*\lambdaF*\Hoehe+\FEv)}
		\newcommand{\PunktEG}{(0.333*5+0.333*\Shift+0.666*\lambdaG*5+0.666*\lambdaG*\Shift+0.666*10+0.666*\Shift-0.666*10*\lambdaG-0.666*\lambdaG*\Shift+\EGh,0.333*\Tiefe+0.333*\lambdaE*\Hoehe+0.666*\Tiefe+0.666*\Hoehe+\EGv)}
		\newcommand{\PunktGE}{(0.666*5+0.666*\Shift+0.333*\lambdaG*5+0.333*\lambdaG*\Shift+0.333*10+0.333*\Shift-0.333*10*\lambdaG-0.333*\lambdaG*\Shift+\GEh,0.666*\Tiefe+0.666*\lambdaE*\Hoehe+0.333*\Tiefe+0.333*\Hoehe+\GEv)}
		\newcommand{\PunktFH}{(0.333*10+0.666*\lambdaH*10+0.666*10+0.666*\Shift-0.666*\lambdaH*10-0.666*\lambdaH*\Shift+\FHh,0.333*\lambdaF*\Hoehe+0.666*\Tiefe-0.666*\lambdaH*\Tiefe+0.666*\Hoehe+\FHv)}
		\newcommand{\PunktHF}{(0.666*10+0.333*\lambdaH*10+0.333*10+0.333*\Shift-0.333*\lambdaH*10-0.333*\lambdaH*\Shift+\HFh,0.666*\lambdaF*\Hoehe+0.333*\Tiefe-0.333*\lambdaH*\Tiefe+0.333*\Hoehe+\HFv)}
		\newcommand{\PunktGH}{(0.333*\lambdaG*5+0.333*\lambdaG*\Shift+0.333*10+0.333*\Shift-0.333*10*\lambdaG-0.333*\lambdaG*\Shift+0.666*\lambdaH*10+0.666*10+0.666*\Shift-0.666*\lambdaH*10-0.666*\lambdaH*\Shift+\GHh,0.333*\Tiefe+0.333*\Hoehe+0.666*\Tiefe-0.666*\lambdaH*\Tiefe+0.666*\Hoehe+\GHv)}
		\newcommand{\PunktHG}{(0.666*\lambdaG*5+0.666*\lambdaG*\Shift+0.666*10+0.666*\Shift-0.666*10*\lambdaG-0.666*\lambdaG*\Shift+0.333*\lambdaH*10+0.333*10+0.333*\Shift-0.333*\lambdaH*10-0.333*\lambdaH*\Shift+\HGh,0.666*\Tiefe+0.666*\Hoehe+0.333*\Tiefe-0.333*\lambdaH*\Tiefe+0.333*\Hoehe+\HGv)}
		\newcommand{\PunktAAtilde}{(\lambdaA*0+\Shift-\lambdaA*\Shift+0.666*\AtildeShifth+\AAtildeh,\lambdaA*0+\Tiefe-\lambdaA*\Tiefe+0.666*\AtildeShiftv+\AAtildev)}
		\newcommand{\PunktAtildeA}{(\lambdaA*0+\Shift-\lambdaA*\Shift+0.333*\AtildeShifth+\AtildeAh,\lambdaA*0+\Tiefe-\lambdaA*\Tiefe+0.333*\AtildeShiftv+\AtildeAv)}
		\newcommand{\PunktBBtilde}{(\lambdaB*5+0.666*\BtildeShifth+\BBtildeh,0.666*\BtildeShiftv+\BBtildev)}
		\newcommand{\PunktBtildeB}{(\lambdaB*5+0.333*\BtildeShifth+\BtildeBh,0.333*\BtildeShiftv+\BtildeBv)}
		\newcommand{\PunktHHtilde}{(\lambdaH*10+10+\Shift-\lambdaH*10-\lambdaH*\Shift+0.666*\HtildeShifth+\HHtildeh,\Tiefe-\lambdaH*\Tiefe+\Hoehe+0.666*\HtildeShiftv+\HHtildev)}
		\newcommand{\PunktHtildeH}{(\lambdaH*10+10+\Shift-\lambdaH*10-\lambdaH*\Shift+0.333*\HtildeShifth+\HtildeHh,\Tiefe-\lambdaH*\Tiefe+\Hoehe+0.333*\HtildeShiftv+\HtildeHv)}
		\newcommand{\PunktGGtilde}{(\lambdaG*5+\lambdaG*\Shift+10+\Shift-10*\lambdaG-\lambdaG*\Shift+0.666*\GtildeShifth+\GGtildeh,\Tiefe+\Hoehe+0.666*\GtildeShiftv+\GGtildev)}
		\newcommand{\PunktGtildeG}{(\lambdaG*5+\lambdaG*\Shift+10+\Shift-10*\lambdaG-\lambdaG*\Shift+0.333*\GtildeShifth+\GtildeGh,\Tiefe+\Hoehe+0.333*\GtildeShiftv+\GtildeGv)}
		
		\node[left] at (-1.0,0.0) {$t^*_0=t_{n-1}$};
		\node[left] at (-1.0,\Hoehe) {$t^*_3=t_n$};
		\node[left] at (-1.0,1.6) {$t^*_1$};
		\node[left] at (-1.0,5.2) {$t^*_2$};
		\node[left] at (-1.0,4) {$t_l$};
		\newcommand{\FdM}{blue}
		

		\node[below, color=blue] at (5.6,3.2) {$S_{\lin}$};


			\newcommand{\FF}{30}
			\fill[\FdM!\FF] \PunktA -- \PunktB -- \PunktB .. controls \PunktDB and \PunktBD .. \PunktD -- %
			\PunktD .. controls \PunktCD and \PunktDC ..  \PunktC .. controls \PunktAC and \PunktCA .. \PunktA;

		\newcommand{\Radius}{2pt}
		\filldraw[fill=black, ultra thick] \PunktA circle (\Radius);
		\filldraw[fill=black, ultra thick] \PunktB circle (\Radius);
		\filldraw[fill=black, ultra thick] \PunktC circle (\Radius);
		\filldraw[fill=black, ultra thick] \PunktD circle (\Radius);
		\filldraw[fill=black, ultra thick] (4.28,1.6) circle (\Radius);
		
		\draw[ thick, color=\FdM] \PunktA -- \PunktB;
		\draw[thick, color=\FdM]       \PunktA .. controls \PunktCA and \PunktAC .. \PunktC;
		
		\draw[thick, color=\FdM]       \PunktB .. controls \PunktDB and \PunktBD .. \PunktD;
		
		\draw[thick, color=\FdM]       \PunktC .. controls \PunktDC and \PunktCD .. \PunktD;
		\newcommand{\steins}{1.6}
		\draw[thick, color=\FdM] (4.28,1.6) -- (\Shift,\Tiefe+\steins);
		\filldraw[fill=red, ultra thick, color =red ] (4.3,4.7) circle (\Radius);
		\filldraw[fill=red, ultra thick, color = red] (4.72,4) circle (\Radius);
		
		\draw[ultra thick] (0.0,0.0) -- (5.0,0.0);
		\draw[dashed] (0.0,0.0) -- (\Shift,\Tiefe);
		\draw[dashed] (\Shift,\Tiefe) -- (5.0,0.0);
		
		\draw[ultra thin, color = black!50] (0.0,\steins) -- (5.0,\steins);
		\draw[dashed, ultra thin, color = black!50] (0.0,\steins) -- (\Shift,\Tiefe+\steins);
		\draw[dashed, ultra thin, color = black!50] (\Shift,\Tiefe+\steins) -- (5.0,\steins);
		\newcommand{\stzwei}{5.2}
		\draw[ultra thin, color = black!50] (0.0,\stzwei) -- (5.0,\stzwei);
		\draw[dashed,ultra thin, color = black!50] (0.0,\stzwei) -- (\Shift,\Tiefe+\stzwei);
		\draw[ultra thin,dashed, color = black!50] (\Shift,\Tiefe+\stzwei) -- (5.0,\stzwei);
		\newcommand{\stdrei}{4}
		\draw[ultra thin, color = black!50] (0.0,\stdrei) -- (5.0,\stdrei);
		\draw[dashed,ultra thin, color = black!50] (0.0,\stdrei) -- (\Shift,\Tiefe+\stdrei);
		\draw[dashed,ultra thin, color = black!50] (\Shift,\Tiefe+\stdrei) -- (5.0,\stdrei);

		\draw[thick, color=red] (4.72,4) -- (4.3,4.7);
		\draw[ultra thick] (0.0,0.0+\Hoehe) -- (5.0,0.0+\Hoehe);
		\draw[ultra thick] (0.0,0.0+\Hoehe) -- (\Shift,\Tiefe+\Hoehe);
		\draw[ultra thick] (\Shift,\Tiefe+\Hoehe) -- (5.0,0.0+\Hoehe);
		\draw[ultra thick] (0.0,0.0) -- (0.0,0.0+\Hoehe);
		\draw[ultra thick] (5.0,0.0) -- (5.0,0.0+\Hoehe);
		\draw[dashed]      (\Shift,\Tiefe) -- (\Shift,\Tiefe+\Hoehe);
		
		\node[above] at (0.85,9) {$K\times I_n$};
		\node[above, color = red] at (6.5,4) {$K\cap\Gamma_{\lin}(t_l)$};
	\end{tikzpicture}
	\caption{We calculate integrals on $S_{h}$ by using quadrature rules on sufficient spatial planar intersections $K\cap \Gamma_{\lin}(t_l)$.}
	\label{tstern}
\end{figure}
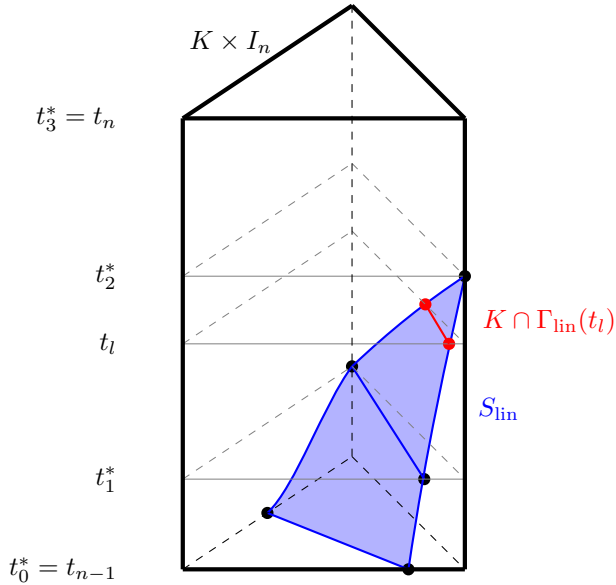

We discuss several  implementation aspects of the discrete problem \cref{discreteproblem}. We first note that all components needed for the implementation of this discretization are available in the open source package \cite{ngsxfem}, an add-on to the finite element library \texttt{NGSolve}$\backslash$\texttt{Netgen} \cite{Schoeberl1997,ngsolve}. The test functions in \cref{discreteproblem} are in the temporal discontinuous finite element space $V_{h,\Theta}^{k_\spa,k_\ti}$, which allows us to implement the method as a \emph{time stepping algorithm}, solving \cref{discreteproblem} time slab after time slab. The solution of the $(n-1)$-th time slab is weakly transferred to $n$-th time slab by the term $(R_+^{n-1} u^{n-1}_-\circ \Theta_h^{n-1}\circ(\Theta_h^{n})^{-1},v_+^{n-1})_{\Gamma_h^n(t_{n-1})}$, $n\in \{1,\dots,N\}$, see \Cref{Remtransfer}. We consider the bilinear form in one fixed time slab  $n$. By $N_n$ we denote the number of degrees of freedom in  that time slab, i.e. the dimension of  the finite element space $V_{h,\Theta}^{k_\spa,k_\ti}\restrict{Q_{\Theta,n}^{S}}$. Let the corresponding basis functions be denoted by $\varphi^n_i$, $i=1,\dots,N_n$. The implementation of the higher order spatial and space-time surface integrals in the discrete bilinear form \cref{discreteformendetail} is done using  quadrature on $\Gamma_{\lin}(\tau)$ at appropriate temporal quadrature points $\tau\in [0,T]$. We illustrate this for the term 
\begin{equation}
	\int_{S_h^n}\frac{1}{\alpha_h}\mathring{\varphi^n_i} \varphi^n_j \dif \sigma_h = \int_{S_h^n}\frac{1}{\alpha_h}\left(\bw_S\cdot \nabla_{S_h} \varphi^n_i\right) \varphi^n_j \dif \sigma_h.\label{impl_one}
\end{equation}
Other integrals in  \cref{discreteproblem} can be treated in a similar way. 
The quadrature of the volume integral in the stabilization bilinear form \cref{NormalStabi} is simpler because it involves integrals over (deformed) prisms. We use the integrals transformation \cref{ho_integral_trafo} and the chain rule to transform \cref{impl_one} to an integral over the tensor product surface $\Gamma_{\lin}(t)\times I_n$.
\begin{align}
	&\int_{S_h^n}\frac{1}{\alpha_h}\left(\bw_S\cdot \nabla_{S_h} \varphi^n_i\right) \varphi^n_j \dif \sigma_h = \int_{t_{n-1}}^{t_n}\int_{\Gamma^n_h(t)}\frac{\sqrt{1+V_h^2}}{\alpha_h}\left(\bw_S\cdot \nabla_{S_h} \varphi^n_i\right) \varphi^n_j \dif s_h\dif t \nonumber\\
	&\qquad=\int_{t_{n-1}}^{t_n}\int_{\Gamma_{\lin}(t)}J_{\alpha}\left(\left(\bP_{S_h}\bw_S\right)\circ \Theta_h^n\cdot \nabla_{(\bx,t)} (\varphi^n_i)\circ\Theta_h^n\right) (\varphi^n_j\circ\Theta_h^n) \dif s_{\lin}\dif t\nonumber\\
	&\qquad=\int_{t_{n-1}}^{t_n}\int_{\Gamma_{\lin}(t)}J_{\alpha}\left(\left(\bP_{S_h}\bw_S\right)\circ \Theta_h^n \cdot D_{(\bx,t)}^{-\intercal}\Theta_h^n\nabla_{(\bx,t)} (\varphi^n_i\circ\Theta_h^n)\right) (\varphi^n_j\circ\Theta_h) \dif s_{\lin}\dif t\nonumber\\
	&\qquad \eqqcolon\int_{t_{n-1}}^{t_n}\int_{\Gamma_{\lin}(t)}g\dif s_{\lin}\dif t,\label{g_example}
\end{align}
where
\begin{equation*}
	J_\alpha\coloneqq \frac{\abs{\det(D\Theta_{h,\spa}^n)}\norm{D^{-\intercal}\Theta_{h,\spa}^n\bn_{\lin}}\sqrt{1+\left(V_h\circ\Theta_h^n\right)^2}}{\alpha_h\circ\Theta_h^n} \quad \text{on } S^n_{\lin}.
\end{equation*}
We briefly comment on computational aspects of the functions appearing in \cref{g_example}. We need the basis functions $\varphi^n_i, \varphi^n_j$ composed with the mesh deformation $\Theta_h^n$. These compositions are finite element functions on the undeformed mesh, as seen in the definition of the parametric finite element space $V_{h,\Theta}^{k_\spa, k_\ti}$ in \cref{Vh_Theta_Def}. This means that $\varphi^n_i\circ \Theta_h^n \in V_h^{k_\spa,k_\ti}\restrict{Q^S_{n}}$. Hence, we only need basis functions of $V_h^{k_\spa, k_\ti}$, cf. \cref{st_fe_space}, to calculate the discretization matrix. The projection $\bP_{S_h}$  depends on the normal $\bn_{S_h}$. Formula for this normal  and  for $V_h\circ \Theta_h^n$  are given  \cref{defintion_nsh_higherorder} and \cref{Vh_def_ho} and are suitable for an efficient implementation. The functions $J_\alpha$ and $D^{-\intercal}_{(\bx,t)}\Theta_h^n$ are defined on $\Gamma_{\lin}(t)$ as well. The sufficiently smooth given functions $\alpha_h$ and $\bw$ are evaluated using $\Theta_h^n$ explicitly. For an efficient and easy implementation of $\Theta_h^n$ it is important that this is a finite element vector function. We refer to \cite[Chapter 5]{preuss2018higher}, where implementation aspects of the parametric mapping $\Theta_h^{n}\in (V_{h,\Theta}^{k_\spa,k_\ti}\restrict{Q_{\Theta,n}^{S}})^4$ are discussed. 
Hence, for a given $t \in I_n$, all functions in the integrand of \cref{g_example} can be efficiently computed on $\Gamma_{\lin}(t)$. We outline how space-time quadrature for the double integral in \eqref{g_example} can be implemented. 
We write
\begin{equation}
	\int_{t_{n-1}}^{t_n}\int_{\Gamma_{\lin}(t)} g \dif s_{\lin}\dif t=\sum_{K\in \T_n^{\Gamma}}\int_{t_{n-1}}^{t_n}\int_{K\cap \Gamma_{\lin}(t)} g \dif s_{\lin}\dif t.\label{g_quadra}
\end{equation}
Similar to \cite[Section 5]{HLP2022} and \cite[Chapter 5]{preuss2018higher} we  split the temporal integral into parts, where the cut topology does not change. 
Let $K_V$ the set of vertices of any $K\in \T$. We define the set of temporal nodes where the surface cuts a vertex of a tetrahedron $K$, enriched with the initial and final times $t_{n-1}$ and $t_n$, i.e. 
\begin{equation*}
	\Sigma_K\coloneqq \left\lbrace t\in [t_{n-1},t_n]: \hat{\phi}_h(\bx_V,t)=0, \bx_V\in K_V\right\rbrace\cup \{t_{n-1}, t_n\}.
\end{equation*}
The roots of the one-dimensional polynomial $t\mapsto\hat{\phi}_h(\bx_V,t)$ can be found by using e.g. the bisection method. We sort the temporal nodes in $\Sigma_K$. Let $\abs{\Sigma_K}$ denote the cardinality of $\Sigma_K$. We define $t_j^*\in \Sigma_K$, $j\in \{0,\dots,\abs{\Sigma_K}-1\}$, such that 
\begin{equation*}
	t_{n-1}=t_0^*\leq \dots \leq t_{\abs{\Sigma_K}-1}^*=t_n.
\end{equation*}
Within $[t_j^*,t_{j+1}^*]$ the cut topology does not change, i.e. the surface $\Gamma_{\lin}(t)$ cuts the same edges and faces of $K$. As a consequence, the integrand $g$ is smooth on ${K\cap \Gamma_{\lin}(t)\times [t_j^*,t_{j+1}^*]}$. We refer to \Cref{tstern} for an illustration. Hence, the mapping
\begin{equation}
	t\mapsto \int_{K\cap \Gamma_{\lin}(t)} g(\bx,t) \dif s_{\lin},\quad t\in  [t_j^*,t_{j+1}^*]\label{mapping_smooth_quad}
\end{equation}
is smooth. If we do not take care of the cut topology, the mapping \cref{mapping_smooth_quad} is in general not smooth, which may lead to a significantly larger quadrature error. Such a  simpler approach in which topology changes are ignored is used in \cite{hansbo2016cut,Zahedi2017}. In \cite[Section 6]{HLP2022}  numerical investigations with these two approaches are presented in a moving domain setting, showing that the naive quadrature approach performs significantly worse. Motivated by these results we take the cut topology into account here. 
We write \cref{g_quadra} as
\begin{equation*}
	\int_{t_{n-1}}^{t_n}\int_{\Gamma_{\lin}(t)} g \dif s_{\lin}\dif t=\sum_{K\in \T_n^{\Gamma}}\sum_{j=0}^{\abs{\Sigma_K}-2}\int_{t_j^*}^{t_{j+1}^*}\int_{K\cap \Gamma_{\lin}(t)} g \dif s_{\lin}\dif t.
\end{equation*}
Depending on $k_\ti$ and $k_{g,\ti}$ we take quadrature nodes $t_l\in [t_j^*,t_{j+1}^*]$, $l=0,\dots, L-1$, $L\in \mathbb{N}$, with corresponding weights $\omega_l$. The number $L$ is chosen sufficiently large, such that the numerical integration does not influence the order of optimal discretization errors. The choice of $L$ is not straightforward, see \cite[Subsection 6.3.1.1]{Diss_Sass}. We thus get a quadrature rule for the time integral of the following form
\begin{equation}
	\int_{t_{n-1}}^{t_n}\int_{\Gamma_{\lin}(t)} g \dif s_{\lin}\dif t \approx \sum_{K\in \T_n^{\Gamma}}\sum_{j=0}^{\abs{\Sigma_K}-2}\sum_{l=0}^L\omega_l\int_{K\cap \Gamma_{\lin}(t_l)} g \dif s_{\lin}.\label{quadrature_integral}
\end{equation}

For a detailed discussion on the implementation of the spatial parametric mappings $\Theta_{h,\tau_m^n}^{n}$ we refer to \cite{lehrenfeld2016high}.
The domain of integration of the remaining spatial integral is the zero level of the piecewise linear level set function $\hat{\phi}_h(\cdot,t_l)$. The whole space-time integral is broken down to spatial integrals with respect to the low order geometry $\Gamma_{\lin}(t_l)$. Every set $K\cap \Gamma_{\lin}(t_l)$ is either a quadrilateral or a triangle. Quadrilaterals are divided into triangles before we use a quadrature rule of sufficient exactness degree ($\approx 2k_s$) on each triangle to approximate integrals on $K\cap \Gamma_{\lin}(t_l)$. For more details concerning the quadrature method outlined above we refer to \cite{HLP2022,Diss_Sass}.

\section{Numerical experiments}\label{section_numerical_experiments}
In this section we present and discuss several numerical experiments. All methods are  implemented in \texttt{ngsxfem} \cite{ngsxfem}, an add-on package to the finite element library \texttt{NGSolve}$\backslash$\texttt{Netgen}, see \cite{Schoeberl1997,ngsolve}. To solve the arising linear systems, we use the direct solver \texttt{Pardiso} which is part of the \texttt{Intel MKL} library \cite{intel}. To visualize our results we use the software \texttt{ParaView}, see \cite{ahrens2005paraview}. In \Cref{section_alpha_exps} we illustrate the influence of the factor $\alpha_h$ on convergence results by considering a piecewise linear and piecewise higher order variant. In \Cref{section_smooth_geometry} we present results for two different smoothly evolving surfaces. We then show results of  numerical experiments for an evolving surface with topological singularities in \Cref{section_top_singular}. Both spatially two- and three-dimensional settings are considered. 

We  explain general settings used in the experiments. We solve the discrete problem \cref{discreteproblem}, and set the diffusion coefficient to $\mu_d=1$. We consider $\beta\in \{0,\frac12,1\}$ and the stabilization parameter is taken as $\xi\in \{h, h^{-1}\}$. All test cases employ velocity fields that are not divergence-free, i.e. $\Div_{\Gamma}\bw\neq 0$. 
As input for the mesh deformation we use the nodal space-time interpolation $\phi_h=I_h(\phi)\in V_h^{k_{g,\spa},k_{g,\ti}}$ of a given level set function $\phi$,  resulting in the discrete space-time manifold $S_h$, cf. \cref{sh_def}. In \Cref{section_smooth_geometry} we prescribe a smooth solution $u$ in the neighborhood $U$ and calculate the corresponding right hand side $f$. Different from \Cref{prob_cont} we allow a non-zero initial condition $u_0$. In \Cref{section_top_singular} we are not able to compute a right-hand side to a prescribed solution $u$, since there are geometric singularities. Instead, we take a nonzero initial solution $u_0$ on $\Gamma(0)$ and a source term $f=0$. As  initial condition for the discrete problem \cref{discreteproblem} we take the parametric interpolation
\begin{equation*}
	(u_h)^0_-(\bx):=\mathcal{I}_{k_\spa}(u(\Theta_{h,0}(\cdot),0))(\Theta_{h,0}^{-1}(\bx)), \quad \bx\in \Gamma_h(0).
\end{equation*}
The initial condition is weakly enforced via the integral over the deformed $\Gamma_h(0)$ contained in the bilinear form defined in \cref{discreteformendetail}. The right-hand side is taken as $f_h=f$ on $S_h$, which can be shown to be sufficiently accurate, see e.g. \cite[Remark 4.43]{Diss_Sass} or \cite[Remark 5]{grande2014eulerian} in similar settings.
We now explain the error measures that we use.
For this we  introduce the jump brackets on the discrete time slab boundaries
\begin{equation}
	[v]_h^n\coloneqq v_+^n-v_-^n\circ \Theta_h^n \circ (\Theta_h^{n+1})^{-1}\in L^2(\Gamma_{h}^{n+1}(t_n)),\quad n\in\{0,\dots,{N-1}\},\label{discrete_jumps}
\end{equation}
and define the energy norm
\begin{align}
	\begin{aligned}\label{energy_norm_ho_def}
		\triplenorm v \triplenorm_{h}^2\coloneqq &\max_{n\in\{1,\dots,N\}}\norm{v_-^n}^2_{L^2(\Gamma_{h}^n(t_n))}+ {\sum_{n=1}^N}\norm{[v]_h^{n-1}}^2_{L^2(\Gamma_{h}^n(t_{n-1}))}+\norm{v}^2_{L^2(S_h)}\\
		&\qquad+ \norm{\nabla_{\Gamma_h}v}_{L^2(S_h)}^2 
		+ \xi\norm{\bn_h\cdot \nabla v}_{L^2(Q_{\Theta}^{S})}^2. 
	\end{aligned}
\end{align}
For measuring the discretization accuracy we will use the energy norm error $\triplenorm u^e - u_h \triplenorm_{h}$. For this measure to be useful it is important that $u^e$  is chosen as the  constant extension of $u$ in spatial normal direction.  Otherwise the bulk error $\xi^\frac12 \norm{\bn_h\cdot \nabla (u^e-u_h)}_{L^2(Q_{\Theta}^{S})}$ does not necessarily converge with optimal order. In \Cref{section_alpha_exps,section_smooth_geometry} such a constant extension is algebraically computable and $u^e$ is chosen accordingly. In \Cref{section_twisting_torus} we consider a more complicated geometry for which this constant normal extension is not readily available.  In that case we consider the modified energy norm
\begin{equation}
	\triplenorm u^e- u_h \triplenorm_{h,*}^2\coloneqq \triplenorm u^e- u_h \triplenorm_{h}^2-\xi\norm{\bn_h\cdot\nabla (u^e-u_h) }_{L^2(Q_{\Theta}^{S})}^2\label{energy_norm_wonormal}
\end{equation}
which requires $u^e$ values only on $S_h$. In this case it is not necessery to use a constant  $u$ extension in normal direction. Instead we use  a natural, smooth extension $u^e$.   We also measure the error in a discrete $L^\infty L^2$-norm 
\begin{equation*}
	\norm{u^e-u_h}_{L^\infty L^2}\coloneqq \max_{n\in\{1,\dots,N\}}\norm{u^e-u_h}_{L^2(\Gamma^n_h(t_n))}.
\end{equation*}
We use the following simple spatial mesh refinement strategy.  We start with  a coarse unstructured triangular/tetrahedral quasi-uniform triangulation of $\Omega$, with an initial mesh size   denoted by $h_{\mathrm{init}}$. We refine the mesh uniformly to obtain the spatial refinement levels $l_\spa\in \{0,\dots,8\}$, such that $h=h_{\mathrm{init}}2^{-l_\spa}$. Using a marked-edge bisection method, we only refine simplices that are intersected by the surface. Taking an initial uniform time step size $\Delta t_{\mathrm{init}}$, we halve this value for each temporal refinement level $l_\ti\in \{0,\dots,7\}$, i.e. $\Delta t=\Delta t_{\mathrm{init}}2^{-l_\ti}$. To assess the order of convergence in space, time and space-time, we calculate the spatial, temporal and space-time experimental order of convergence, which we denote as $\mathrm{eoc}_\spa$, $\mathrm{eoc}_\ti$ and $\mathrm{eoc}_{\spati}$ respectively. To do so we take two errors $e_{l}$ and $e_{l+1}$, $l\in \{l_\ti,l_\spa\}$, of subsequent refinement levels $l$ and $l+1$ and calculate $\mathrm{eoc}=\frac{1}{\log(2)}\log\left(\frac{e_{l}}{e_{l+1}}\right)$. For the calculation of $\mathrm{eoc}_{\spati}$ we consider one refinement in both space and time as subsequent refinement.
To compute $\mathrm{eoc}_\spa$ and $\mathrm{eoc}_\ti$ we consider the errors of the finest temporal and the finest spatial discretizations used in the specific experiment, respectively.

\subsection{On the factor $\alpha_h$ in the discrete bilinear form}\label{section_alpha_exps}
\begin{figure}[!htbp]
	\begin{minipage}{\textwidth}
		\vspace{-0.5cm}
	\includegraphics[width=\textwidth]{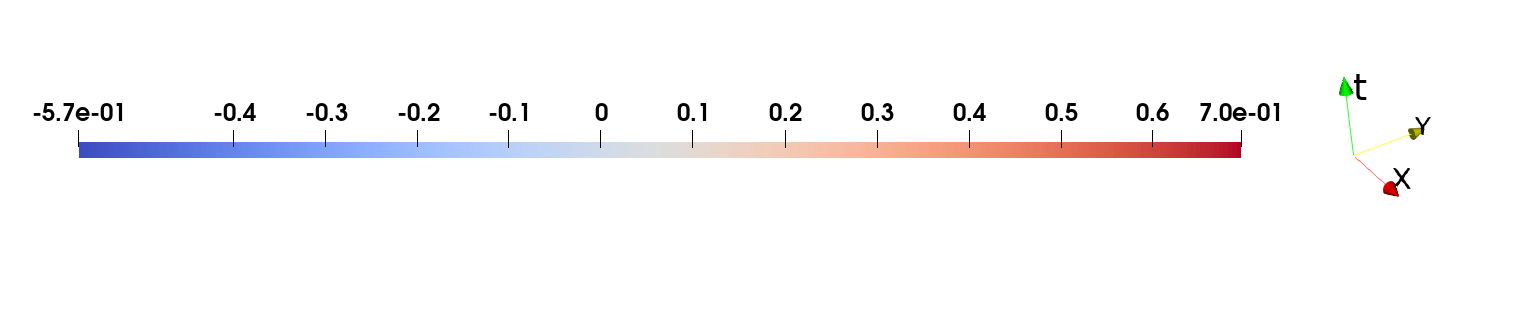}   
	\vspace{-1.5cm}
\end{minipage}
	\begin{minipage}{\textwidth}
	\includegraphics[width=\textwidth]{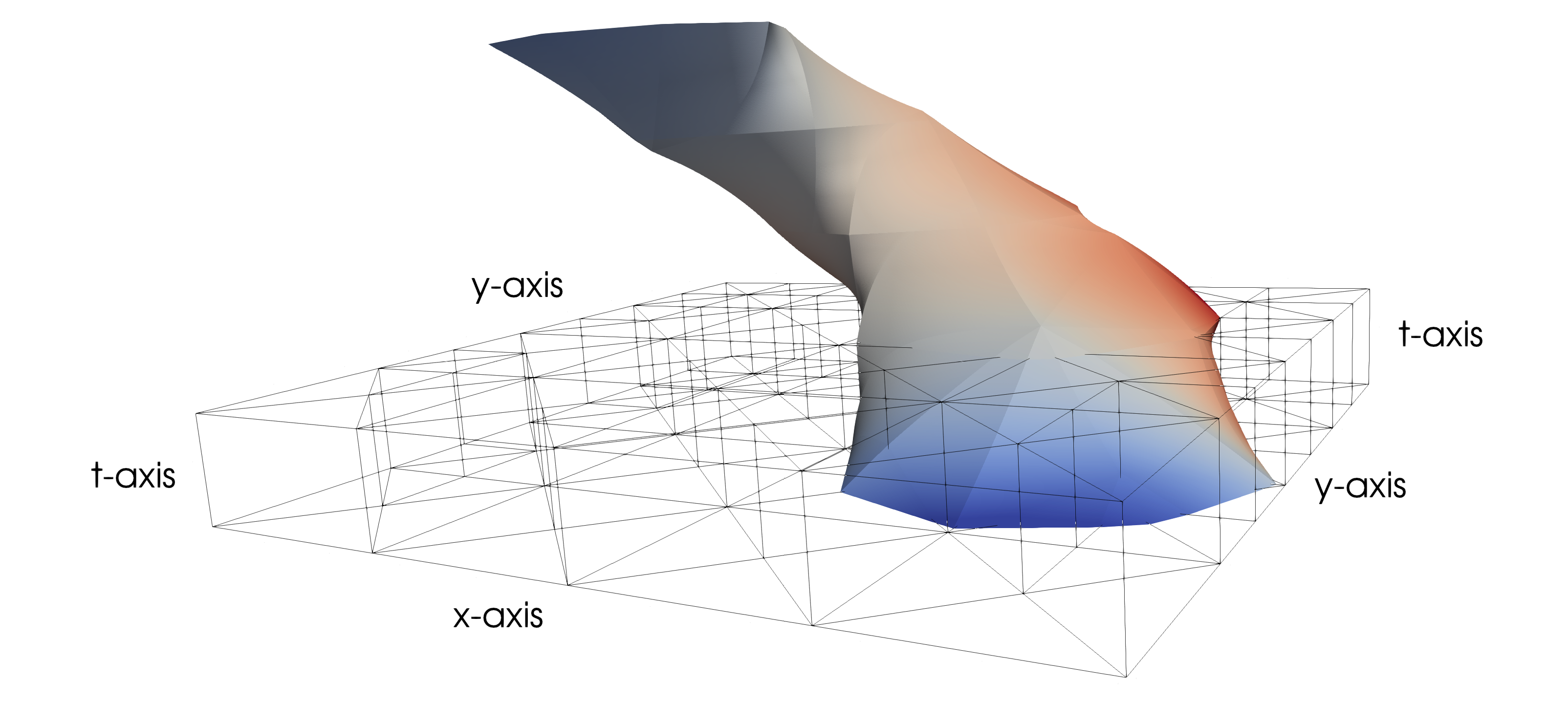}    
\end{minipage}
	\caption{Moving circle: a sketch of the space-time manifold $S_h$ in the space-time prismatic triangulation with $N=4$ time slabs. Here, we have $h=2^{-2}$ and $k_{g_\spa}=k_{g,\ti}=1$. We see the discrete solution $u_h$ on a circle moving through the triangular background mesh.}
	\label{figure_circle_background}
\end{figure}

We consider a shrinking circle that is moving through the background mesh with a time- and space-dependent velocity field. The center of the circle moves along a semicircle in the mesh. As outer domain we take $\Omega=(-1,1)^2$ and as time interval we choose $[0,1]$. On $Q=\Omega\times [0,1]$ the evolution of the surface $\Gamma(t)$ is described by the zero level of
\begin{equation*}
	\phi(\bx,t)=\sqrt{\left(x-\frac{\cos(\pi t)}{2}\right)^2+\left(y-\frac{\sin(\pi t)}{2}\right)^2}-\frac{9e^{-\frac{t}{4}}}{20}, \quad (\bx,t)=(x,y,t)\in Q.
\end{equation*}
On $Q$ we construct a normal velocity field that transports $\Gamma(t)$ by
\begin{equation*}
	\bw_n=-\ddt{\phi}{t}\nabla\phi. 
\end{equation*}
Note that $\norm{\nabla\phi(\bx,t)}=1$ for all $(\bx,t)\in Q$, since $\phi(\cdot, t)$ is a signed distance function to $\Gamma(t)$. We do not want to restrict to the special case of a problem with a strictly normal velocity field $\bw$. Therefore we add a tangential part 
$\bw_{\mathrm{tang}}=\frac{1}{2}\left((\nabla\phi)_2, -(\nabla\phi)_1\right)^\intercal
$ to $\bw_n$, where $(\nabla\phi)_i$ denotes the $i$th component of $\nabla \phi$.
Since $\bw_n\cdot \nabla \phi=-\ddt{\phi}{t}$ and $\bw_{\mathrm{tang}}\cdot \nabla \phi=0$ we conclude that $\bw=\bw_n+\bw_{\mathrm{tang}}$ satisfies the level set equation $\ddt{\phi}{t}+\bw\cdot\nabla\phi=0$. We use this $\bw$ for the experiments and prescribe the exact solution as
\begin{equation*}
	u(\bx,t)=y(x^2+1)e^{-t}, \quad (\bx,t)\in S.
\end{equation*}
We take an initial mesh size of $h_{\mathrm{init}}=2^{-3}$ and an initial time step size as $\Delta t_{\mathrm{init}}=2^{-3}$. In \Cref{figure_circle_background} we illustrate the space-time manifold $S_h$ in the space-time prismatic background mesh for relatively coarse discretization parameters.
In \Cref{fig_alpha_h_left} we present numerical results for $k=1$. Note that for $k=1$ the mesh-transformation reduces to the identity: $\Theta_h^n={\rm id}$.   A seemingly natural choice is the approximation $\alpha_h=\sqrt{1+V_h^2}$, where $V_h$ is a piecewise bi-linear approximation of the normal velocity field $\bw\cdot \bn$, cf. \cref{Vh_def_ho}. It is natural in the sense that $\sqrt{1+V_h^2}$ is the transformation factor between an integral over $S_h^n$, i.e. $\int_{S_h^n}\cdot \dif \sigma_h$, and the corresponding iterated integral $\int_{t_{n-1}}^{t_n}\int_{\Gamma_h(t)}\cdot \dif s_h\dif t$, cf. \Cref{theo_discrete_trafo}. In the continuous setting this factor is $\alpha$, which we approximate by $\alpha_h$. In \Cref{fig_alpha_h_right} results for the piecewise bi-quadratic variant $\alpha_h=\hat P_{h,n}(\sqrt{1+\tilde{V}_h^2})$, cf. \Cref{RemB}, are shown. In \Cref{fig_alpha_h_ho} we  present corresponding results for $k=3$. 
\begin{figure}[!htbp]
	\begin{subfigure}[t]{0.5\textwidth}
		\centering
		\begin{tikzpicture}[scale=0.73]
			\def\vara{0.5}
			\begin{semilogyaxis}[xlabel={Refinement level $l_\ti=l_{\spa}$}, ylabel={$\triplenorm{u^e-u_h}\triplenorm_h$},legend style={ at={(0.5,1.02)}, anchor=south, legend columns =2}, legend cell align=left]
				\addplot table [x=energylist_x, y=-, col sep=comma] {Experiments/EnergyError_DiagKonvTable_Diss1_alpha_h_beta0_stabiwith_h_SpatialNormalVolume_userh_1_tmax_9_tstart_3_xmax_9_xstart_3_ks_1_kt_1_kgs_1_kgt_1_gp_0_tempquadorder_4_tempquadorderDiffu_4.csv};
									\addplot table [x=energylist_x, y=-, col sep=comma] {Experiments/EnergyError_DiagKonvTable_Diss1_alpha_h_beta0_stabiwith_1h_SpatialNormalVolume_userh_1_tmax_9_tstart_3_xmax_9_xstart_3_ks_1_kt_1_kgs_1_kgt_1_gp_0_tempquadorder_4_tempquadorderDiffu_4.csv};
				\addplot +[mark repeat =2, mark phase =1]table [x=energylist_x, y=-, col sep=comma] {Experiments/EnergyError_DiagKonvTable_Diss1_alpha_h_beta05_stabiwith_h_SpatialNormalVolume_userh_1_tmax_9_tstart_3_xmax_9_xstart_3_ks_1_kt_1_kgs_1_kgt_1_gp_0_tempquadorder_4_tempquadorderDiffu_4.csv};
					\addplot +[mark repeat =2, mark phase =2]table [x=energylist_x, y=-, col sep=comma] {Experiments/EnergyError_DiagKonvTable_Diss1_alpha_h_beta05_stabiwith_1h_SpatialNormalVolume_userh_1_tmax_9_tstart_3_xmax_9_xstart_3_ks_1_kt_1_kgs_1_kgt_1_gp_0_tempquadorder_4_tempquadorderDiffu_4.csv};

				\addplot +[mark repeat =2, mark phase =1]table [x=energylist_x, y=-, col sep=comma] {Experiments/EnergyError_DiagKonvTable_Diss1_alpha_h_beta1_stabiwith_h_SpatialNormalVolume_userh_1_tmax_9_tstart_3_xmax_9_xstart_3_ks_1_kt_1_kgs_1_kgt_1_gp_0_tempquadorder_4_tempquadorderDiffu_4.csv};
				\addplot +[mark repeat =2, mark phase =2]table [x=energylist_x, y=-, col sep=comma] {Experiments/EnergyError_DiagKonvTable_Diss1_alpha_h_beta1_stabiwith_1h_SpatialNormalVolume_userh_1_tmax_9_tstart_3_xmax_9_xstart_3_ks_1_kt_1_kgs_1_kgt_1_gp_0_tempquadorder_4_tempquadorderDiffu_4.csv};
				\addplot[dotted,line width=0.75pt] coordinates { 
					(0,\vara) (1,\vara*0.5) (2,\vara*0.25) (3,\vara*0.125) (4,\vara*0.0625)(5,\vara*0.03125) (6,\vara*0.015625)
				};
				\legend{{$\beta=0, \xi=h$}, {$\beta=0, \xi=h^{-1}$},{$\beta=\frac12, \xi=h$},{$\beta=\frac12, \xi=h^{-1}$}, {$\beta=1, \xi=h$},{$\beta=1, \xi=h^{-1}$}, $\mathcal{O}(h)$}
			\end{semilogyaxis}
		\end{tikzpicture}
		\caption{$\alpha_h=(1+V_h^2)^{1/2}$}
		\label{fig_alpha_h_left} 
	\end{subfigure}\hfill
	\begin{subfigure}[t]{0.5\textwidth}
		\centering
		\begin{tikzpicture}[scale=0.73]
			\def\vara{0.4}
			\begin{semilogyaxis}[xlabel={Refinement level $l_\ti=l_{\spa}$}, ylabel={$\triplenorm{u^e-u_h}\triplenorm_h$},legend style={ at={(0.5,1.02)}, anchor=south, legend columns =2}, legend cell align=left]
				\addplot +[mark repeat =3, mark phase =1]table [x=energylist_x, y=-, col sep=comma] {Experiments/EnergyError_DiagKonvTable_Diss1_alpha_h_interpol_ho_oswald_beta0_stabiwith_h_SpatialNormalVolume_userh_1_tmax_9_tstart_3_xmax_9_xstart_3_ks_1_kt_1_kgs_1_kgt_1_gp_0_tempquadorder_4_tempquadorderDiffu_4.csv};
			\addplot +[mark repeat =3, mark phase =1]table [x=energylist_x, y=-, col sep=comma] {Experiments/EnergyError_DiagKonvTable_Diss1_alpha_h_interpol_ho_oswald_beta0_stabiwith_1h_SpatialNormalVolume_userh_1_tmax_9_tstart_3_xmax_9_xstart_3_ks_1_kt_1_kgs_1_kgt_1_gp_0_tempquadorder_4_tempquadorderDiffu_4.csv};
			\addplot +[mark repeat =3, mark phase =2]table [x=energylist_x, y=-, col sep=comma] {Experiments/EnergyError_DiagKonvTable_Diss1_alpha_h_interpol_ho_oswald_beta05_stabiwith_h_SpatialNormalVolume_userh_1_tmax_9_tstart_3_xmax_9_xstart_3_ks_1_kt_1_kgs_1_kgt_1_gp_0_tempquadorder_4_tempquadorderDiffu_4.csv};
			\addplot +[mark repeat =3, mark phase =2]table [x=energylist_x, y=-, col sep=comma] {Experiments/EnergyError_DiagKonvTable_Diss1_alpha_h_interpol_ho_oswald_beta05_stabiwith_1h_SpatialNormalVolume_userh_1_tmax_9_tstart_3_xmax_9_xstart_3_ks_1_kt_1_kgs_1_kgt_1_gp_0_tempquadorder_4_tempquadorderDiffu_4.csv};
			\addplot +[mark repeat =3, mark phase =3]table [x=energylist_x, y=-, col sep=comma] {Experiments/EnergyError_DiagKonvTable_Diss1_alpha_h_interpol_ho_oswald_beta1_stabiwith_h_SpatialNormalVolume_userh_1_tmax_9_tstart_3_xmax_9_xstart_3_ks_1_kt_1_kgs_1_kgt_1_gp_0_tempquadorder_4_tempquadorderDiffu_4.csv};
			\addplot +[mark repeat =3, mark phase =3]table [x=energylist_x, y=-, col sep=comma] {Experiments/EnergyError_DiagKonvTable_Diss1_alpha_h_interpol_ho_oswald_beta1_stabiwith_1h_SpatialNormalVolume_userh_1_tmax_9_tstart_3_xmax_9_xstart_3_ks_1_kt_1_kgs_1_kgt_1_gp_0_tempquadorder_4_tempquadorderDiffu_4.csv};
			\addplot[dotted,line width=0.75pt] coordinates { 
				(0,\vara) (1,\vara*0.5) (2,\vara*0.25) (3,\vara*0.125) (4,\vara*0.0625)(5,\vara*0.03125) (6,\vara*0.015625)
			};
				\legend{{$\beta=0, \xi=h$}, {$\beta=0, \xi=h^{-1}$},{$\beta=\frac12, \xi=h$},{$\beta=\frac12, \xi=h^{-1}$}, {$\beta=1, \xi=h$},{$\beta=1, \xi=h^{-1}$},  $\mathcal{O}(h)$}
		\end{semilogyaxis}
	\end{tikzpicture}
\caption{$\alpha_h=\hat P_{h,n}((1+\tilde{V}_h^2)^{1/2})$}
		\label{fig_alpha_h_right} 
	\end{subfigure}
	\caption{ Moving circle: convergence results varying $\alpha_h$, $\beta$ and $\xi$ in the bilinear forms for $k=1$.}
	\label{fig_alpha_h}
\end{figure}
\begin{figure}[!htbp]
		\centering
	\begin{subfigure}[b]{0.5\textwidth}
		\centering
		\begin{tikzpicture}[scale=0.73]
			\def\varb{0.02}
			\def\varc{0.01}
			\begin{semilogyaxis}[xlabel={Refinement level $l_\ti=l_{\spa}$}, ylabel={$\triplenorm{u^e-u_h}\triplenorm_h$},legend style={ at={(0.5,1.02)}, anchor=south, legend columns =2}, legend cell align=left]
				\addplot table [x=energylist_x, y=-, col sep=comma] {Experiments/EnergyError_DiagKonvTable_Diss1_alpha_h_beta0_stabiwith_h_SpatialNormalVolume_userh_1_tmax_9_tstart_3_xmax_9_xstart_3_ks_3_kt_3_kgs_3_kgt_3_gp_0_tempquadorder_12_tempquadorderDiffu_12.csv};
				\addplot table [x=energylist_x, y=-, col sep=comma] {Experiments/EnergyError_DiagKonvTable_Diss1_alpha_h_beta0_stabiwith_1h_SpatialNormalVolume_userh_1_tmax_9_tstart_3_xmax_9_xstart_3_ks_3_kt_3_kgs_3_kgt_3_gp_0_tempquadorder_12_tempquadorderDiffu_12.csv};
				\addplot +[mark repeat =2, mark phase =1]table [x=energylist_x, y=-, col sep=comma] {Experiments/EnergyError_DiagKonvTable_Diss1_alpha_h_beta05_stabiwith_h_SpatialNormalVolume_userh_1_tmax_9_tstart_3_xmax_9_xstart_3_ks_3_kt_3_kgs_3_kgt_3_gp_0_tempquadorder_12_tempquadorderDiffu_12.csv};
				\addplot+[mark repeat =2, mark phase =2] table [x=energylist_x, y=-, col sep=comma] {Experiments/EnergyError_DiagKonvTable_Diss1_alpha_h_beta05_stabiwith_1h_SpatialNormalVolume_userh_1_tmax_9_tstart_3_xmax_9_xstart_3_ks_3_kt_3_kgs_3_kgt_3_gp_0_tempquadorder_12_tempquadorderDiffu_12.csv};
				\addplot+[mark repeat =2, mark phase =1] table [x=energylist_x, y=-, col sep=comma] {Experiments/EnergyError_DiagKonvTable_Diss1_alpha_h_beta1_stabiwith_h_SpatialNormalVolume_userh_1_tmax_9_tstart_3_xmax_9_xstart_3_ks_3_kt_3_kgs_3_kgt_3_gp_0_tempquadorder_12_tempquadorderDiffu_12.csv};
				\addplot+[mark repeat =2, mark phase =2] table [x=energylist_x, y=-, col sep=comma] {Experiments/EnergyError_DiagKonvTable_Diss1_alpha_h_beta1_stabiwith_1h_SpatialNormalVolume_userh_1_tmax_9_tstart_3_xmax_9_xstart_3_ks_3_kt_3_kgs_3_kgt_3_gp_0_tempquadorder_12_tempquadorderDiffu_12.csv};
				\addplot[dotted,line width=0.75pt] coordinates { 
	(0,\varb) (1,\varb*0.5*0.5) (2,\varb*0.25*0.25) (3,\varb*0.125*0.125) (4,\varb*0.0625*0.0625)(5,\varb*0.03125*0.03125) (6,\varb*0.015625*0.015625)
};
	\addplot[dashed,line width=0.75pt] coordinates { 
	(0,\varc) (1,\varc*0.5*0.5*0.5) (2,\varc*0.25*0.25*0.25) (3,\varc*0.125*0.125*0.125) (4,\varc*0.0625*0.0625*0.0625)(5,\varc*0.03125*0.03125*0.03125) (6,\varc*0.015625*0.015625*0.015625)
};
				\legend{{$\beta=0, \xi=h$}, {$\beta=0, \xi=h^{-1}$},{$\beta=\frac12, \xi=h$},{$\beta=\frac12, \xi=h^{-1}$}, {$\beta=1, \xi=h$},{$\beta=1, \xi=h^{-1}$},  $\mathcal{O}(h^2)$,$\mathcal{O}(h^3)$}
			\end{semilogyaxis}
		\end{tikzpicture}
\caption{$\alpha_h=(1+V_h^2)^{1/2}$}
		\label{fig_alpha_h_ho_left} 
	\end{subfigure}\hfill
	\begin{subfigure}[b]{0.5\textwidth}
		\centering
		\begin{tikzpicture}[scale=0.73]
			\def\varb{0.2}
			\def\varc{0.08}
			\begin{semilogyaxis}[xlabel={Refinement level $l_t=l_{\bx}$}, ylabel={$\triplenorm{u^e-u_h}\triplenorm_h$},legend style={ at={(0.5,1.02)}, anchor=south, legend columns =2}, legend cell align=left]
				\addplot table [x=energylist_x, y=-, col sep=comma] {Experiments/EnergyError_DiagKonvTable_Diss1_alpha_h_interpol_ho_oswald_beta0_stabiwith_h_SpatialNormalVolume_userh_1_tmax_9_tstart_3_xmax_9_xstart_3_ks_3_kt_3_kgs_3_kgt_3_gp_0_tempquadorder_12_tempquadorderDiffu_12.csv};
				\addplot +[mark repeat =4, mark phase =1]table [x=energylist_x, y=-, col sep=comma] {Experiments/EnergyError_DiagKonvTable_Diss1_alpha_h_interpol_ho_oswald_beta0_stabiwith_1h_SpatialNormalVolume_userh_1_tmax_9_tstart_3_xmax_9_xstart_3_ks_3_kt_3_kgs_3_kgt_3_gp_0_tempquadorder_12_tempquadorderDiffu_12.csv};
				\addplot table [x=energylist_x, y=-, col sep=comma] {Experiments/EnergyError_DiagKonvTable_Diss1_alpha_h_interpol_ho_oswald_beta05_stabiwith_h_SpatialNormalVolume_userh_1_tmax_9_tstart_3_xmax_9_xstart_3_ks_3_kt_3_kgs_3_kgt_3_gp_0_tempquadorder_12_tempquadorderDiffu_12.csv};
				\addplot +[mark repeat =4, mark phase =2]table [x=energylist_x, y=-, col sep=comma] {Experiments/EnergyError_DiagKonvTable_Diss1_alpha_h_interpol_ho_oswald_beta05_stabiwith_1h_SpatialNormalVolume_userh_1_tmax_9_tstart_3_xmax_9_xstart_3_ks_3_kt_3_kgs_3_kgt_3_gp_0_tempquadorder_12_tempquadorderDiffu_12.csv};
				\addplot +[mark repeat =4, mark phase =1]table [x=energylist_x, y=-, col sep=comma] {Experiments/EnergyError_DiagKonvTable_Diss1_alpha_h_interpol_ho_oswald_beta1_stabiwith_h_SpatialNormalVolume_userh_1_tmax_9_tstart_3_xmax_9_xstart_3_ks_3_kt_3_kgs_3_kgt_3_gp_0_tempquadorder_12_tempquadorderDiffu_12.csv};
				\addplot +[mark repeat =4, mark phase =3]table [x=energylist_x, y=-, col sep=comma] {Experiments/EnergyError_DiagKonvTable_Diss1_alpha_h_interpol_ho_oswald_beta1_stabiwith_1h_SpatialNormalVolume_userh_1_tmax_9_tstart_3_xmax_9_xstart_3_ks_3_kt_3_kgs_3_kgt_3_gp_0_tempquadorder_12_tempquadorderDiffu_12.csv};
				\addplot[dashed,line width=0.75pt] coordinates { 
				(0,\varc) (1,\varc*0.5*0.5*0.5) (2,\varc*0.25*0.25*0.25) (3,\varc*0.125*0.125*0.125) (4,\varc*0.0625*0.0625*0.0625)(5,\varc*0.03125*0.03125*0.03125) (6,\varc*0.015625*0.015625*0.015625)
			};
				\legend{{$\beta=0, \xi=h$}, {$\beta=0, \xi=h^{-1}$},{$\beta=\frac12, \xi=h$},{$\beta=\frac12, \xi=h^{-1}$}, {$\beta=1, \xi=h$},{$\beta=1, \xi=h^{-1}$}, $\mathcal{O}(h^3)$}
			\end{semilogyaxis}
		\end{tikzpicture}
		\caption{$\alpha_h=\hat P_{h,n}((1+\tilde{V}_h^2)^{1/2})$}
		\label{fig_alpha_h_ho_right} 
	\end{subfigure}
	\caption{Moving circle: convergence results varying $\alpha_h$, $\beta$ and $\xi$ in the bilinear forms for $k=3$.}
\label{fig_alpha_h_ho}
\end{figure}
In \Cref{fig_alpha_h_ho_left,fig_alpha_h_left} we observe optimal order convergence for the methods with $\beta=0$ and one order less for the methods with $\beta\in\{\frac12,1\}$. A heuristic explanation for this difference in order of convergence is as follows. 
Note that the methods with $\beta\in \{\frac12,1\}$ contain a time derivative of the test function and due to this one needs partial integration in the consistency error analysis. Hence, jumps of $\alpha_h$ occur in the consistency error, see \Cref{PI_Theorem}. In the method with $\beta=0$ partial integration is not needed in the consistency error analysis (only in the stability analysis). The relatively poor $\alpha^e$ approximation $\sqrt{1+V_h^2}$, which  is  discontinuous across surface element faces, causes large jumps of $\alpha_h$ that lead  to a suboptimal order term in the consistency error. This does not happen for the case $\beta =0$. For the more accurate approximation $\alpha_h=\hat P_{h,n}(\sqrt{1+\tilde{V}^2_h})$ the jump terms on interior boundary faces of surface elements, i.e. on $\mathcal{F}_I^n$, vanish. Furthermore, corresponding jumps between time slabs are sufficiently small to allow an optimal order consistency error. 

If not stated differently, we restrict to $\beta=0$, $\xi=h$ and $\alpha_h=\sqrt{1+V_h^2}$ in the experiments below.

\subsection{Smoothly evolving surfaces}\label{section_smooth_geometry}
In this section we consider  spatially two-  and three-dimensional examples, both  with smooth geometries.
\subsubsection{Moving circle}
We take the same $\Omega, \phi$, $\bw$, $T$ and $u$ as in \Cref{section_alpha_exps}. We choose the initial mesh size  $h_{\mathrm{init}}=2^{-2}$ and the initial time step size  $\Delta t_{\mathrm{init}}=2^{-2}$. In \Cref{fig_circling} and \Cref{table_circling_energy,table_circling_l2} we depict numerical results in the energy norm and in the $L^\infty L^2$-norm for 7 levels of refinement. 
\begin{figure}[!htbp]
	\begin{subfigure}[t]{0.5\textwidth}
		\centering
		\begin{tikzpicture}[scale=0.73]
			\def\vara{0.6}
			\def\varb{0.1}
			\def\varc{0.01}
			\def\vard{0.001}
			\def\vare{0.0001}
			\begin{semilogyaxis}[xlabel={Refinement level $l_\ti=l_{\spa}$}, ylabel={$\triplenorm{u^e-u_h}\triplenorm_h$},legend style={ at={(0.5,1.02)}, anchor=south, legend columns =2}, legend cell align=left]
				\addplot table [x=energylist_x, y=-, col sep=comma] {Experiments/EnergyError_DiagKonvTable_Diss1_alpha_h_beta0_stabiwith_h_SpatialNormalVolume_userh_1_tmax_9_tstart_2_xmax_9_xstart_2_ks_1_kt_1_kgs_1_kgt_1_gp_0_tempquadorder_4_tempquadorderDiffu_4.csv};
				\addplot table [x=energylist_x, y=-, col sep=comma] {Experiments/EnergyError_DiagKonvTable_Diss1_alpha_h_beta0_stabiwith_h_SpatialNormalVolume_userh_1_tmax_9_tstart_2_xmax_9_xstart_2_ks_2_kt_2_kgs_2_kgt_2_gp_0_tempquadorder_8_tempquadorderDiffu_8.csv};
				\addplot table [x=energylist_x, y=-, col sep=comma] {Experiments/EnergyError_DiagKonvTable_Diss1_alpha_h_beta0_stabiwith_h_SpatialNormalVolume_userh_1_tmax_9_tstart_2_xmax_9_xstart_2_ks_3_kt_3_kgs_3_kgt_3_gp_0_tempquadorder_12_tempquadorderDiffu_12.csv};
				\addplot table [x=energylist_x, y=-, col sep=comma] {Experiments/EnergyError_DiagKonvTable_Diss1_alpha_h_beta0_stabiwith_h_SpatialNormalVolume_userh_1_tmax_9_tstart_2_xmax_9_xstart_2_ks_4_kt_4_kgs_4_kgt_4_gp_0_tempquadorder_16_tempquadorderDiffu_16.csv};
				\addplot table [x=energylist_x, y=-, col sep=comma] {Experiments/EnergyError_DiagKonvTable_Diss1_alpha_h_beta0_stabiwith_h_SpatialNormalVolume_userh_1_tmax_9_tstart_2_xmax_9_xstart_2_ks_5_kt_5_kgs_5_kgt_5_gp_0_tempquadorder_20_tempquadorderDiffu_20.csv};
				\addplot[dotted,line width=0.75pt] coordinates { 
					(0,\vara) (1,\vara*0.5) (2,\vara*0.25) (3,\vara*0.125) (4,\vara*0.0625)(5,\vara*0.03125) (6,\vara*0.015625)(7,\vara*0.0078125)
				};
				\addplot[dotted,line width=0.75pt] coordinates { 
			(0,\varb) (1,\varb*0.5*0.5) (2,\varb*0.25*0.25) (3,\varb*0.125*0.125) (4,\varb*0.0625*0.0625)(5,\varb*0.03125*0.03125) (6,\varb*0.015625*0.015625)(7,\varb*0.0078125*0.0078125)
		};
				\addplot[dotted,line width=0.75pt] coordinates { 
		(0,\varc) (1,\varc*0.5*0.5*0.5) (2,\varc*0.25*0.25*0.25) (3,\varc*0.125*0.125*0.125) (4,\varc*0.0625*0.0625*0.0625)(5,\varc*0.03125*0.03125*0.03125) (6,\varc*0.015625*0.015625*0.015625)(7,\varc*0.0078125*0.0078125*0.0078125)
	};
			\addplot[dotted,line width=0.75pt] coordinates { 
	(0,\vard) (1,\vard*0.5*0.5*0.5*0.5) (2,\vard*0.25*0.25*0.25*0.25) (3,\vard*0.125*0.125*0.125*0.125) (4,\vard*0.0625*0.0625*0.0625*0.0625)(5,\vard*0.03125*0.03125*0.03125*0.03125) (6,\vard*0.015625*0.015625*0.015625*0.015625)(7,\vard*0.0078125*0.0078125*0.0078125*0.0078125)
};
			\addplot[dotted,line width=0.75pt] coordinates { 
	(0,\vare) (1,\vare*0.5*0.5*0.5*0.5*0.5) (2,\vare*0.25*0.25*0.25*0.25*0.25) (3,\vare*0.125*0.125*0.125*0.125*0.125) (4,\vare*0.0625*0.0625*0.0625*0.0625*0.0625)(5,\vare*0.03125*0.03125*0.03125*0.03125*0.03125) 
};
				\legend{{$k=1$}, {$k=2$},{$k=3$},{$k=4$}, {$k=5$}, $\mathcal{O}(h^k)$}
			\end{semilogyaxis}
		\end{tikzpicture}
		\caption{Energy norm errors}
		\label{fig_circling_energy} 
	\end{subfigure}\hfill
	\begin{subfigure}[t]{0.5\textwidth}
		\centering
		\begin{tikzpicture}[scale=0.73]

		\def\varb{0.6}
		\def\varc{0.05}
		\def\vard{0.01}
		\def\vare{0.0001}
		\def\varf{0.00001}
		\begin{semilogyaxis}[xlabel={Refinement level $l_\ti=l_{\spa}$}, ylabel={$\norm{u^e-u_h}_{L^\infty L^2}$},legend style={ at={(0.5,1.02)}, anchor=south, legend columns =2}, legend cell align=left]
			\addplot table [x=lxt, y=-, col sep=comma] {Experiments/L2Error_DiagKonvTable_Diss1_alpha_h_beta0_stabiwith_h_SpatialNormalVolume_userh_1_tmax_9_tstart_2_xmax_9_xstart_2_ks_1_kt_1_kgs_1_kgt_1_gp_0_tempquadorder_4_tempquadorderDiffu_4.csv};
			\addplot table [x=lxt, y=-, col sep=comma] {Experiments/L2Error_DiagKonvTable_Diss1_alpha_h_beta0_stabiwith_h_SpatialNormalVolume_userh_1_tmax_9_tstart_2_xmax_9_xstart_2_ks_2_kt_2_kgs_2_kgt_2_gp_0_tempquadorder_8_tempquadorderDiffu_8.csv};
			\addplot table [x=lxt, y=-, col sep=comma] {Experiments/L2Error_DiagKonvTable_Diss1_alpha_h_beta0_stabiwith_h_SpatialNormalVolume_userh_1_tmax_9_tstart_2_xmax_9_xstart_2_ks_3_kt_3_kgs_3_kgt_3_gp_0_tempquadorder_12_tempquadorderDiffu_12.csv};
			\addplot table [x=lxt, y=-, col sep=comma] {Experiments/L2Error_DiagKonvTable_Diss1_alpha_h_beta0_stabiwith_h_SpatialNormalVolume_userh_1_tmax_9_tstart_2_xmax_9_xstart_2_ks_4_kt_4_kgs_4_kgt_4_gp_0_tempquadorder_16_tempquadorderDiffu_16.csv};
			\addplot table [x=lxt, y=-, col sep=comma] {Experiments/L2Error_DiagKonvTable_Diss1_alpha_h_beta0_stabiwith_h_SpatialNormalVolume_userh_1_tmax_9_tstart_2_xmax_9_xstart_2_ks_5_kt_5_kgs_5_kgt_5_gp_0_tempquadorder_20_tempquadorderDiffu_20.csv};
			\addplot[dashed,line width=0.75pt] coordinates { 
				(0,\varb) (1,\varb*0.5*0.5) (2,\varb*0.25*0.25) (3,\varb*0.125*0.125) (4,\varb*0.0625*0.0625)(5,\varb*0.03125*0.03125) (6,\varb*0.015625*0.015625)(7,\varb*0.0078125*0.0078125)
			};
			\addplot[dashed,line width=0.75pt] coordinates { 
				(0,\varc) (1,\varc*0.5*0.5*0.5) (2,\varc*0.25*0.25*0.25) (3,\varc*0.125*0.125*0.125) (4,\varc*0.0625*0.0625*0.0625)(5,\varc*0.03125*0.03125*0.03125) (6,\varc*0.015625*0.015625*0.015625)(7,\varc*0.0078125*0.0078125*0.0078125)
			};
			\addplot[dashed,line width=0.75pt] coordinates { 
				(0,\vard) (1,\vard*0.5*0.5*0.5*0.5) (2,\vard*0.25*0.25*0.25*0.25) (3,\vard*0.125*0.125*0.125*0.125) (4,\vard*0.0625*0.0625*0.0625*0.0625)(5,\vard*0.03125*0.03125*0.03125*0.03125) (6,\vard*0.015625*0.015625*0.015625*0.015625)(7,\vard*0.0078125*0.0078125*0.0078125*0.0078125)
			};
			\addplot[dashed,line width=0.75pt] coordinates { 
				(0,\vare) (1,\vare*0.5*0.5*0.5*0.5*0.5) (2,\vare*0.25*0.25*0.25*0.25*0.25) (3,\vare*0.125*0.125*0.125*0.125*0.125) (4,\vare*0.0625*0.0625*0.0625*0.0625*0.0625)(5,\vare*0.03125*0.03125*0.03125*0.03125*0.03125) (6,\vare*0.015625*0.015625*0.015625*0.015625*0.015625)
			};
					\addplot[dashed,line width=0.75pt] coordinates { 
			(0,\varf) (1,\varf*0.5*0.5*0.5*0.5*0.5*0.5) (2,\varf*0.25*0.25*0.25*0.25*0.25*0.25) (3,\varf*0.125*0.125*0.125*0.125*0.125*0.125) (4,\varf*0.0625*0.0625*0.0625*0.0625*0.0625*0.0625)(5,\varf*0.03125*0.03125*0.03125*0.03125*0.03125*0.03125) 
		};
			\legend{{$k=1$}, {$k=2$},{$k=3$},{$k=4$}, {$k=5$}, $\mathcal{O}(h^{k+1})$}
		\end{semilogyaxis}
	\end{tikzpicture}
		\caption{$L^\infty L^2$-norm errors}
		\label{fig_circling_l2} 
	\end{subfigure}
	\caption{Moving circle: convergence results in different norms.}
	\label{fig_circling}
\end{figure}
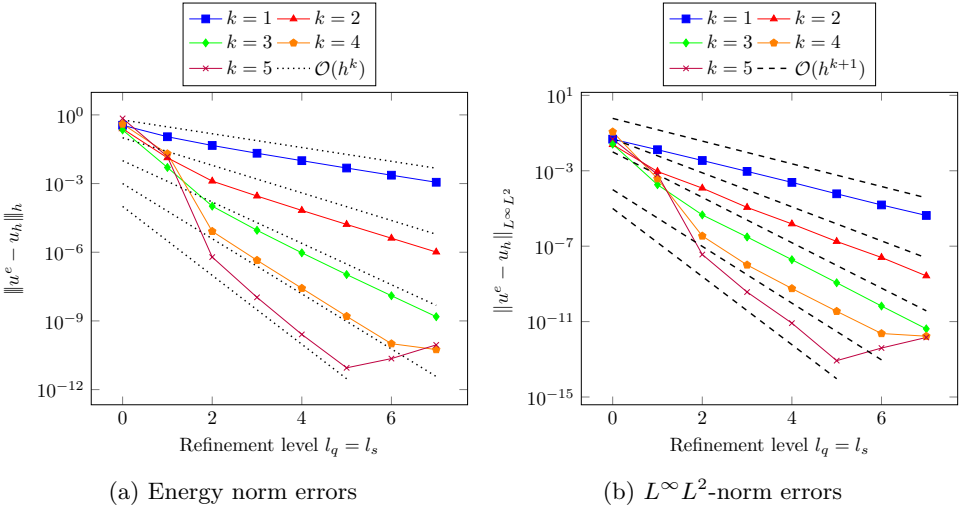
\begin{table}[!htbp]
	\centering
	\resizebox{\columnwidth}{!}{
		\begin{filecontents*}{Experiments/EnergyError_FullKonvTable_Diss1_alpha_h_beta0_stabiwith_h_SpatialNormalVolume_userh_1_tmax_9_tstart_2_xmax_9_xstart_2_ks_3_kt_3_kgs_3_kgt_3_gp_0_tempquadorder_12_tempquadorderDiffu_12.csv}
$l_\ti\downarrow\backslash l_\spa\rightarrow$,0,1,2,3,4,5,6,7,$\text{eoc}_{\ti}$
0,\cellcolor{gray!50} 2.20E-1,6.48E-2,2.45E-2,1.37E-2,1.24E-2,1.12E-2,5.97E-3,6.79E-3,-
1,2.72E-1,\cellcolor{gray!50} 5.05E-3,4.33E-4,2.73E-4,1.98E-4,1.52E-4,1.25E-4,1.10E-4,5.95
2,1.20E+0,1.40E-3,\cellcolor{gray!50} 1.03E-4,3.24E-5,2.07E-5,1.43E-5,1.07E-5,8.67E-6,3.67
3,1.54E+0,1.19E-3,7.53E-5,\cellcolor{gray!50} 9.09E-6,2.69E-6,1.64E-6,1.09E-6,7.88E-7,3.46
4,1.12E+0,1.11E-3,6.52E-5,6.32E-6,\cellcolor{gray!50} 9.42E-7,2.35E-7,1.31E-7,8.73E-8,3.17
5,8.10E-1,1.09E-3,6.18E-5,4.98E-6,6.23E-7,\cellcolor{gray!50} 1.06E-7,2.20E-8,1.08E-8,3.01
6,1.66E+0,1.08E-3,6.06E-5,4.37E-6,4.30E-7,6.95E-8,\cellcolor{gray!50} 1.26E-8,2.28E-9,2.24
7,2.06E+2,1.07E-3,6.02E-5,4.15E-6,3.38E-7,4.37E-8,8.24E-9,\cellcolor{gray!50} 1.53E-9,0.576
$\text{eoc}_{\spa}$,-,17.6,4.15,3.86,3.62,2.95,2.41,2.43, 
$\text{eoc}_{\spati}$ \cellcolor{gray!50},-,5.45,5.62,3.5,3.27,3.15,3.07,3.04, 
\end{filecontents*}
\csvreader[
		no head, 
		table head=\toprule,
		tabular={rrrrrrrrrr}, 
		late after line =,
		before line={\ifthenelse{\equal{\csvcoli}{$\text{eoc}_{\spa}$}\or\equal{\csvcolx}{-}}{\\\midrule}{\\}},
		before first line=,
		table foot=\\\bottomrule,
		]
		{Experiments/EnergyError_FullKonvTable_Diss1_alpha_h_beta0_stabiwith_h_SpatialNormalVolume_userh_1_tmax_9_tstart_2_xmax_9_xstart_2_ks_3_kt_3_kgs_3_kgt_3_gp_0_tempquadorder_12_tempquadorderDiffu_12.csv}{}{\csvcoli & \csvcolii & \csvcoliii & \csvcoliv & \csvcolv & \csvcolvi & \csvcolvii & \csvcolviii & \csvcolix  &\csvcolx}	
	}
	\caption{Moving circle: the errors $\triplenorm{u^e-u_h}\triplenorm_h$ for $k=3$.}
	\label{table_circling_energy}
\end{table}
\begin{table}[!htbp]
	\centering
	\resizebox{\columnwidth}{!}{
		\begin{filecontents*}{Experiments/L2Error_FullKonvTable_Diss1_alpha_h_beta0_stabiwith_h_SpatialNormalVolume_userh_1_tmax_9_tstart_2_xmax_9_xstart_2_ks_3_kt_3_kgs_3_kgt_3_gp_0_tempquadorder_12_tempquadorderDiffu_12.csv}
$l_\ti\downarrow\backslash l_\spa\rightarrow$,0,1,2,3,4,5,6,7,$\text{eoc}_{\ti}$
0,\cellcolor{gray!50} 2.57E-2,5.93E-3,5.81E-4,4.28E-4,3.68E-4,3.53E-4,3.46E-4,3.42E-4,-
1,3.20E-2,\cellcolor{gray!50} 1.83E-4,2.30E-5,1.62E-5,1.44E-5,1.38E-5,1.37E-5,1.37E-5,4.64
2,2.95E-1,1.64E-4,\cellcolor{gray!50} 4.56E-6,1.01E-6,6.43E-7,5.88E-7,5.84E-7,5.86E-7,4.55
3,2.27E-1,1.62E-4,4.46E-6,\cellcolor{gray!50} 3.08E-7,5.89E-8,3.38E-8,2.28E-8,2.12E-8,4.79
4,7.46E-1,1.62E-4,4.45E-6,3.12E-7,\cellcolor{gray!50} 1.92E-8,3.65E-9,2.25E-9,1.16E-9,4.19
5,4.77E-1,1.62E-4,4.45E-6,3.11E-7,1.92E-8,\cellcolor{gray!50} 1.15E-9,2.12E-10,1.44E-10,3.01
6,1.04E+0,1.62E-4,4.45E-6,3.16E-7,1.93E-8,1.15E-9,\cellcolor{gray!50} 6.68E-11,1.23E-11,3.55
7,9.49E+1,1.62E-4,4.45E-6,3.16E-7,1.93E-8,1.16E-9,6.67E-11,\cellcolor{gray!50} 4.13E-12,1.57
$\text{eoc}_{\spa}$,-,19.2,5.19,3.82,4.03,4.06,4.12,4.01, 
$\text{eoc}_{\spati}$ \cellcolor{gray!50},-,7.13,5.33,3.89,4.0,4.06,4.11,4.02, 
\end{filecontents*}
\csvreader[
		no head, 
		table head=\toprule,
		tabular={rrrrrrrrrr}, 
		late after line =,
		before line={\ifthenelse{\equal{\csvcoli}{$\text{eoc}_{\spa}$}\or\equal{\csvcolx}{-}}{\\\midrule}{\\}},
		before first line=,
		table foot=\\\bottomrule,
		]
		{Experiments/L2Error_FullKonvTable_Diss1_alpha_h_beta0_stabiwith_h_SpatialNormalVolume_userh_1_tmax_9_tstart_2_xmax_9_xstart_2_ks_3_kt_3_kgs_3_kgt_3_gp_0_tempquadorder_12_tempquadorderDiffu_12.csv}{}{\csvcoli & \csvcolii & \csvcoliii & \csvcoliv & \csvcolv & \csvcolvi & \csvcolvii & \csvcolviii & \csvcolix  &\csvcolx}	
	}
	\caption{Moving circle: the errors $\norm{u^e-u_h}_{L^\infty L^2}$ for $k=3$.}
	\label{table_circling_l2}
\end{table}
 For $k\in \{1,2,3,4,5\}$ we observe in \Cref{fig_circling} optimal order convergence for both the energy norm error and the $L^\infty L^2$-error until error values small than $10^{-10}$ are reached. 
 For $k=4$, the error in the $L^\infty L^2$-norm shows a suboptimal rate of convergence with order approximately $4$ (instead of the optimal order 5).
 In further experiments (results not shown here) we observed that an optimal order convergence can be obtained by increasing the order of the geometry approximation, i.e. taking $k_{g,s}=k_{g,t}=5$. Taking only a higher order approximation of the geometric quantity $\alpha_h$ is not sufficient here, cf. \Cref{section_alpha_exps}. 
To further study the higher order methods, we present more detailed space-time convergence  results for the case $k=3$ in \Cref{table_circling_l2,table_circling_energy}. We observe optimal order convergence on the space-time diagonal.
In these tables the first column gives somewhat unstable results. This is caused by the mesh deformation. When deforming a too coarse mesh, the distance between an undeformed and the corresponding deformed tetrahedron might become so large that undesirable effects occur. For example, faces of the deformed tetrahedra might touch each other. To avoid such effects, one uses a strategy in which certain tetrahedra are not allowed to deform, leading to an only piecewise linear approximation of the surface for the under-resolved case. 

\subsubsection{Twisting torus}\label{section_twisting_torus}
\begin{figure}[ht!]
	\begin{minipage}{\textwidth}
		\vspace{-1cm}
		\includegraphics[width=\textwidth]{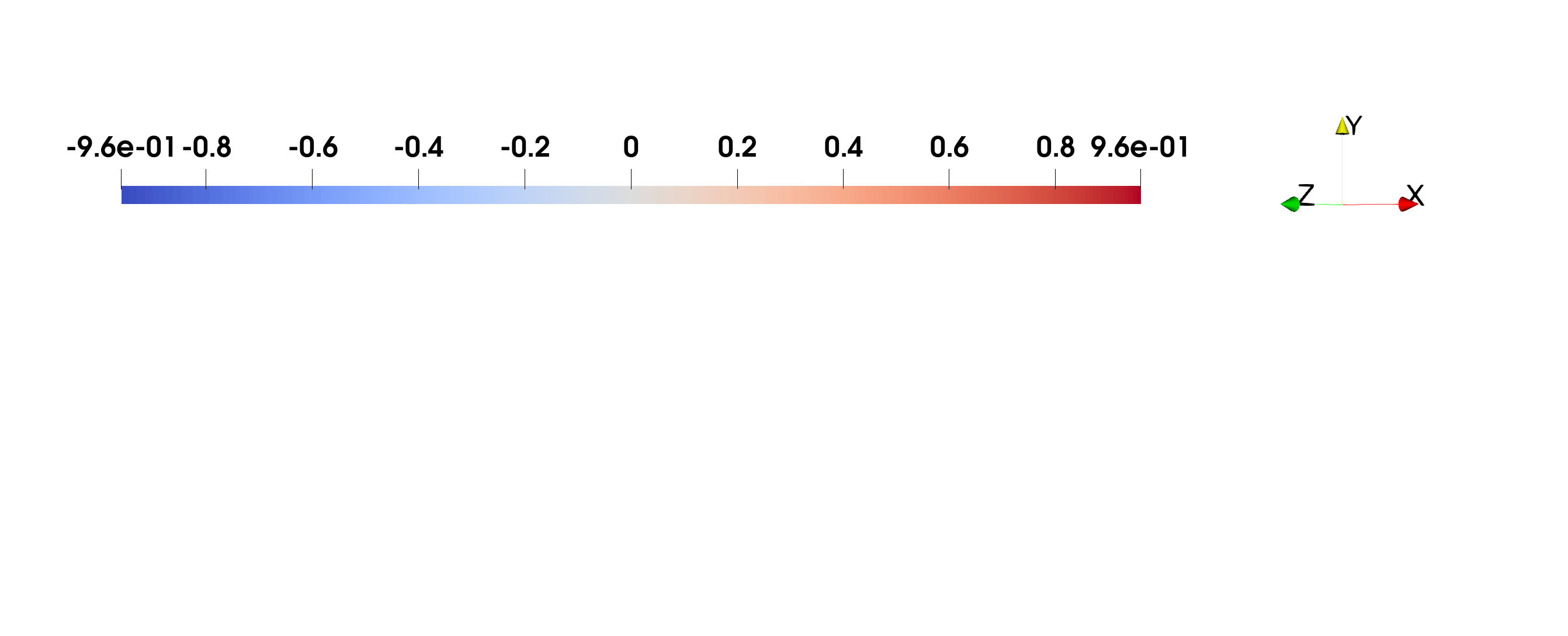}    
		\vspace{-3.5cm}
	\end{minipage}
	\begin{minipage}{0.33\textwidth}
		\includegraphics[width=\textwidth]{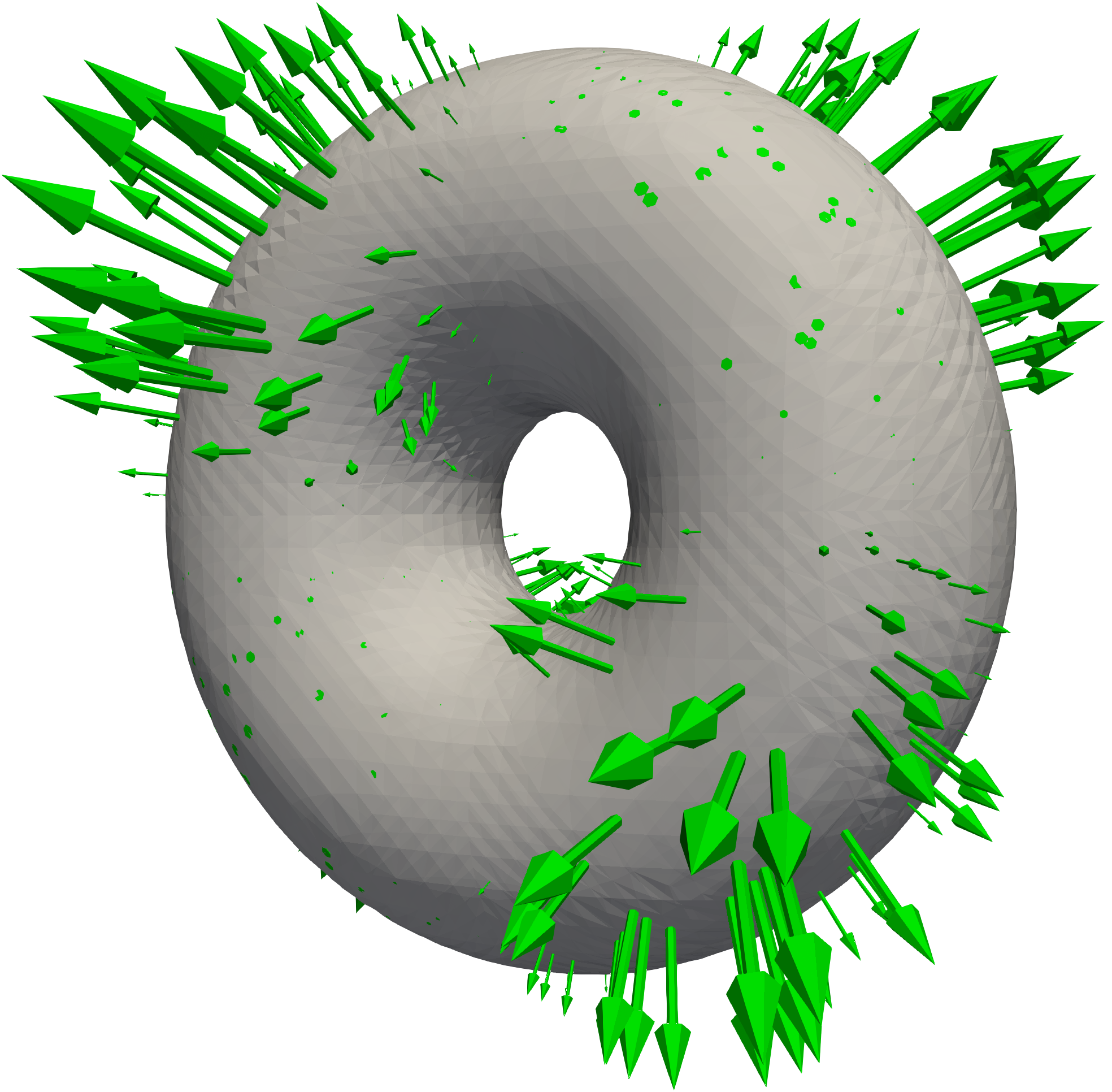}    
	\end{minipage}\hfill
	\begin{minipage}{0.33\textwidth}
		\includegraphics[width=\textwidth]{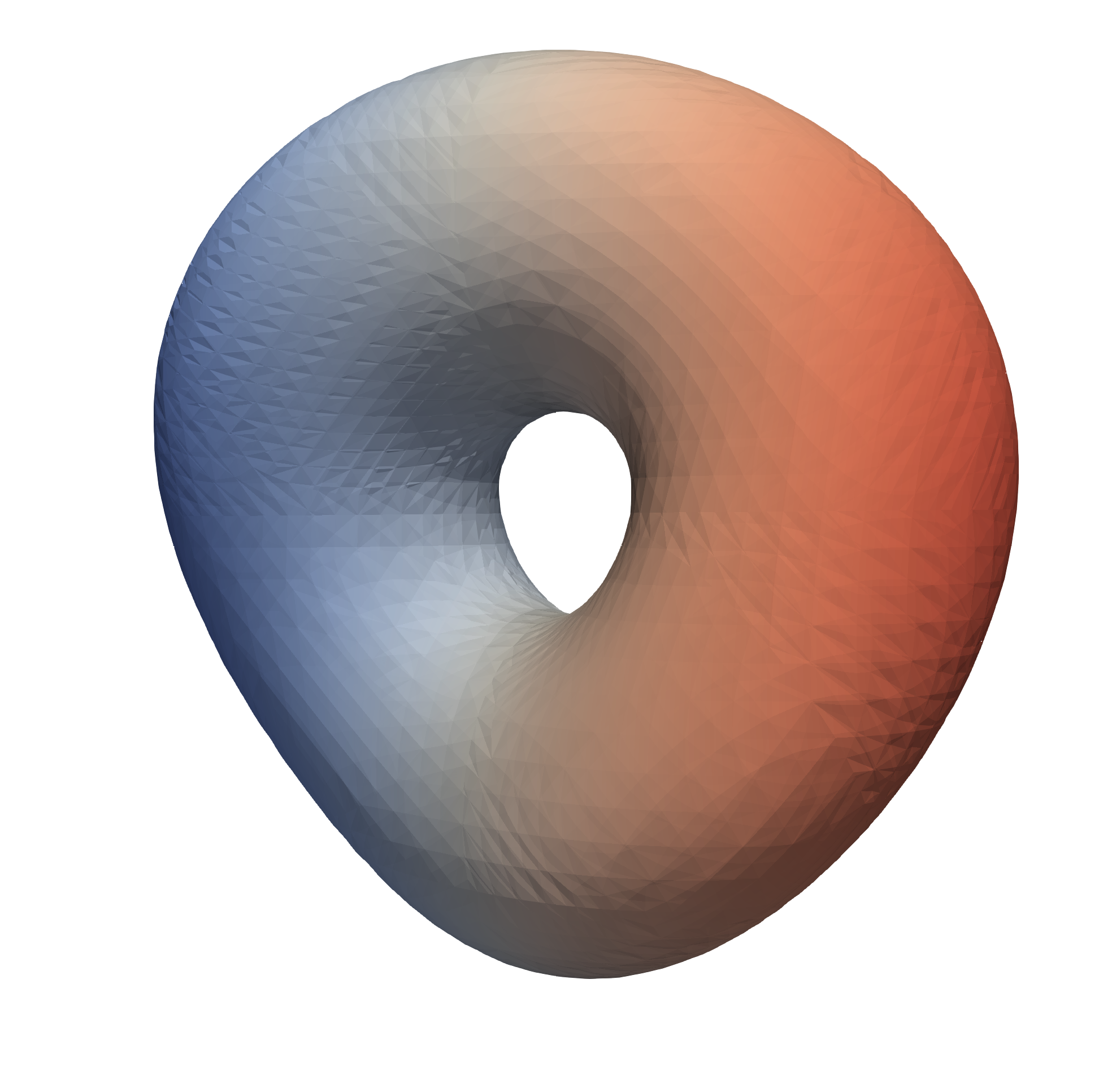}     
	\end{minipage}\hfill
	\begin{minipage}{0.33\textwidth}
		\includegraphics[width=\textwidth]{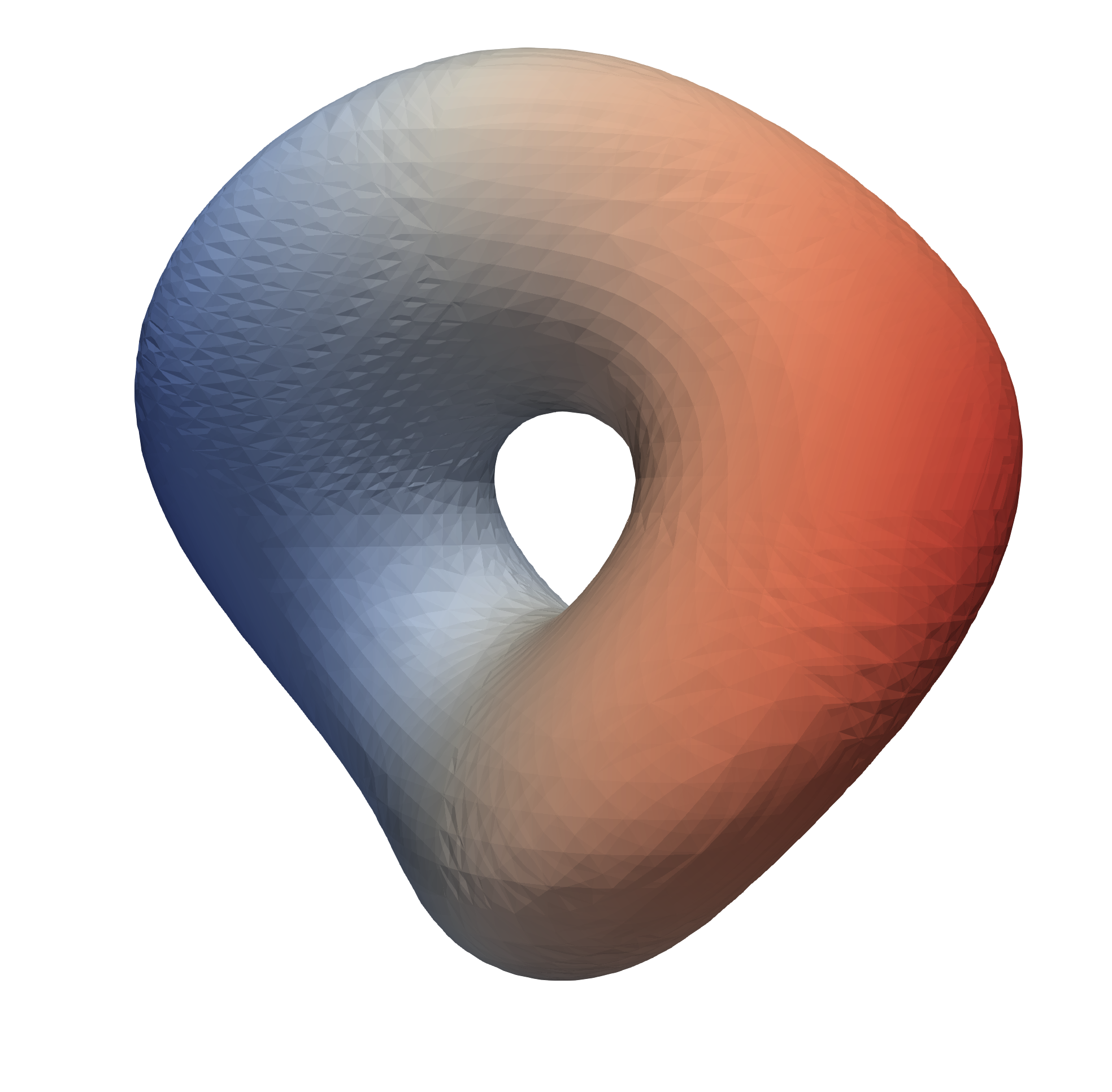}     
	\end{minipage}
	\begin{minipage}{0.33\textwidth}
		\includegraphics[width=\textwidth]{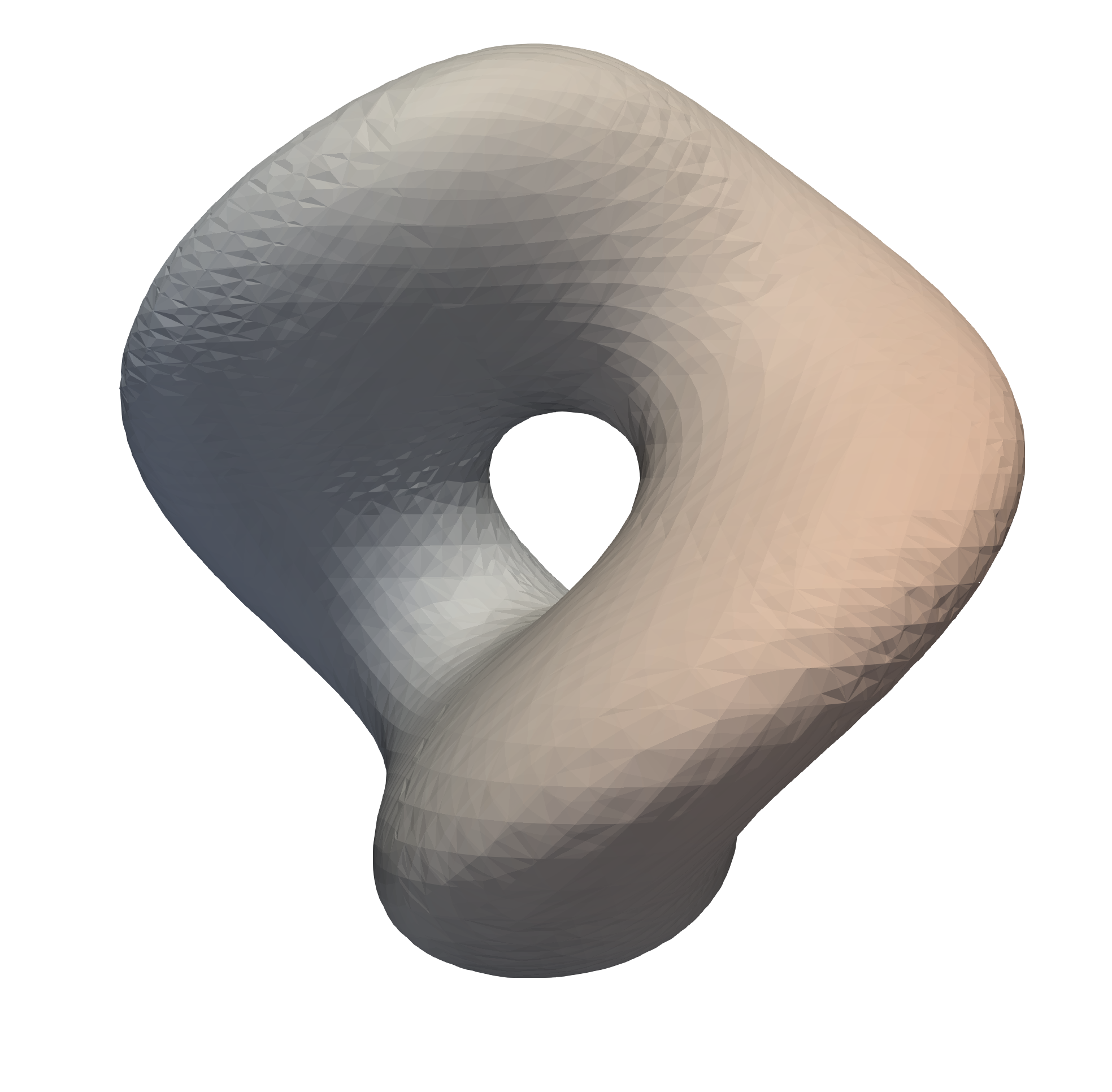}     
	\end{minipage}\hfill
	\begin{minipage}{0.33\textwidth}
		\includegraphics[width=\textwidth]{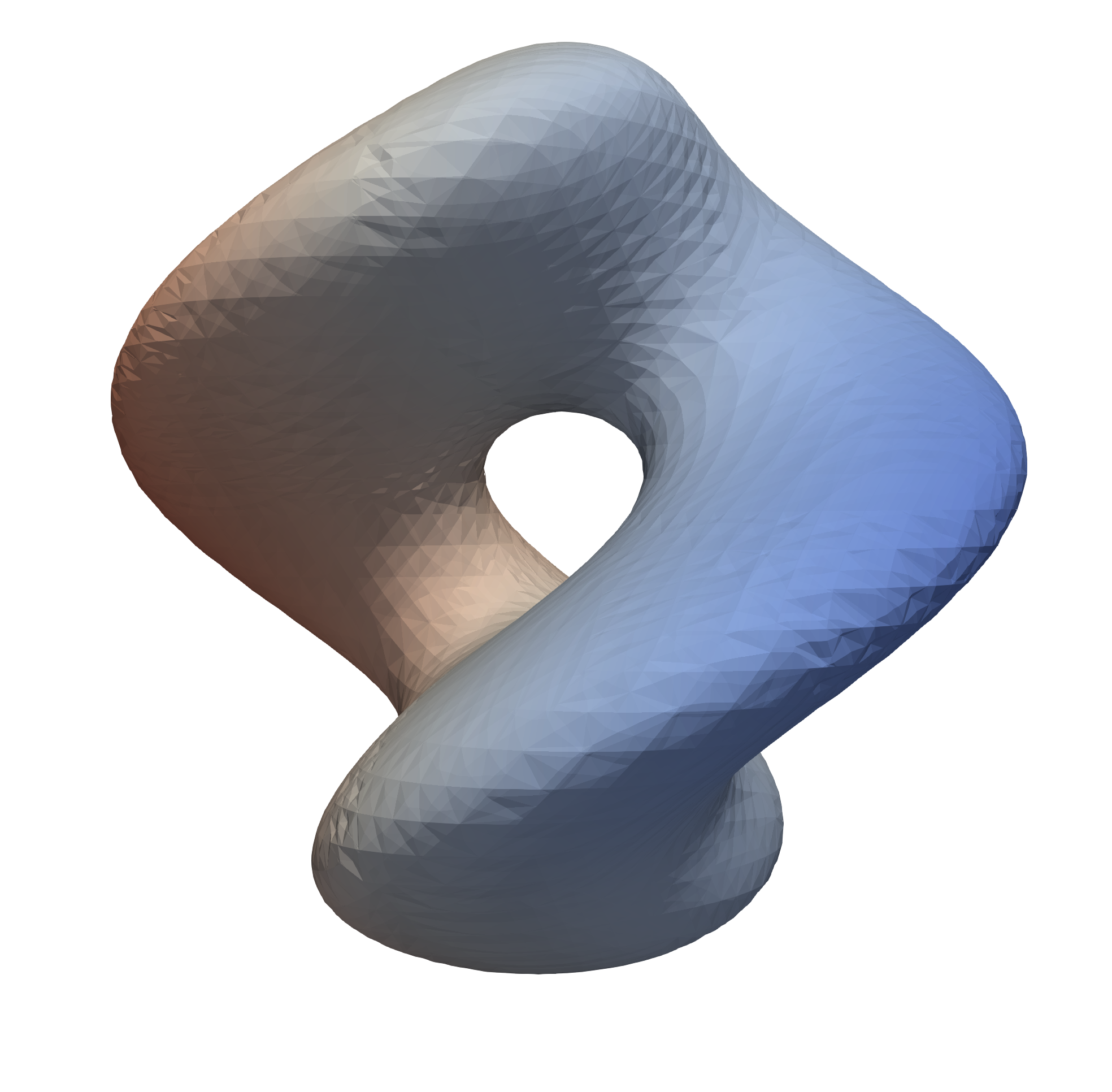}     
	\end{minipage}\hfill
	\begin{minipage}{0.33\textwidth}
		\includegraphics[width=\textwidth]{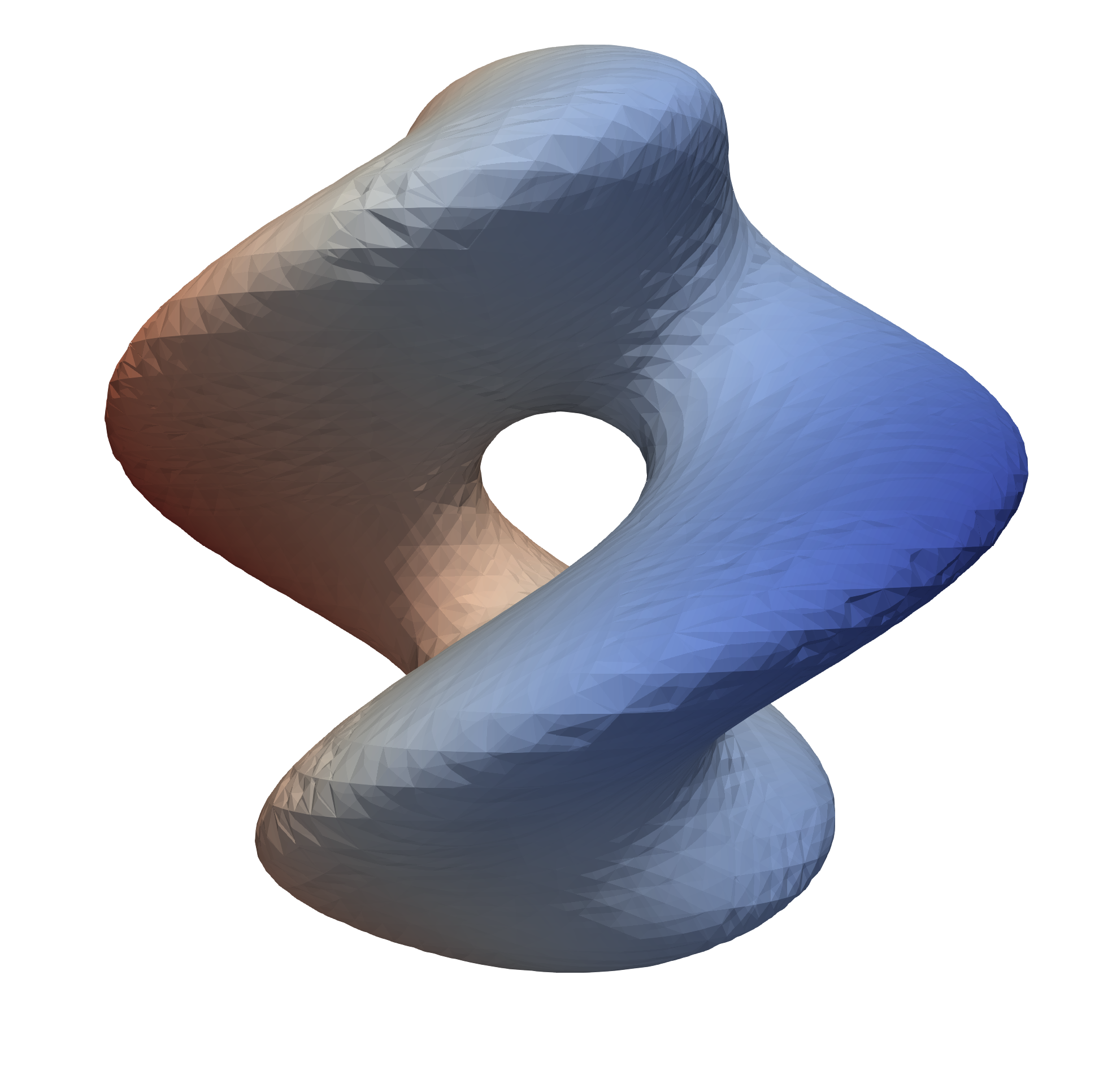}     
	\end{minipage}
	\begin{minipage}{0.33\textwidth}
		\includegraphics[width=\textwidth]{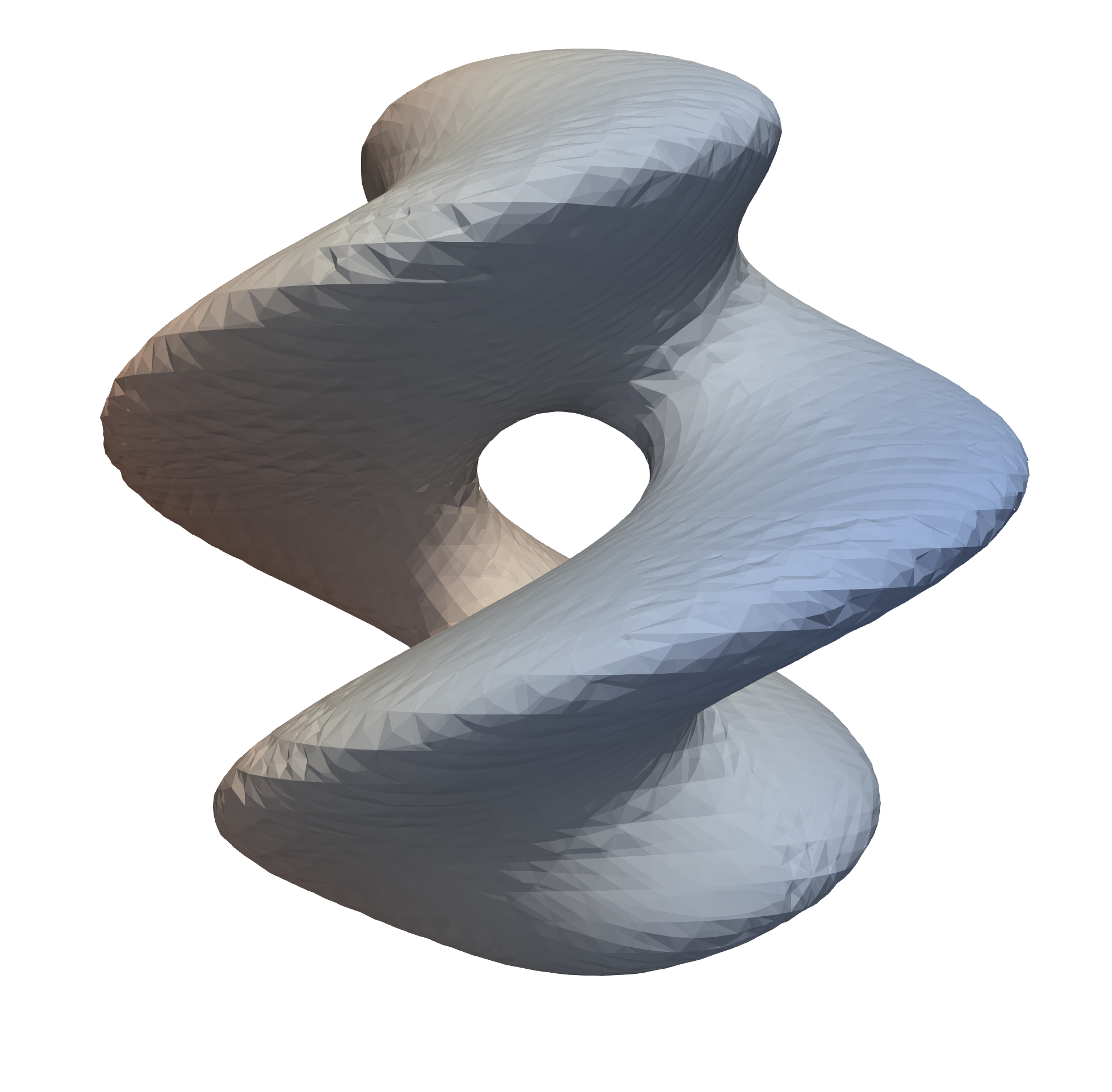}     
	\end{minipage}\hfill
	\begin{minipage}{0.33\textwidth}
		\includegraphics[width=\textwidth]{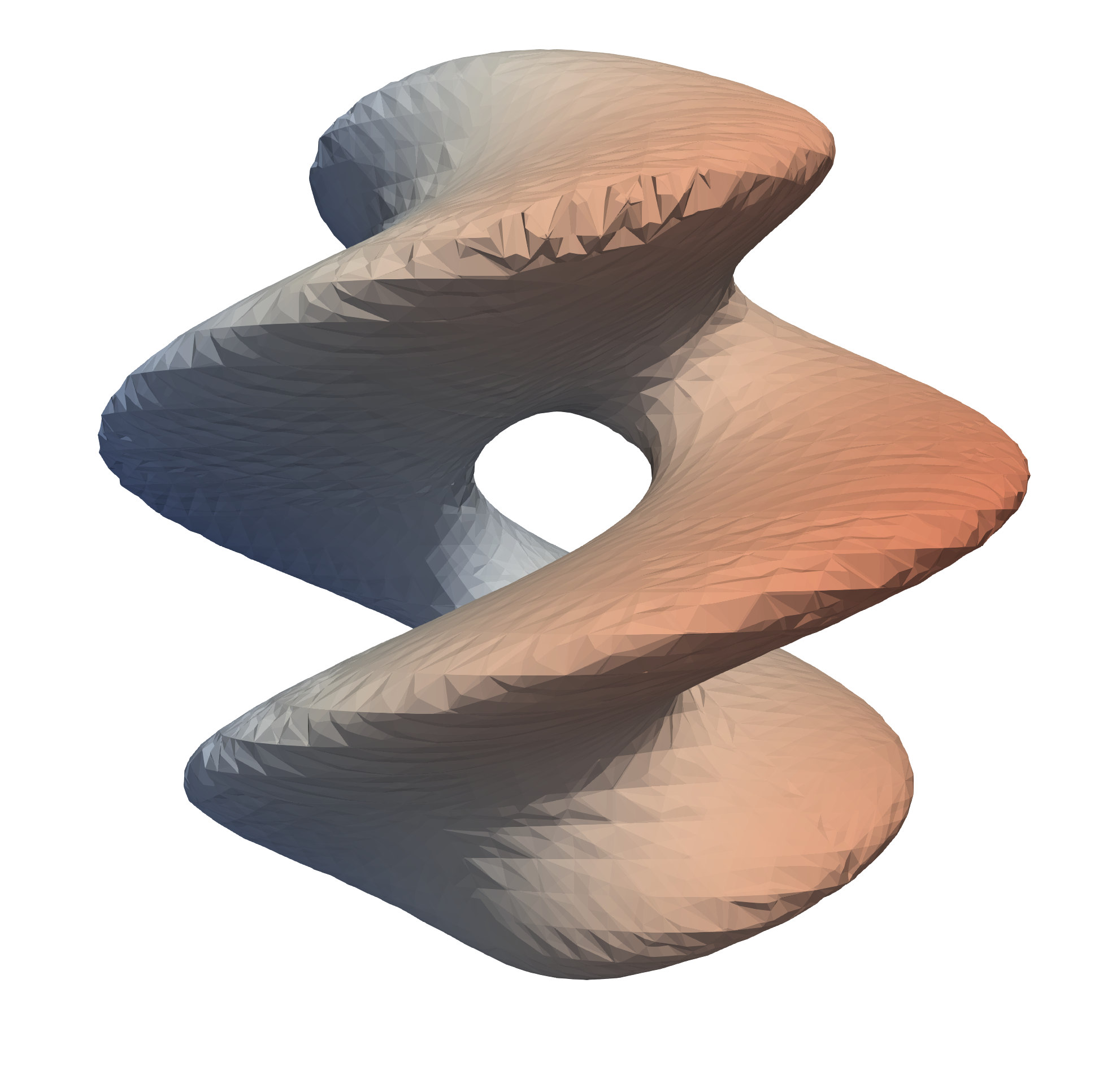}     
	\end{minipage}\hfill
	\begin{minipage}{0.33\textwidth}
		\includegraphics[width=\textwidth]{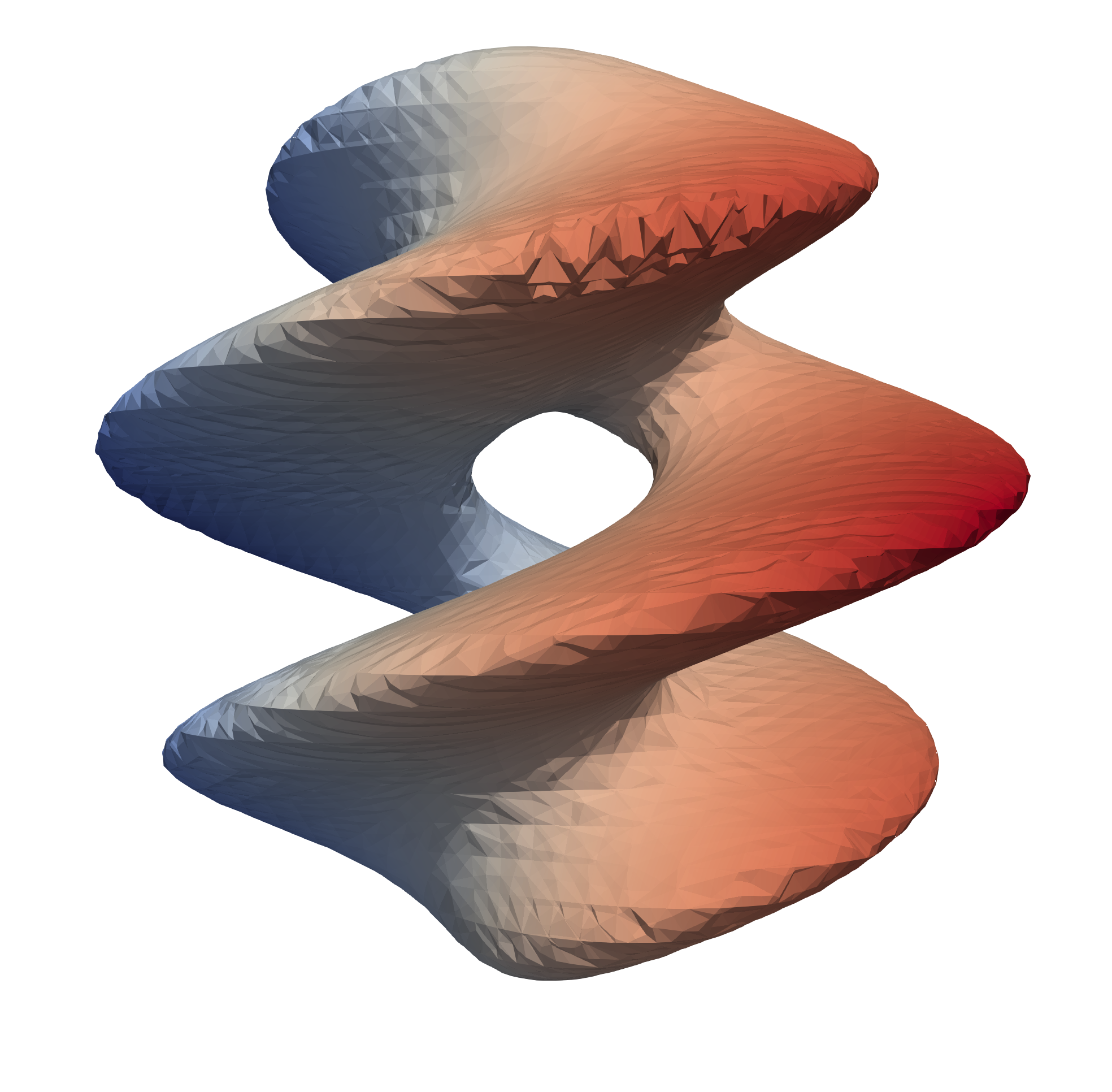}     
	\end{minipage}
	\caption{A sketch of the evolution of the twisted torus for $k=1$. The discrete solution $u_h$ at times $t=\frac{iT}{8}$, $i=0,\dots,8$, is depicted in the images from top left to bottom right. Here, we have $\Delta t=5\cdot 2^{-5}$, $h=2^{-5}$ and $T=5$. The green arrows represent the velocity field $\bw$. A full animation can be found \href{https://doi.org/10.5281/zenodo.7385373}{\textcolor{blue}{here}}.}
	\label{figure_torus_twist}
\end{figure}

We consider a torus, formed with the $y$-axis as axis of rotation. The torus is twisted over time. We take the background domain $\Omega=(-1,1)^3$ and the time interval $[0,T]$, $T\in \{\frac52, 5\}$. On $Q=\Omega\times [0,5]$ the evolution of the surface $\Gamma(t)$ is described by the zero level of 
\begin{align*}
	\phi(\bx,t)&=\left(\omega_1^2(\bx,t)+\omega_2^2(\bx,t)+\omega^2_3(\bx,t)+\frac{9}{50}\right)^2-\frac{121}{100}\left(\omega_1^2(\bx,t)+\omega^2_3(\bx,t)\right),\\
	(\bx,t)&=(x,y,z,t)\in Q,
\end{align*}
where $\omega_1,\omega_2,\omega_2:Q\rightarrow \R$ are defined as
\begin{equation*}
	\begin{pmatrix}\omega_1\\\omega_2\\\omega_3\end{pmatrix}(\bx,t)\coloneqq \begin{pmatrix}\cos(ty) & 0 &-\sin(ty)\\\sin(ty)& 0& \cos(ty)\\ 0&1&0\end{pmatrix}\begin{pmatrix}x\\y\\z\end{pmatrix}.
\end{equation*}
On $Q$ we use a normal velocity field that transports $\Gamma(t)$ by 
\begin{equation*}
	\bw=-\ddt{\phi}{t}\frac{\nabla\phi}{\norm{\nabla \phi}^2}. 
\end{equation*}
We prescribe the exact solution as 
\begin{equation*}
	u(\bx,t)= u^e(\bx,t)=x\sin\left(\frac{\pi t}{2}\right), \quad(\bx,t)=(x,y,z,t)\in Q.
\end{equation*}
We take the initial mesh size as $h_{\mathrm{init}}=2^{-1}$ and the initial time step size as $\Delta t_{\mathrm{init}}=2^{-3}$. In \Cref{figure_torus_twist} we give a sketch of the geometry and the discrete solution on the torus. In \Cref{fig_twisting,twisting_t25} we present errors in the surface norm $\triplenorm \cdot\triplenorm_{h,*}$ and in the $L^2L^\infty$-norm.

\begin{figure}[!htbp]
	\begin{subfigure}[t]{0.5\textwidth}
		\centering
		\begin{tikzpicture}[scale=0.73]
			\def\vara{5}
			\def\varb{10}
			\def\varc{10}
			\def\vard{20}
			\begin{semilogyaxis}[xlabel={Refinement level $l_\ti=l_{\spa}$}, ylabel={$\triplenorm{u^e-u_h}\triplenorm_{h,*}$},legend style={ at={(0.5,1.02)}, anchor=south, legend columns =2}, legend cell align=left]
				\addplot table [x=energylistwonormal_x, y=-, col sep=comma] {Experiments/TwistingTorus_T25/EnergyErrorWoNormal_DiagKonvTable_Diss3_alpha_h_beta0_stabiwith_h_SpatialNormalVolume_userh_1_tmax_9_tstart_3_xmax_7_xstart_1_ks_1_kt_1_kgs_1_kgt_1_gp_0_tempquadorder_4_tempquadorderDiffu_4.csv};
				\addplot table [x=energylistwonormal_x, y=-, col sep=comma] {Experiments/TwistingTorus_T25/EnergyErrorWoNormal_DiagKonvTable_Diss3_alpha_h_beta0_stabiwith_h_SpatialNormalVolume_userh_1_tmax_8_tstart_3_xmax_6_xstart_1_ks_2_kt_2_kgs_2_kgt_2_gp_0_tempquadorder_8_tempquadorderDiffu_8.csv};
				\addplot table [x=energylistwonormal_x, y=-, col sep=comma] {Experiments/TwistingTorus_T25/EnergyErrorWoNormal_DiagKonvTable_Diss3_alpha_h_beta0_stabiwith_h_SpatialNormalVolume_userh_1_tmax_7_tstart_3_xmax_5_xstart_1_ks_3_kt_3_kgs_3_kgt_3_gp_0_tempquadorder_12_tempquadorderDiffu_12.csv};
				\addplot table [x=energylistwonormal_x, y=-, col sep=comma] {Experiments/TwistingTorus_T25/EnergyErrorWoNormal_DiagKonvTable_Diss3_alpha_h_beta0_stabiwith_h_SpatialNormalVolume_userh_1_tmax_7_tstart_3_xmax_5_xstart_1_ks_4_kt_4_kgs_4_kgt_4_gp_0_tempquadorder_16_tempquadorderDiffu_16.csv};
				\addplot[dotted,line width=0.75pt] coordinates { 
					(1,\vara*0.5) (2,\vara*0.25) (3,\vara*0.125) (4,\vara*0.0625)(5,\vara*0.03125) (6,\vara*0.015625)};
				\addplot[dotted,line width=0.75pt] coordinates { 
					(1,\varb*0.5*0.5) (2,\varb*0.25*0.25) (3,\varb*0.125*0.125) (4,\varb*0.0625*0.0625)(5,\varb*0.03125*0.03125) 	};
								\addplot[dotted,line width=0.75pt] coordinates { 
									(0,\varc) (1,\varc*0.5*0.5*0.5) (2,\varc*0.25*0.25*0.25) (3,\varc*0.125*0.125*0.125) (4,\varc*0.0625*0.0625*0.0625)};
								\addplot[dotted,line width=0.75pt] coordinates { 
									(1,\vard*0.5*0.5*0.5*0.5) (2,\vard*0.25*0.25*0.25*0.25) (3,\vard*0.125*0.125*0.125*0.125) (4,\vard*0.0625*0.0625*0.0625*0.0625)	};
				\legend{{$k=1$},{$k=2$},{$k=3$},{$k=4$}, $\mathcal{O}(h^k)$}
			\end{semilogyaxis}
		\end{tikzpicture}
		\caption{Surface energy norm errors}
		\label{twisting_t25_energy} 
	\end{subfigure}\hfill
	\begin{subfigure}[t]{0.5\textwidth}
		\centering
		\begin{tikzpicture}[scale=0.73]	
			\def\varb{10}
			\def\varc{2}
			\def\vard{3}
			\def\vare{3}
			\begin{semilogyaxis}[xlabel={Refinement level $l_\ti=l_{\spa}$}, ylabel={$\norm{u^e-u_h}_{L^\infty L^2}$},legend style={ at={(0.5,1.02)}, anchor=south, legend columns =2}, legend cell align=left]
				\addplot table [x=lxt, y=-, col sep=comma] {Experiments/TwistingTorus_T25/L2Error_DiagKonvTable_Diss3_alpha_h_beta0_stabiwith_h_SpatialNormalVolume_userh_1_tmax_9_tstart_3_xmax_7_xstart_1_ks_1_kt_1_kgs_1_kgt_1_gp_0_tempquadorder_4_tempquadorderDiffu_4.csv};
				\addplot table [x=lxt, y=-, col sep=comma] {Experiments/TwistingTorus_T25/L2Error_DiagKonvTable_Diss3_alpha_h_beta0_stabiwith_h_SpatialNormalVolume_userh_1_tmax_8_tstart_3_xmax_6_xstart_1_ks_2_kt_2_kgs_2_kgt_2_gp_0_tempquadorder_8_tempquadorderDiffu_8.csv};
				\addplot table [x=lxt, y=-, col sep=comma] {Experiments/TwistingTorus_T25/L2Error_DiagKonvTable_Diss3_alpha_h_beta0_stabiwith_h_SpatialNormalVolume_userh_1_tmax_7_tstart_3_xmax_5_xstart_1_ks_3_kt_3_kgs_3_kgt_3_gp_0_tempquadorder_12_tempquadorderDiffu_12.csv};
				\addplot table [x=lxt, y=-, col sep=comma] {Experiments/TwistingTorus_T25/L2Error_DiagKonvTable_Diss3_alpha_h_beta0_stabiwith_h_SpatialNormalVolume_userh_1_tmax_7_tstart_3_xmax_5_xstart_1_ks_4_kt_4_kgs_4_kgt_4_gp_0_tempquadorder_16_tempquadorderDiffu_16.csv};
				\addplot[dashed,line width=0.75pt] coordinates { 
					(1,\varb*0.5*0.5) (2,\varb*0.25*0.25) (3,\varb*0.125*0.125) (4,\varb*0.0625*0.0625)(5,\varb*0.03125*0.03125)(6,\varb*0.015625*0.015625)	};
				\addplot[dashed,line width=0.75pt] coordinates { 
					(1,\varc*0.5*0.5*0.5) (2,\varc*0.25*0.25*0.25) (3,\varc*0.125*0.125*0.125) (4,\varc*0.0625*0.0625*0.0625) (5,\varc*0.03125*0.03125*0.03125)	};
								\addplot[dashed,line width=0.75pt] coordinates { 
								(1,\vard*0.5*0.5*0.5*0.5) (2,\vard*0.25*0.25*0.25*0.25) (3,\vard*0.125*0.125*0.125*0.125) (4,\vard*0.0625*0.0625*0.0625*0.0625)};
								\addplot[dashed,line width=0.75pt] coordinates { 
									 (1,\vare*0.5*0.5*0.5*0.5*0.5) (2,\vare*0.25*0.25*0.25*0.25*0.25) (3,\vare*0.125*0.125*0.125*0.125*0.125)  (4,\vare*0.0625*0.0625*0.0625*0.0625*0.0625)	};
				\legend{{$k=1$},{$k=2$},{$k=3$},{$k=4$}, $\mathcal{O}(h^{k+1})$}
			\end{semilogyaxis}
		\end{tikzpicture}
		\caption{$L^\infty L^2$-norm errors}
		\label{twisting_t25_l2} 
	\end{subfigure}
	\caption{Twisting torus: convergence results for $T=\frac52$ in different norms.}
	\label{twisting_t25}
\end{figure}
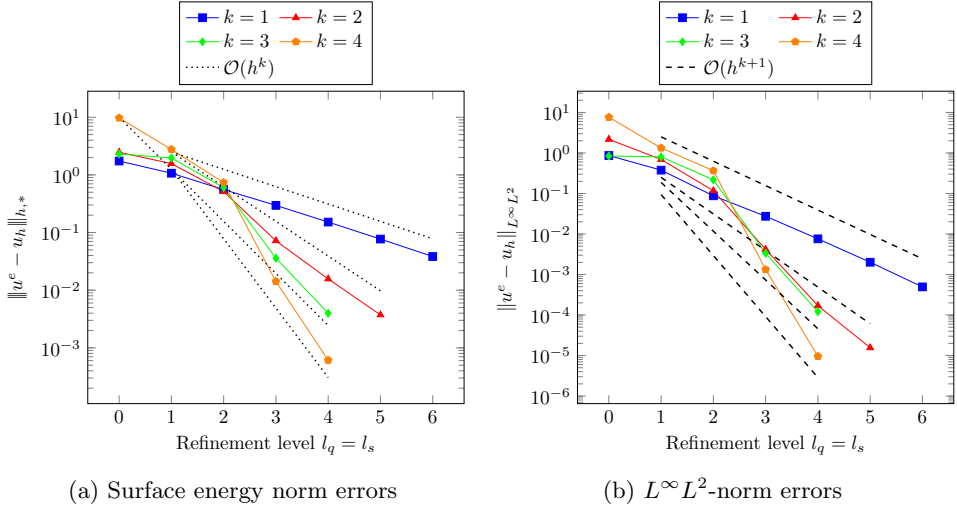
\begin{figure}[!htbp]
	\begin{subfigure}[t]{0.5\textwidth}
		\centering
		\begin{tikzpicture}[scale=0.73]
			\def\vara{5}
			\def\varb{10}
			\def\varc{10}
			\def\vard{10}
			\begin{semilogyaxis}[xlabel={Refinement level $l_\ti=l_{\spa}$}, ylabel={$\triplenorm{u^e-u_h}\triplenorm_{h,*}$},legend style={ at={(0.5,1.02)}, anchor=south, legend columns =2}, legend cell align=left]
			\addplot+[mark repeat =2, mark phase =1] table [x=energylistwonormal_x, y=-, col sep=comma] {Experiments/EnergyErrorWoNormal_DiagKonvTable_Diss3_alpha_h_beta0_stabiwith_h_SpatialNormalVolume_userh_1_tmax_8_tstart_3_xmax_6_xstart_1_ks_1_kt_1_kgs_2_kgt_2_gp_0_tempquadorder_4_tempquadorderDiffu_4.csv};
				\addplot +[mark repeat =2, mark phase =2]table [x=energylistwonormal_x, y=-, col sep=comma] {Experiments/EnergyErrorWoNormal_DiagKonvTable_Diss3_alpha_h_beta0_stabiwith_h_SpatialNormalVolume_userh_1_tmax_8_tstart_3_xmax_6_xstart_1_ks_2_kt_2_kgs_3_kgt_3_gp_0_tempquadorder_8_tempquadorderDiffu_8.csv};
				\addplot +[mark repeat =2, mark phase =1]table [x=energylistwonormal_x, y=-, col sep=comma] {Experiments/EnergyErrorWoNormal_DiagKonvTable_Diss3_alpha_h_beta0_stabiwith_h_SpatialNormalVolume_userh_1_tmax_7_tstart_3_xmax_5_xstart_1_ks_3_kt_3_kgs_4_kgt_4_gp_0_tempquadorder_12_tempquadorderDiffu_12.csv};
				\addplot[dotted,line width=0.75pt] coordinates { 
					(1,\vara*0.5) (2,\vara*0.25) (3,\vara*0.125) (4,\vara*0.0625)(5,\vara*0.03125) };
				\addplot[dashed,line width=0.75pt] coordinates { 
					(1,\varb*0.5*0.5) (2,\varb*0.25*0.25) (3,\varb*0.125*0.125) (4,\varb*0.0625*0.0625)	(5,\varb*0.03125*0.03125) };
			\legend{{$k=1$}, {$k=2$},{$k=3$}, $\mathcal{O}(h)$, $\mathcal{O}(h^2)$}
			\end{semilogyaxis}
		\end{tikzpicture}
		\caption{Surface energy norm errors}
		\label{fig_twisting_energy} 
	\end{subfigure}\hfill
	\begin{subfigure}[t]{0.5\textwidth}
		\centering
		\begin{tikzpicture}[scale=0.73]	
			\def\varb{80}
			\def\varc{90}
			\def\vard{1}
			\def\vare{1}
			\begin{semilogyaxis}[xlabel={Refinement level $l_\ti=l_{\spa}$}, ylabel={$\norm{u^e-u_h}_{L^\infty L^2}$},legend style={ at={(0.5,1.02)}, anchor=south, legend columns =2}, legend cell align=left]
			\addplot +[mark repeat =2, mark phase =1]table [x=lxt, y=-, col sep=comma] {Experiments/L2Error_DiagKonvTable_Diss3_alpha_h_beta0_stabiwith_h_SpatialNormalVolume_userh_1_tmax_8_tstart_3_xmax_6_xstart_1_ks_1_kt_1_kgs_2_kgt_2_gp_0_tempquadorder_4_tempquadorderDiffu_4.csv};
			\addplot +[mark repeat =2, mark phase =2]table [x=lxt, y=-, col sep=comma] {Experiments/L2Error_DiagKonvTable_Diss3_alpha_h_beta0_stabiwith_h_SpatialNormalVolume_userh_1_tmax_8_tstart_3_xmax_6_xstart_1_ks_2_kt_2_kgs_3_kgt_3_gp_0_tempquadorder_8_tempquadorderDiffu_8.csv};
				\addplot +[mark repeat =2, mark phase =1]table [x=lxt, y=-, col sep=comma] {Experiments/L2Error_DiagKonvTable_Diss3_alpha_h_beta0_stabiwith_h_SpatialNormalVolume_userh_1_tmax_7_tstart_3_xmax_5_xstart_1_ks_3_kt_3_kgs_4_kgt_4_gp_0_tempquadorder_12_tempquadorderDiffu_12.csv};
				\addplot[dotted,line width=0.75pt] coordinates { 
					(1,\varb*0.5*0.5) (2,\varb*0.25*0.25) (3,\varb*0.125*0.125) (4,\varb*0.0625*0.0625)(5,\varb*0.03125*0.03125)	};
				\addplot[dashed,line width=0.75pt] coordinates { 
					(1,\varc*0.5*0.5*0.5) (2,\varc*0.25*0.25*0.25) (3,\varc*0.125*0.125*0.125) (4,\varc*0.0625*0.0625*0.0625) (5,\varc*0.03125*0.03125*0.03125)	};
\legend{{$k=1$}, {$k=2$},{$k=3$}, $\mathcal{O}(h^{2})$,$\mathcal{O}(h^{3})$}
			\end{semilogyaxis}
		\end{tikzpicture}
		\caption{$L^\infty L^2$-norm errors}
		\label{fig_twisting_l2} 
	\end{subfigure}
	\caption{Twisting torus: convergence results for $T=5$ in different norms. We use a higher order geometry: $k_{g,\spa}=k_{g,\ti}=k_{\spa}+1=k_{\ti}+1= k+1$.}
	\label{fig_twisting}
\end{figure}
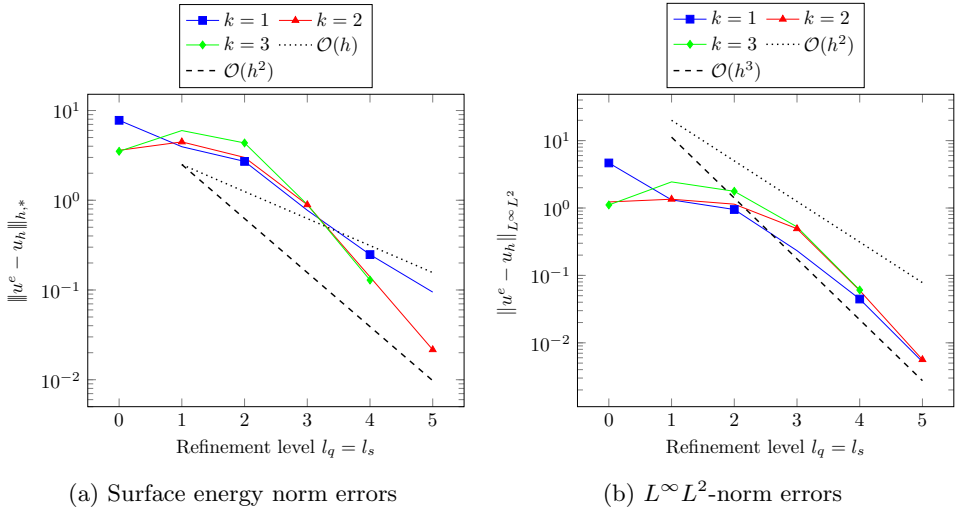
In \Cref{twisting_t25} we observe optimal order convergence for $k\in \{1,2,3\}$ and $T=\frac52$, while for $k=4$ we see faster convergence as expected. This is probably because, for runtime reasons, we have not refined further than $h=2^{-5}$  (level $l=4$)  and the asymptotic range of the space-time convergence order has not  been reached, yet. In \Cref{fig_twisting} we chose a larger temporal end point $T=5$ which leads to much larger curvatures, cf. \Cref{figure_torus_twist}. For $k\in\{1,2\}$ the results show optimal convergence behavior, with in some cases even larger orders. This may caused by   the fact that we use a larger geometry order here ($k_g=k+1$) and due to the very large curvatures the geometric error still dominates on the levels considered. 
The experiments  above indicate that the higher order methods perform well also in the spatially three-dimensional case. 

\begin{figure}[ht!]
		\centering
		\begin{tikzpicture}[scale=0.73]
		\def\varb{30}
		\def\varc{300}
		\begin{groupplot}[
		group style={
			group name=my fancy plots,
			group size=1 by 2,
			xticklabels at=edge bottom,
			vertical sep=0pt
		},
	xlabel={Refinement level $l_{\ti}=l_{\spa}$},
	legend style={ at={(0.5,1.02)}, anchor=south, legend columns =2}, legend cell align=left,
		width=8.5cm,
		xmin=-0.5, xmax=6.5,
		xtick={0,1,2,3,4,5,6}
		]
		
		\nextgroupplot[ymin=0.00001,ymax=100,
		ytick={100,1,0.01,0.0001},
		ymode=log,
		ylabel={$e_{\mass}$},
		axis x line=top, 
		axis y discontinuity=crunch,
		height=6.5cm]

			\addplot +[mark repeat =2, mark phase =1]table [x=mass_list_x, y=-, col sep=comma] {Experiments/Mass_DiagKonvTable_Diss3_zero_alpha_h_interpol_ho_oswald_beta0_stabiwith_h_SpatialNormalVolume_userh_1_tmax_9_tstart_3_xmax_7_xstart_1_ks_1_kt_1_kgs_1_kgt_1_gp_0_tempquadorder_4_tempquadorderDiffu_4.csv};
			\addplot +[mark repeat =2, mark phase =2]table [x=mass_list_x, y=-, col sep=comma] {Experiments/Mass_DiagKonvTable_Diss3_zero_alpha_h_interpol_ho_oswald_beta0_stabiwith_h_SpatialNormalVolume_userh_0_tmax_9_tstart_3_xmax_7_xstart_1_ks_1_kt_1_kgs_1_kgt_1_gp_0_tempquadorder_4_tempquadorderDiffu_4.csv};
			\addplot +[mark repeat =1, mark phase =1]table [x=mass_list_x, y=-, col sep=comma] {Experiments/Mass_DiagKonvTable_Diss3_zero_alpha_h_interpol_ho_oswald_beta1_stabiwith_h_SpatialNormalVolume_userh_1_tmax_9_tstart_3_xmax_7_xstart_1_ks_1_kt_1_kgs_1_kgt_1_gp_0_tempquadorder_4_tempquadorderDiffu_4.csv};
			\addplot +[mark repeat =1, mark phase =1]table [x=mass_list_x, y=-, col sep=comma] {Experiments/Mass_DiagKonvTable_Diss3_zero_alpha_h_interpol_ho_oswald_beta1_stabiwith_h_SpatialNormalVolume_userh_0_tmax_9_tstart_3_xmax_7_xstart_1_ks_1_kt_1_kgs_1_kgt_1_gp_0_tempquadorder_4_tempquadorderDiffu_4.csv};
			\addplot[dotted,line width=0.75pt] coordinates {(0,\varb) (1,\varb*0.5*0.5) (2,\varb*0.25*0.25) (3,\varb*0.125*0.125) (4,\varb*0.0625*0.0625)(5,\varb*0.03125*0.03125)(6,\varb*0.015625*0.015625)};
					\legend{{$\beta=0$, $R=\frac{\bw_S\cdot \bnu_{\partial}}{\alpha_h}$},{$\beta=0$, $R=1$},{$\beta=1$, $R=\frac{\bw_S\cdot \bnu_{\partial}}{\alpha_h}$},{$\beta=1$, $R=1$},$\mathcal{O}(h^{2})$}
		\nextgroupplot[ymin=0,ymax=1e-10,
		ytick={1e-11,1e-13},
		ymode=log,
		axis x line=bottom,
		height=2.5cm]
					\addplot +[mark repeat =2, mark phase =1]table [x=mass_list_x, y=-, col sep=comma] {Experiments/Mass_DiagKonvTable_Diss3_zero_alpha_h_interpol_ho_oswald_beta0_stabiwith_h_SpatialNormalVolume_userh_1_tmax_9_tstart_3_xmax_7_xstart_1_ks_1_kt_1_kgs_1_kgt_1_gp_0_tempquadorder_4_tempquadorderDiffu_4.csv};
		\addplot +[mark repeat =2, mark phase =2]table [x=mass_list_x, y=-, col sep=comma] {Experiments/Mass_DiagKonvTable_Diss3_zero_alpha_h_interpol_ho_oswald_beta0_stabiwith_h_SpatialNormalVolume_userh_0_tmax_9_tstart_3_xmax_7_xstart_1_ks_1_kt_1_kgs_1_kgt_1_gp_0_tempquadorder_4_tempquadorderDiffu_4.csv};
		\addplot +[mark repeat =1, mark phase =1]table [x=mass_list_x, y=-, col sep=comma] {Experiments/Mass_DiagKonvTable_Diss3_zero_alpha_h_interpol_ho_oswald_beta1_stabiwith_h_SpatialNormalVolume_userh_1_tmax_9_tstart_3_xmax_7_xstart_1_ks_1_kt_1_kgs_1_kgt_1_gp_0_tempquadorder_4_tempquadorderDiffu_4.csv};
		\addplot +[mark repeat =1, mark phase =1]table [x=mass_list_x, y=-, col sep=comma] {Experiments/Mass_DiagKonvTable_Diss3_zero_alpha_h_interpol_ho_oswald_beta1_stabiwith_h_SpatialNormalVolume_userh_0_tmax_9_tstart_3_xmax_7_xstart_1_ks_1_kt_1_kgs_1_kgt_1_gp_0_tempquadorder_4_tempquadorderDiffu_4.csv};
		\addplot[dotted,line width=0.75pt] coordinates { (0,\varb) (1,\varb*0.5*0.5) (2,\varb*0.25*0.25) (3,\varb*0.125*0.125) (4,\varb*0.0625*0.0625)(5,\varb*0.03125*0.03125)(6,\varb*0.015625*0.015625)};

		\end{groupplot}
	\end{tikzpicture}
	\caption{Mass conservation results; $k=1$.}
	\label{Twist_mass_loss}
\end{figure}
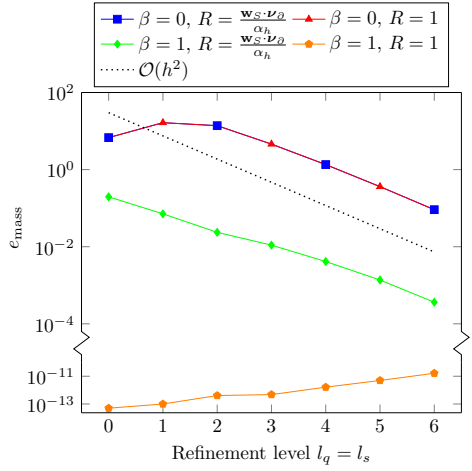

\begin{remark}\label{remark_exact_mass_conservation} \rm 
	We comment on a difference with respect to mass conservation between the methods with $\beta=1$ and $\beta \neq 1$.  
	The method for $\beta=1$ and with $R$ replaced by $1$ is exactly mass conserving for $k_{g,\spa}=1$  in the following sense. Let $u_h$ be the solution of \cref{discreteproblem} with $R=1$ and  $f_h=0$, and define  $i_{\rm mass}^n:=\int_{\Gamma_{\lin}(t_n)}u_{h,-}^n$,  $n=0,\ldots,N$.  Using a piecewise (in time) constant finite element test function
		$ v_h=1$ for  $t\leq t_n$ and $v_h=0$ for $t > t_n$ one easily checks that $i_{\rm mass}^n=i_{\rm mass}^0$ for $n=0,\ldots,N$.
	 Note that we consider $k_{g,\spa}=1$ which implies $\Theta_h^n=\id$ and $\Gamma_h^n(t_n)=\Gamma_h^{n-1}(t_n)=\Gamma_{\lin}(t_n)$. In \Cref{Twist_mass_loss}, for the twisting torus problem with initial condition $u_0(x,y,z)=5+xyz$ and $f_h=0$ we show results for $e_{\rm mass}:= \max_{n=0,\ldots,N}|i_{\rm mass}^n-i_{\rm mass}^0| $. These  confirm the exact mass conservation property.  The results also show that the method with $\beta=1$ leads to significantly better results than the method with $\beta=0$, also for the $R\neq 1$ case used in \eqref{discreteformendetail}. We present results only for $k=1$. For the higher order case with $\beta=1$ we observe  that the discretization errors in the $L^\infty L^2$ norm are of optimal order for $k=3$ but suboptimal for $k=2$. The reason for this (unexpected) behaviour is unclear and has to be investigated further.   
	\end{remark}

\subsection{Space-time surfaces with topological singularities}\label{section_top_singular}
A unique feature of the technique presented in this paper is that it is also applicable to problems in which evolving surfaces with topological singularities occur. In the two examples below we illustrate this robustness property.  Due to the singularities, higher order convergence can not be expected and therefore we restrict to $k=1$. 
\subsubsection{Merging spheres}\label{subsection_smooth_merge}
\begin{figure}[!htbp]
	\begin{minipage}{\textwidth}
		\vspace{-1cm}
		\includegraphics[width=\textwidth]{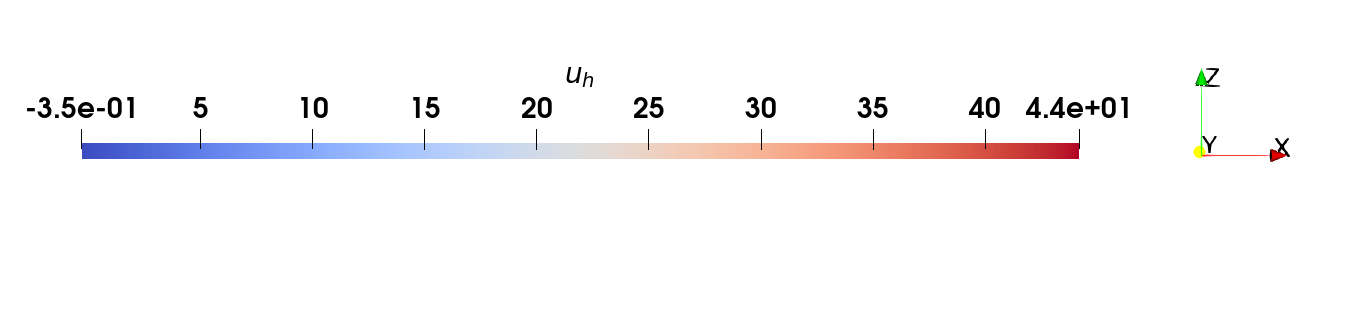}    
		\vspace{-1.5cm}
	\end{minipage}
	\begin{minipage}{0.49\textwidth}
		\includegraphics[width=\textwidth]{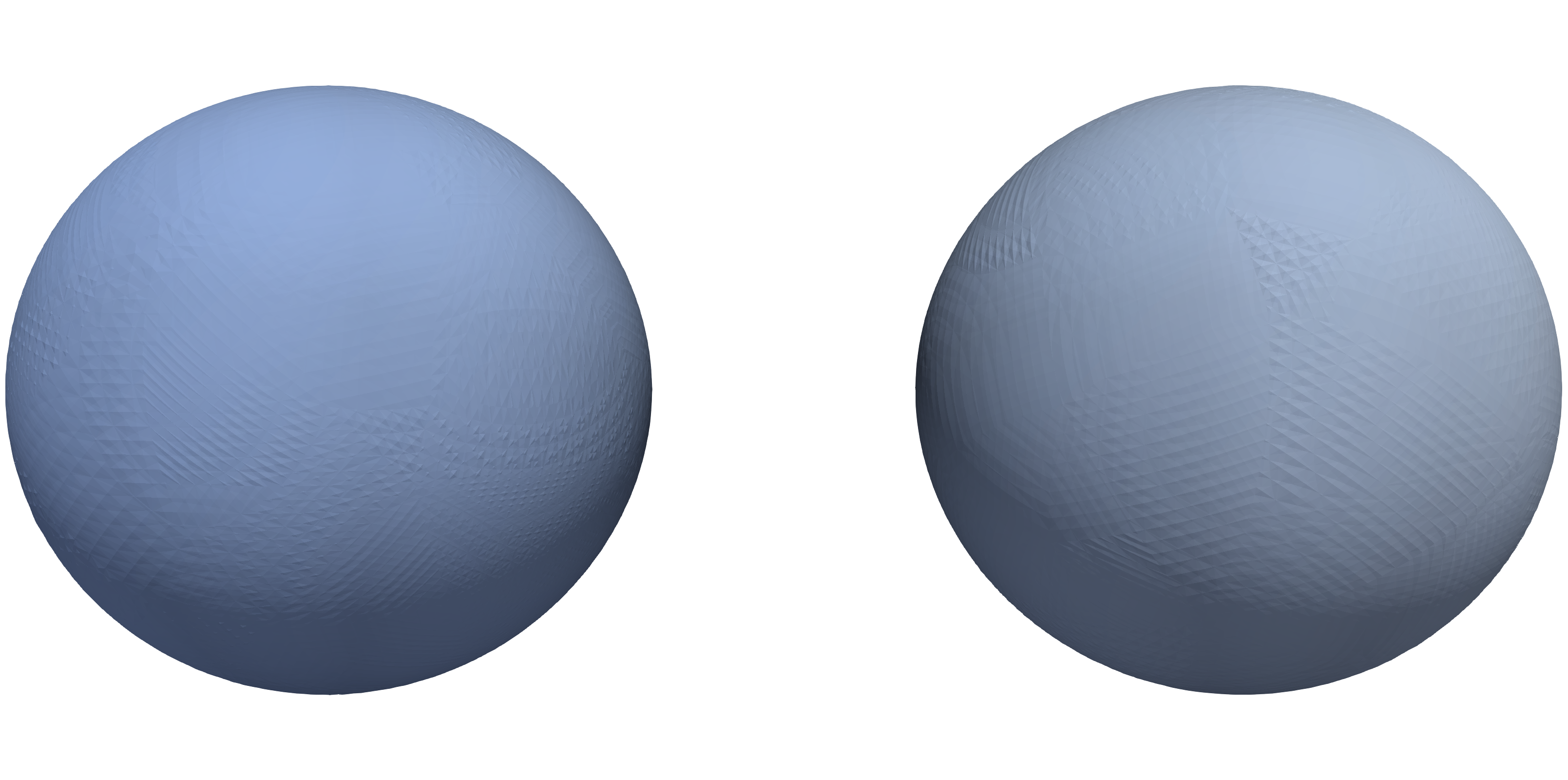}    
	\end{minipage}\hfill
	\begin{minipage}{0.49\textwidth}
		\includegraphics[width=\textwidth]{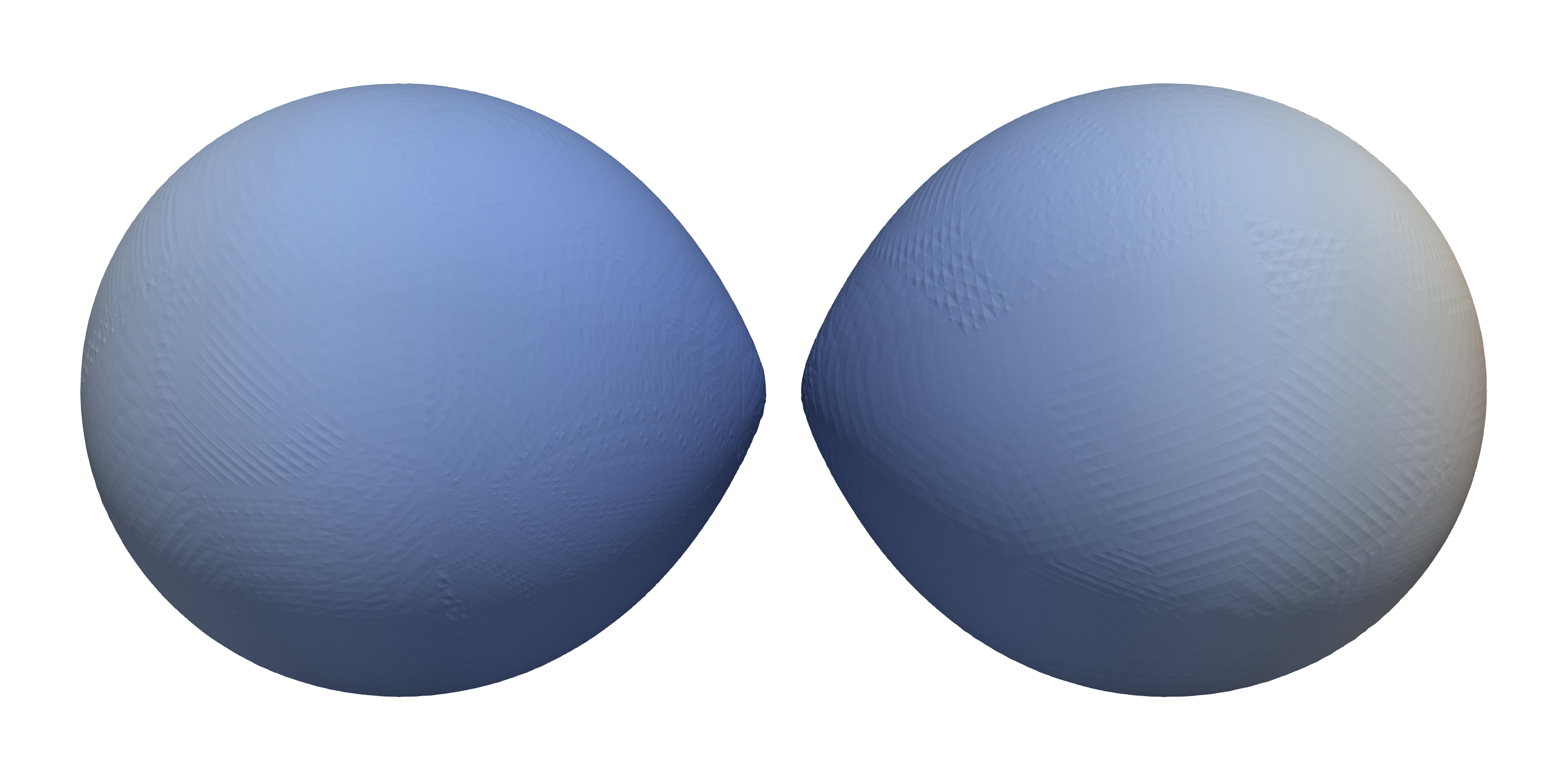}     
	\end{minipage}\hfill
	\begin{minipage}{0.49\textwidth}
		\includegraphics[width=\textwidth]{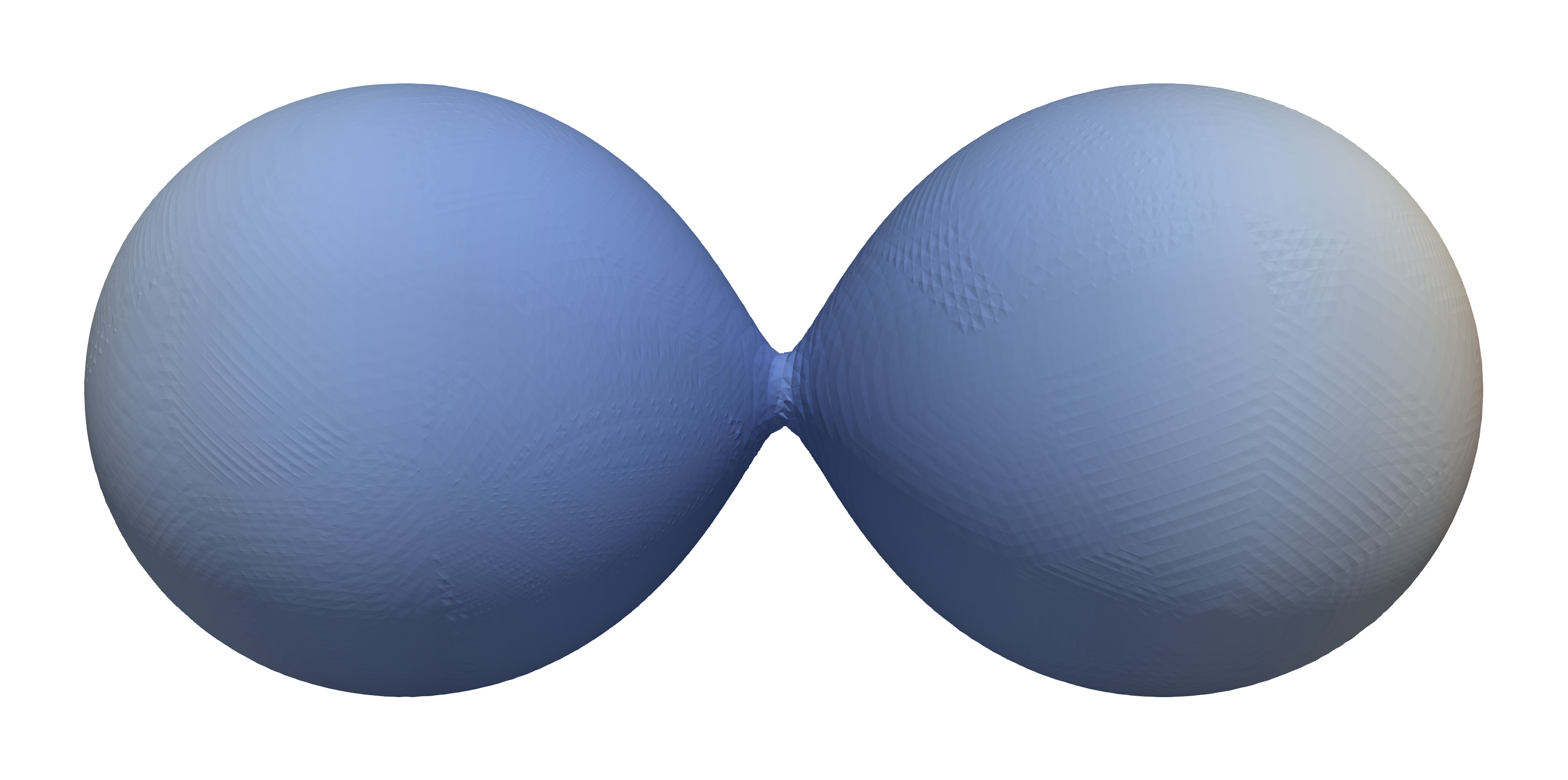}     
	\end{minipage}
	\begin{minipage}{0.49\textwidth}
		\includegraphics[width=\textwidth]{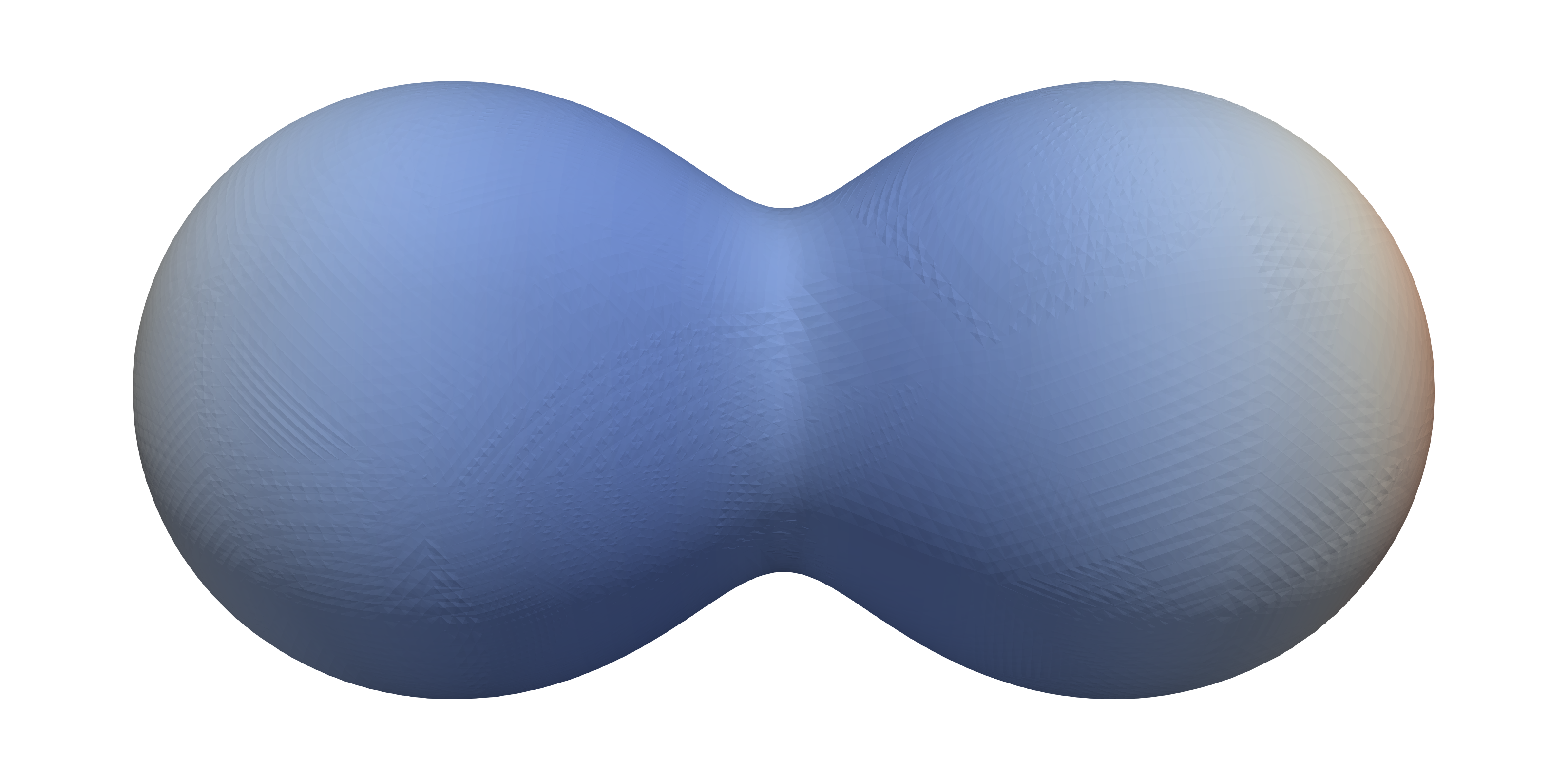}     
	\end{minipage}\hfill
	\begin{minipage}{0.49\textwidth}
		\includegraphics[width=\textwidth]{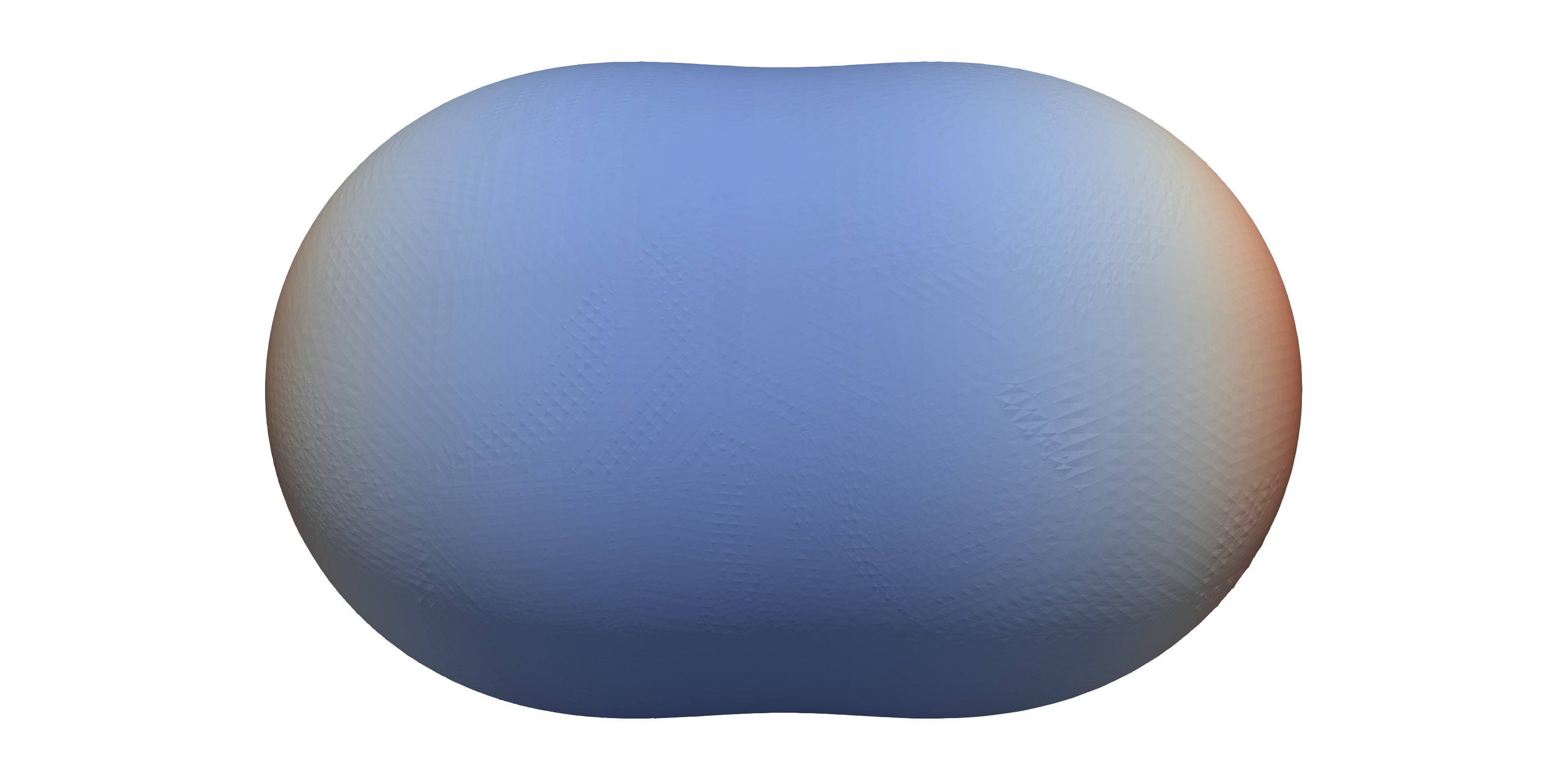}     
	\end{minipage}\hfill
	\begin{minipage}{0.49\textwidth}
		\includegraphics[width=\textwidth]{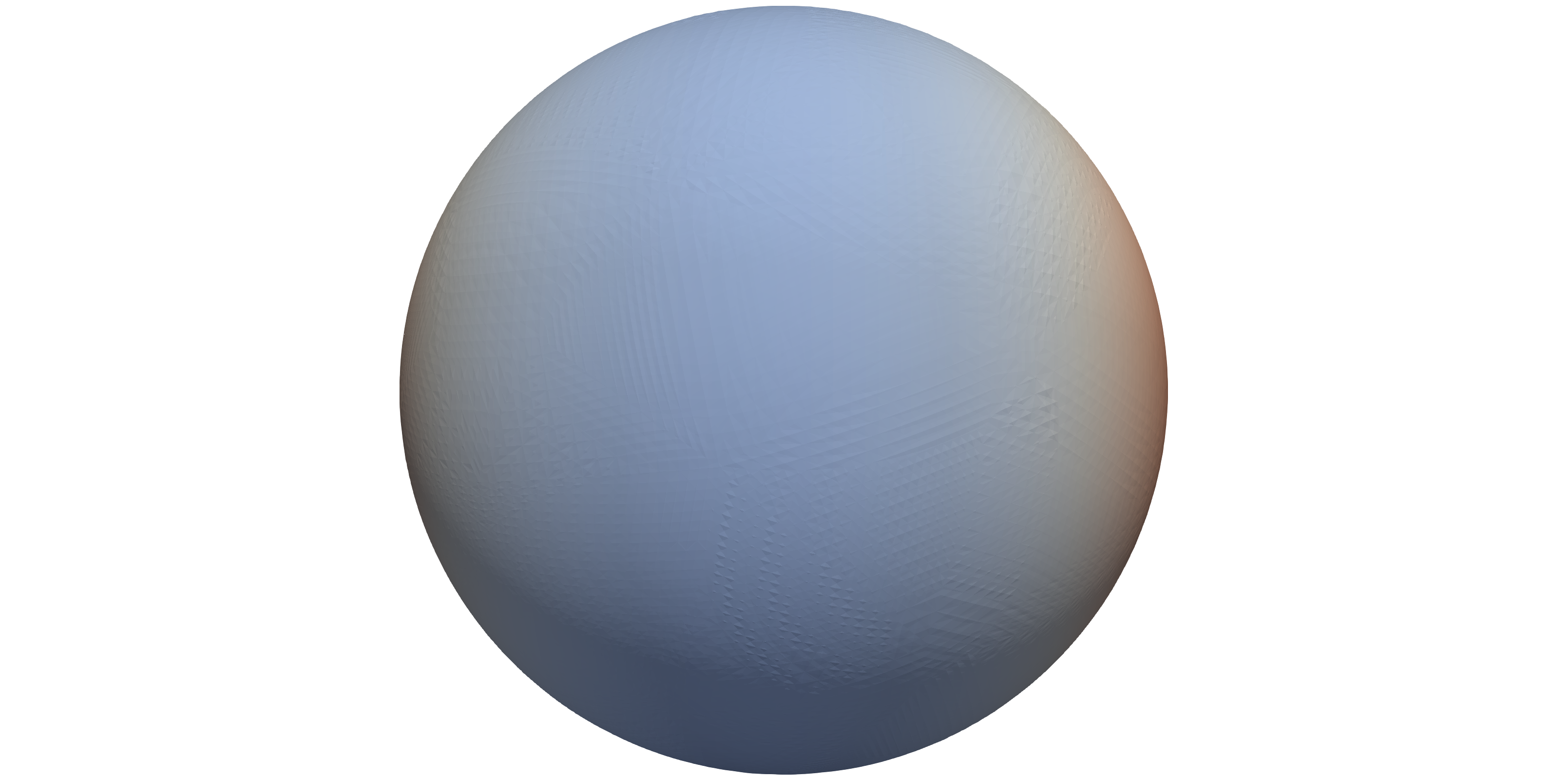}     
	\end{minipage}
	\caption{A sketch of the evolution of the smooth sphere for $k=1$. The discrete solution $u_h$ at times $t\in\{0,\frac{20}{128},\frac{21}{128},\frac{34}{128},\frac{71}{128},1\}$, is depicted in the images from top left to bottom right. Here, we have $\Delta t= 2^{-7}$ and $h=2^{-5}$. The full animation can be found \href{https://doi.org/10.5281/zenodo.7385373}{\textcolor{blue}{here}}.}
	\label{figure_smooth_merge_fine}
\end{figure}

We consider two spheres that merge, similar to \cite{GOReccomas}. We take the background domain $\Omega=(-3,3)^3$ and the time interval $[0,1]$. On ${Q=\Omega\times [0,1]}$ the evolution of the surface $\Gamma(t)$ is described by the zero level of
\begin{equation*}
	\phi(\bx,t)=1-\frac{1}{\norm{\bx-c_+(t)}^3}-\frac{1}{\norm{\bx-c_-(t)}^3},\quad (\bx,t)\in Q,
\end{equation*}
where $c_{\pm}=\pm \frac32(t-1,0,0)^\intercal$. On $Q$ we construct a normal velocity field that transports $\Gamma(t)$ by
\begin{equation*}
	\bw=-\ddt{\phi}{t}\frac{\nabla\phi}{\norm{\nabla \phi}^2}. 
\end{equation*}
We consider the  initial condition 
\begin{equation*}
	u_0(\bx,t)=yz+x+15, \quad(\bx,t)=(x,y,z,t)\in Q
\end{equation*}
and $f=0$ for the right-hand side. We take the initial mesh size as $h_{\mathrm{init}}=2^{-1}$ and the initial time step size as $\Delta t_{\mathrm{init}}=2^{-3}$. 
In  \Cref{figure_smooth_merge_fine} we illustrate the discrete solution at different points in time. 

To evaluate the robustness of the method, we measure the mass evolution over time. We run numerical experiments to assess if and how accurate the discrete solution satisfies a mass conservation property. 
To illustrate the change of total mass and area over time, we measure these quantities at the discrete points $t_n$, i.e. 
\begin{equation*}
	i^n_{\mass}\coloneqq \int_{\Gamma_h^n(t_n)}u_{h,-}^n \dif s_h,\quad i^n_{\surf}\coloneqq \int_{\Gamma_h^n(t)}1 \dif s_h, \quad n=0,\dots,N.
\end{equation*}
In the figures below these quantities are linearly interpolated in each time interval, resulting in the continuous time dependent functions $i_{\mass}$ and $i_{\surf}$. 
We measure the maximum numerical mass change using the quantity
\begin{equation}
	e_{\mass}\coloneqq \max_{n=0,\dots,N}\abs{	i^n_{\mass} - i^0_{\mass}}.\label{discmasschange}
\end{equation}
In \Cref{Smooth_merge_over_time} we illustrate the change of mass and surface area over time. Recall that $l_q$ and $l_s$ denote the refinement level in time and space, respectively. In \Cref{Smooth_merge_mass_loss} we depict the maximal mass error quantity $e_{\rm mass}$. 
\begin{figure}[!htbp]
	\begin{subfigure}{0.5\textwidth}
		\centering
		\begin{tikzpicture}[scale=0.73]
			\def\varb{1}
			\def\vara{1}
			\begin{axis}[xlabel={t}, ylabel={$i_{\mass}$},legend style={ at={(0.5,1.02)}, anchor=south, legend columns =2}, legend cell align=left,/pgf/number format/1000 sep={}]
				\addplot +[mark repeat = 16] table [x=mass_list_x, y=-, col sep=comma] {Experiments/Mass_Vector_lt_5_lx_3_Smooth_Merge_alpha_h_beta0_stabiwith_h_SpatialNormalVolume_userh_1_tmax_5_tstart_5_xmax_3_xstart_3_ks_1_kt_1_kgs_1_kgt_1_gp_0_tempquadorder_4_tempquadorderDiffu_4.csv};
				\addplot +[mark repeat = 32] table [x=mass_list_x, y=-, col sep=comma] {Experiments/Mass_Vector_lt_6_lx_4_Smooth_Merge_alpha_h_beta0_stabiwith_h_SpatialNormalVolume_userh_1_tmax_6_tstart_6_xmax_4_xstart_4_ks_1_kt_1_kgs_1_kgt_1_gp_0_tempquadorder_4_tempquadorderDiffu_4.csv};
				\addplot +[mark repeat = 64] table [x=mass_list_x, y=-, col sep=comma] {Experiments/Mass_Vector_lt_7_lx_5_Smooth_Merge_alpha_h_beta0_stabiwith_h_SpatialNormalVolume_userh_1_tmax_7_tstart_7_xmax_5_xstart_5_ks_1_kt_1_kgs_1_kgt_1_gp_0_tempquadorder_4_tempquadorderDiffu_4.csv};
				\addplot +[mark repeat = 128] table [x=mass_list_x, y=-, col sep=comma] {Experiments/Mass_Vector_lt_8_lx_6_Smooth_Merge_alpha_h_beta0_stabiwith_h_SpatialNormalVolume_userh_1_tmax_8_tstart_8_xmax_6_xstart_6_ks_1_kt_1_kgs_1_kgt_1_gp_0_tempquadorder_4_tempquadorderDiffu_4.csv};
				\addplot[dashed, purple,line width=0.5pt] coordinates { 
					(0.16,385) (0.16,420)};	
				\legend{{$l_\ti=l_{\spa}=2$},{$l_\ti=l_{\spa}=3$},{$l_\ti=l_{\spa}=4$}, {$l_\ti=l_{\spa}=5$}, top. singularity}
			\end{axis}
		\end{tikzpicture}
		\caption{Total mass over time.}
		\label{Smooth_merge_mass} 
	\end{subfigure}\hfill
	\begin{subfigure}{0.5\textwidth}
		\centering
		\begin{tikzpicture}[scale=0.73]
			\def\varb{1}
			\def\vara{1}
			\begin{axis}[xlabel={t}, ylabel={$i_{\surf}$},legend style={ at={(0.5,1.02)}, anchor=south, legend columns =2}, legend cell align=left]
				
				\addplot +[mark repeat = 32] table [x=surface_list_x, y=-, col sep=comma] {Experiments/Surface_Vector_lt_5_lx_3_Smooth_Merge_alpha_h_beta0_stabiwith_h_SpatialNormalVolume_userh_1_tmax_5_tstart_5_xmax_3_xstart_3_ks_1_kt_1_kgs_1_kgt_1_gp_0_tempquadorder_4_tempquadorderDiffu_4.csv};
				\addplot +[mark repeat = 64] table [x=surface_list_x, y=-, col sep=comma] {Experiments/Surface_Vector_lt_6_lx_4_Smooth_Merge_alpha_h_beta0_stabiwith_h_SpatialNormalVolume_userh_1_tmax_6_tstart_6_xmax_4_xstart_4_ks_1_kt_1_kgs_1_kgt_1_gp_0_tempquadorder_4_tempquadorderDiffu_4.csv};
				\addplot +[mark repeat = 128] table [x=surface_list_x, y=-, col sep=comma] {Experiments/Surface_Vector_lt_7_lx_5_Smooth_Merge_alpha_h_beta0_stabiwith_h_SpatialNormalVolume_userh_1_tmax_7_tstart_7_xmax_5_xstart_5_ks_1_kt_1_kgs_1_kgt_1_gp_0_tempquadorder_4_tempquadorderDiffu_4.csv};
				\addplot +[mark repeat = 256] table [x=surface_list_x, y=-, col sep=comma] {Experiments/Surface_Vector_lt_8_lx_6_Smooth_Merge_alpha_h_beta0_stabiwith_h_SpatialNormalVolume_userh_1_tmax_8_tstart_8_xmax_6_xstart_6_ks_1_kt_1_kgs_1_kgt_1_gp_0_tempquadorder_4_tempquadorderDiffu_4.csv};
				\addplot[dashed,purple,line width=0.5pt] coordinates { 
					(0.16,20) (0.16,28)};	
				\legend{{$l_\ti=l_{\spa}=2$},{$l_\ti=l_{\spa}=3$},{$l_\ti=l_{\spa}=4$}, {$l_\ti=l_{\spa}=5$}, top. singularity}
			\end{axis}
		\end{tikzpicture}
		\caption{Surface area over time.}
		\label{Smooth_merge_surface} 
	\end{subfigure}
	\caption{Merging spheres: the evolution of mass and surface area for $k=1$.}
	\label{Smooth_merge_over_time}
\end{figure}
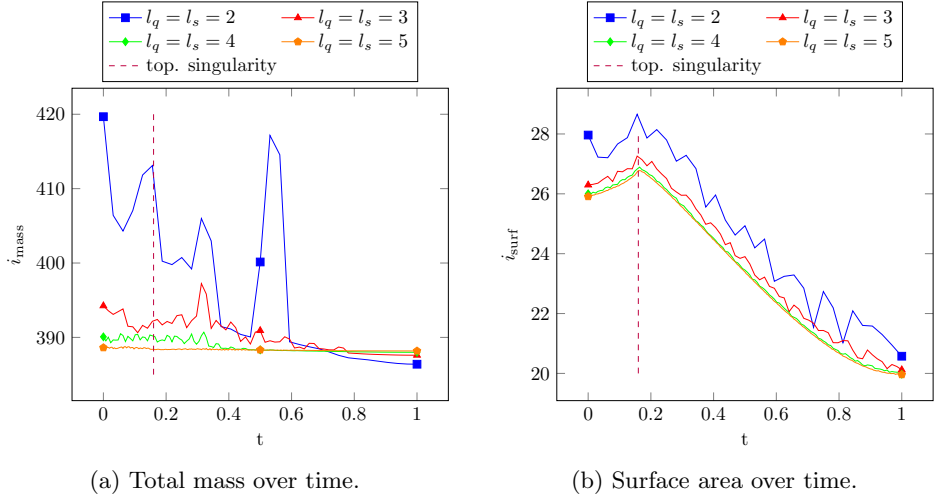
\begin{figure}[!htbp]
	\centering
	\begin{tikzpicture}[scale=0.73]
		\def\varb{800}
		\begin{semilogyaxis}[xlabel={Refinement level $l_{\ti}=l_{\spa}$}, ylabel={$e_{\mass}$},legend style={ at={(0.5,1.02)}, anchor=south, legend columns =2}, legend cell align=left]
			\addplot +[mark repeat =1, mark phase =1]table [x=mass_list_x, y=-, col sep=comma] {Experiments/Mass_DiagKonvTable_Smooth_Merge_alpha_h_beta0_stabiwith_h_SpatialNormalVolume_userh_1_tmax_8_tstart_3_xmax_6_xstart_1_ks_1_kt_1_kgs_1_kgt_1_gp_0_tempquadorder_4_tempquadorderDiffu_4.csv};
			\addplot[dashed,line width=0.75pt] coordinates { 
				(0,\varb) (1,\varb*0.5*0.5) (2,\varb*0.25*0.25) (3,\varb*0.125*0.125) (4,\varb*0.0625*0.0625)(5,\varb*0.03125*0.03125)
			};
			\legend{$k=1$,$\mathcal{O}(h^{2})$}
		\end{semilogyaxis}
	\end{tikzpicture}
	\caption{The maximum mass error \eqref{discmasschange}. Note that due to the large values of $i_{\mass}$, cf. \Cref{Smooth_merge_over_time}, the relative errors have reasonable sizes.}
		\label{Smooth_merge_mass_loss}
	\end{figure}

	Since the velocity field $\bw$ is not divergence free, the surface area changes over time. The surface area increases  until the spheres collide at $t=1-\frac23 2^{\frac13}\approx 0.16$. After this collision the surface evolves to a sphere  and the surface area decreases as a function of time, cf.  \Cref{Smooth_merge_surface}.  We find that the topological singularity at $t\approx 0.16$ has no significant influence on the approximate mass conservation. Under space and time mesh refinement the discrete  mass $i_{\rm mass}(t)$ converges towards a constant function, cf. \Cref{Smooth_merge_mass}. The maximum mass error $e_{\mass}$ decreases with (optimal) second order, see \Cref{Smooth_merge_mass_loss}.

\subsubsection{Disintegrating sphere}
\begin{figure}[!tbp]
	\begin{minipage}{\textwidth}
		\vspace{-0.6cm}
		\includegraphics[width=\textwidth]{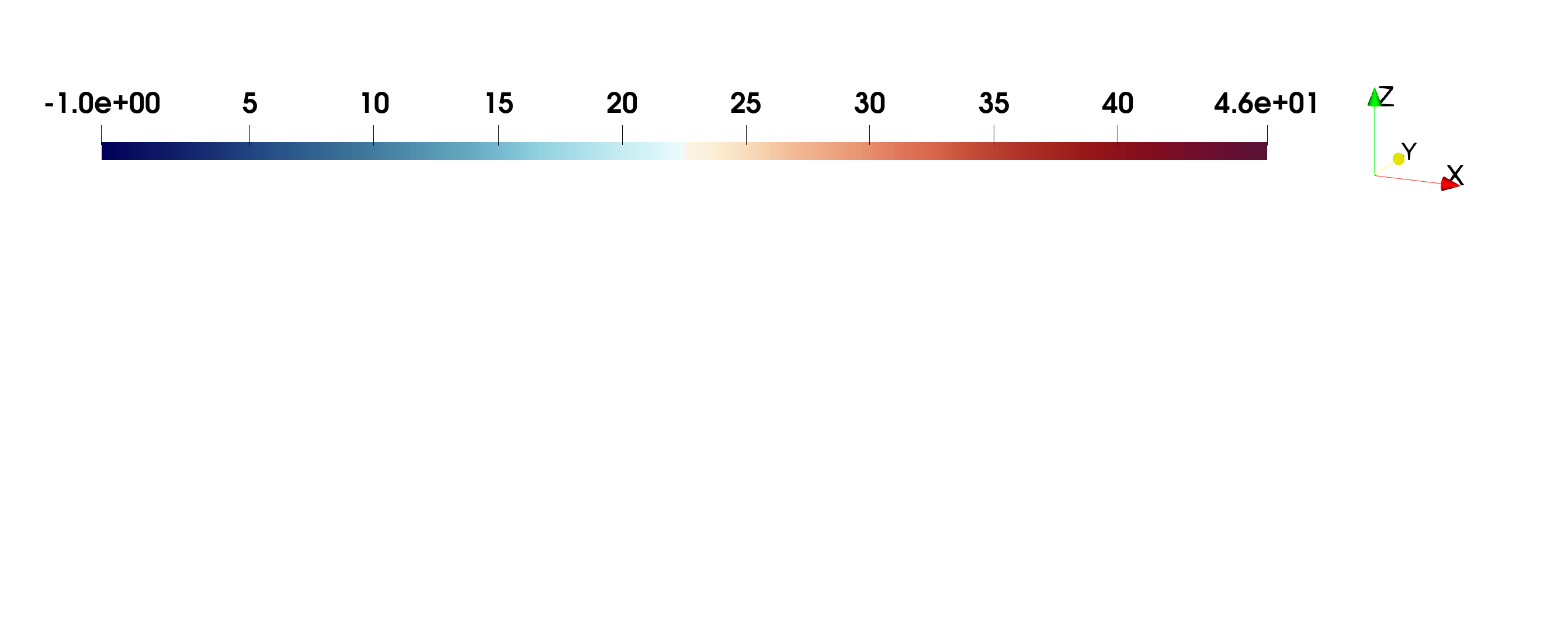}    
		\vspace*{-4.2cm}
	\end{minipage}
	\begin{minipage}{0.33\textwidth}
		\includegraphics[width=\textwidth]{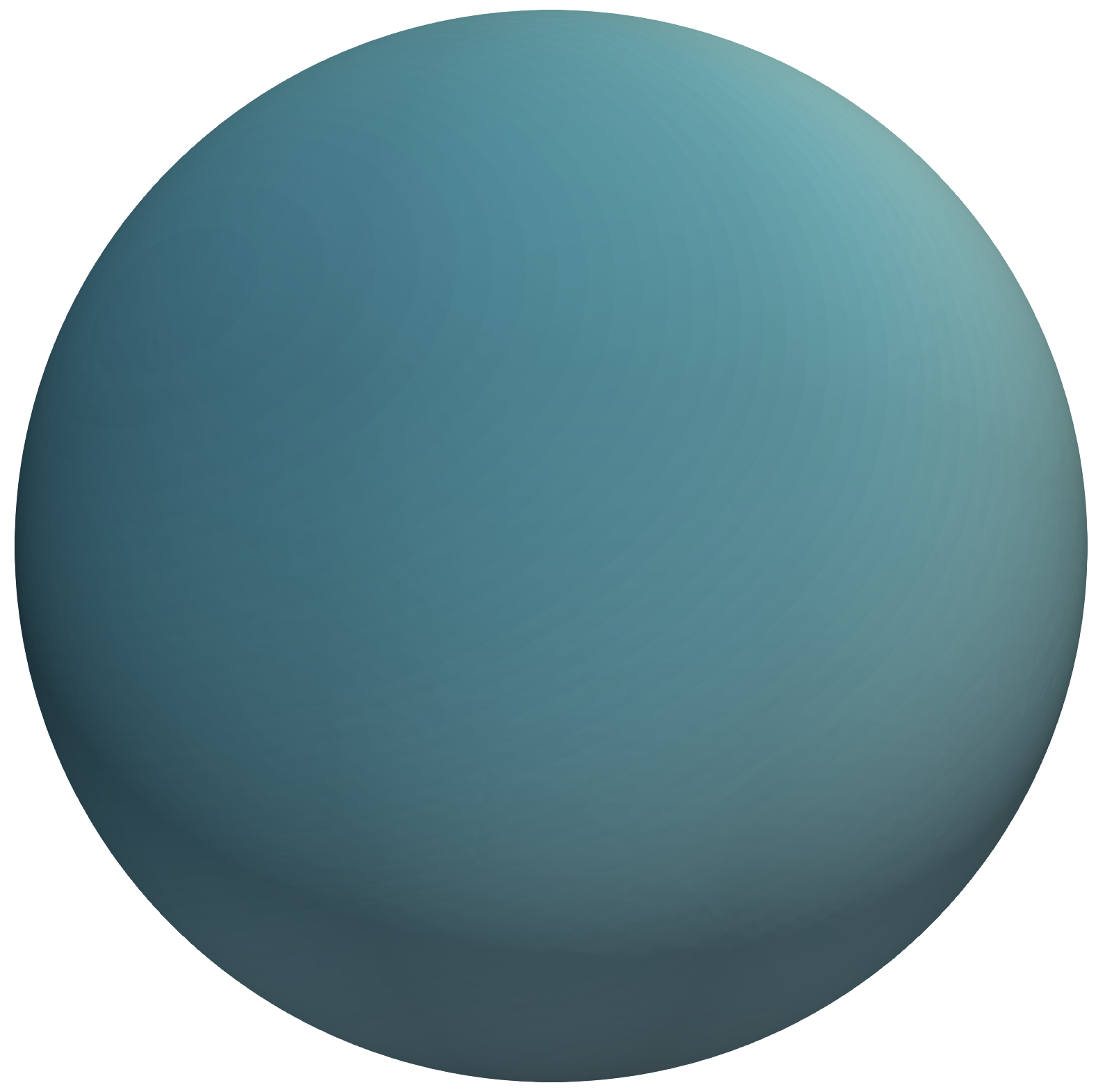}    
	\end{minipage}\hfill
	\begin{minipage}{0.33\textwidth}
		\includegraphics[width=\textwidth]{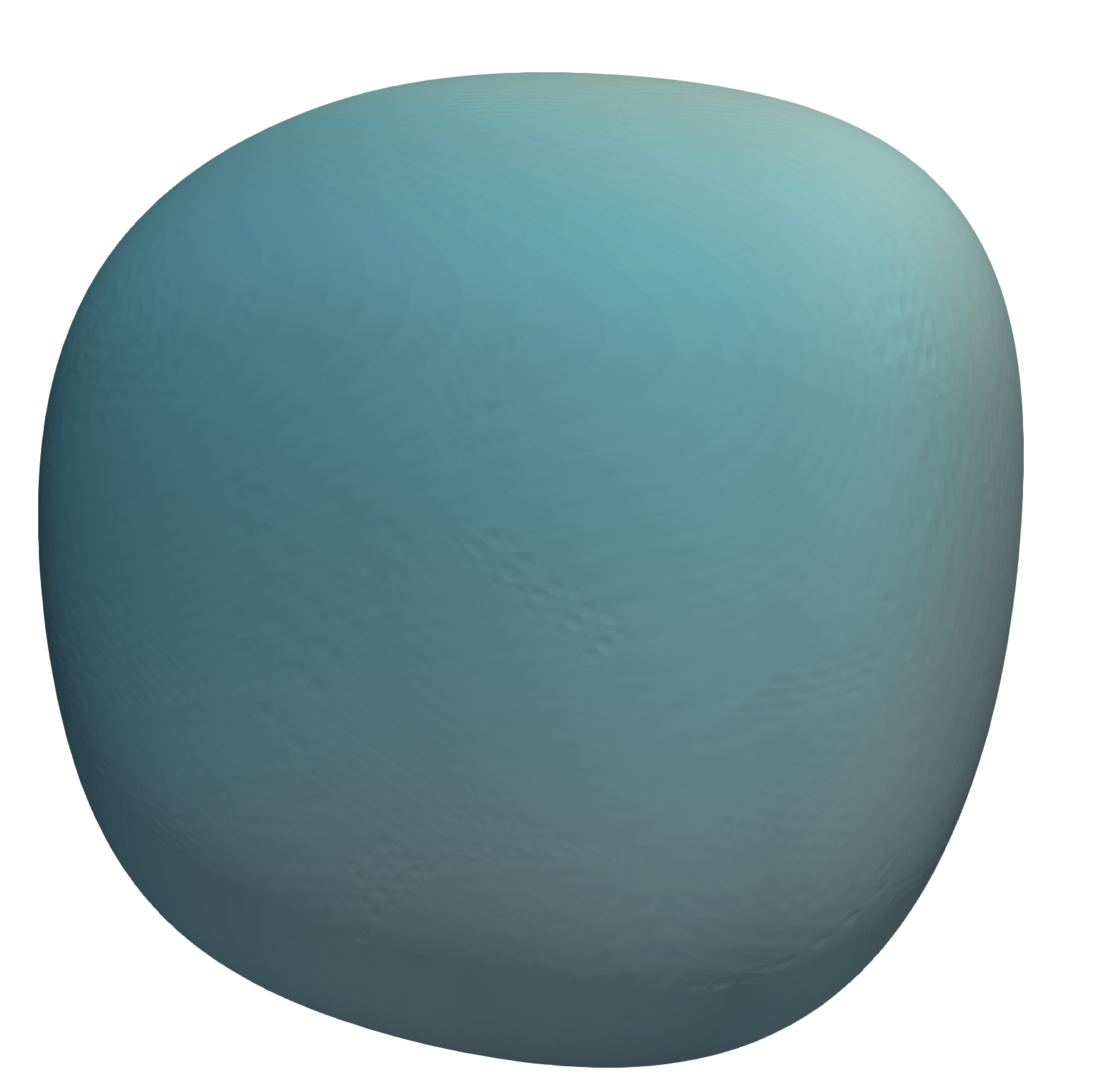}     
	\end{minipage}\hfill
	\begin{minipage}{0.33\textwidth}
		\includegraphics[width=\textwidth]{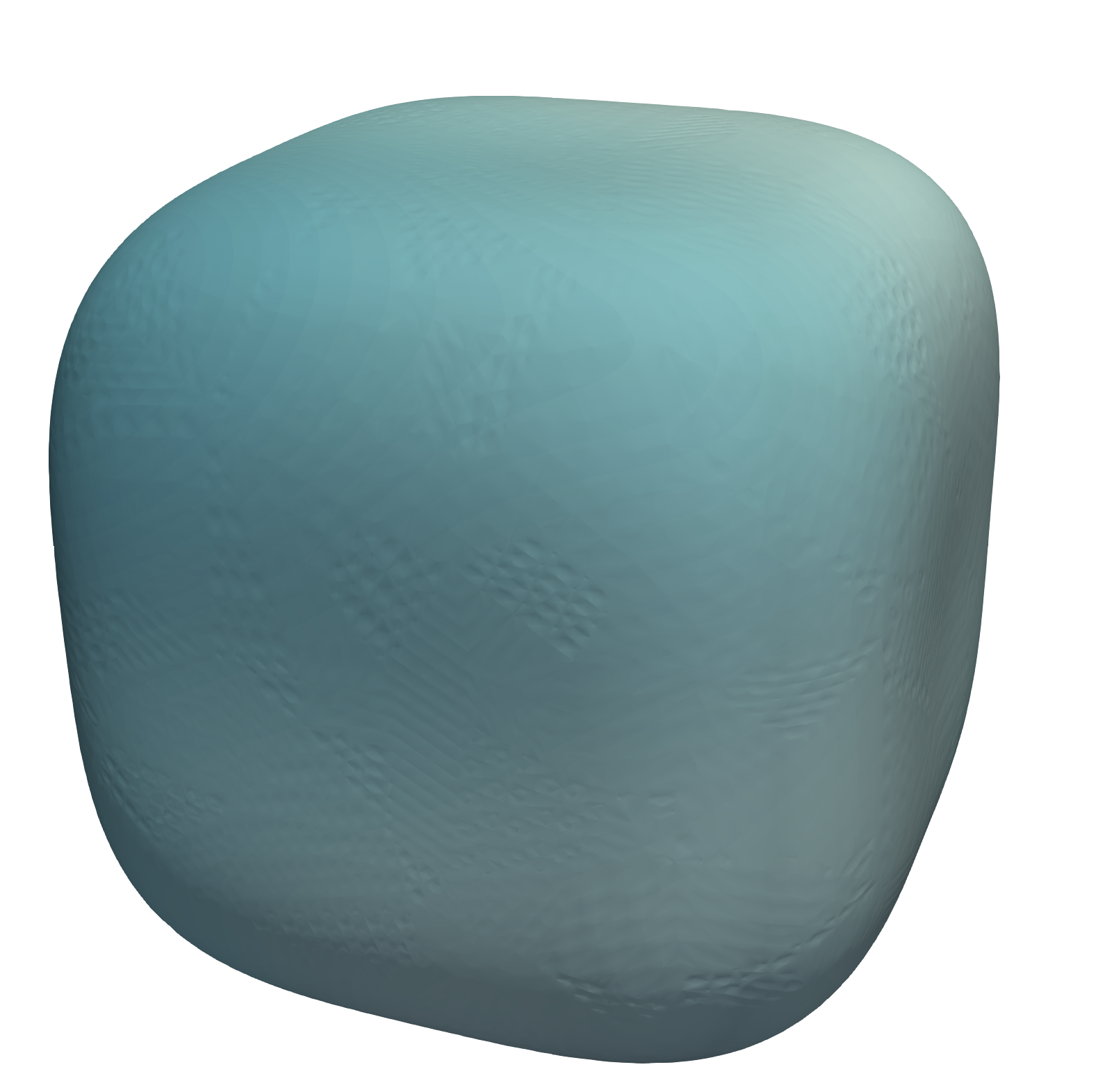}     
	\end{minipage}
	\begin{minipage}{0.33\textwidth}
		\includegraphics[width=\textwidth]{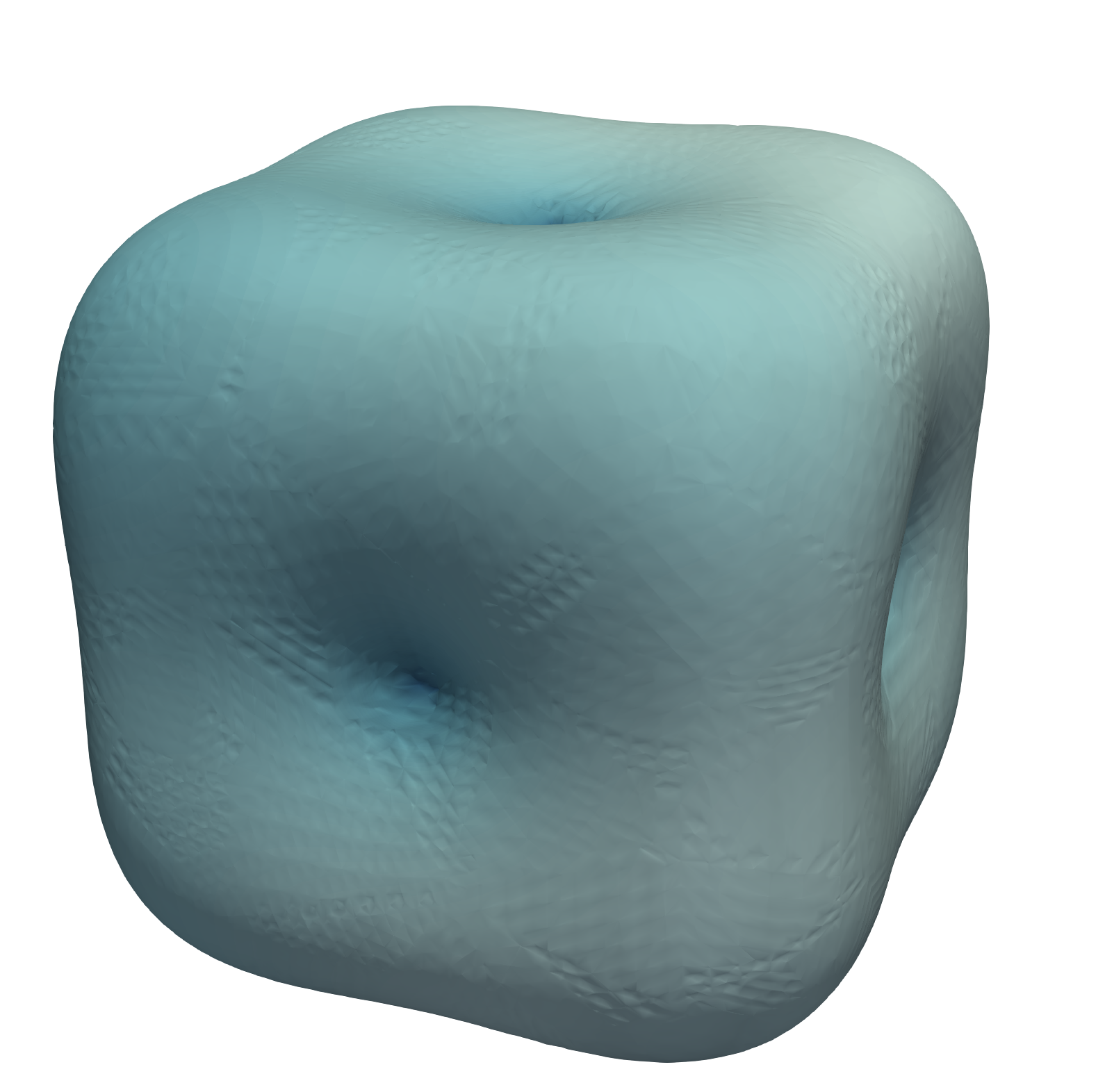}     
	\end{minipage}\hfill
	\begin{minipage}{0.33\textwidth}
		\includegraphics[width=\textwidth]{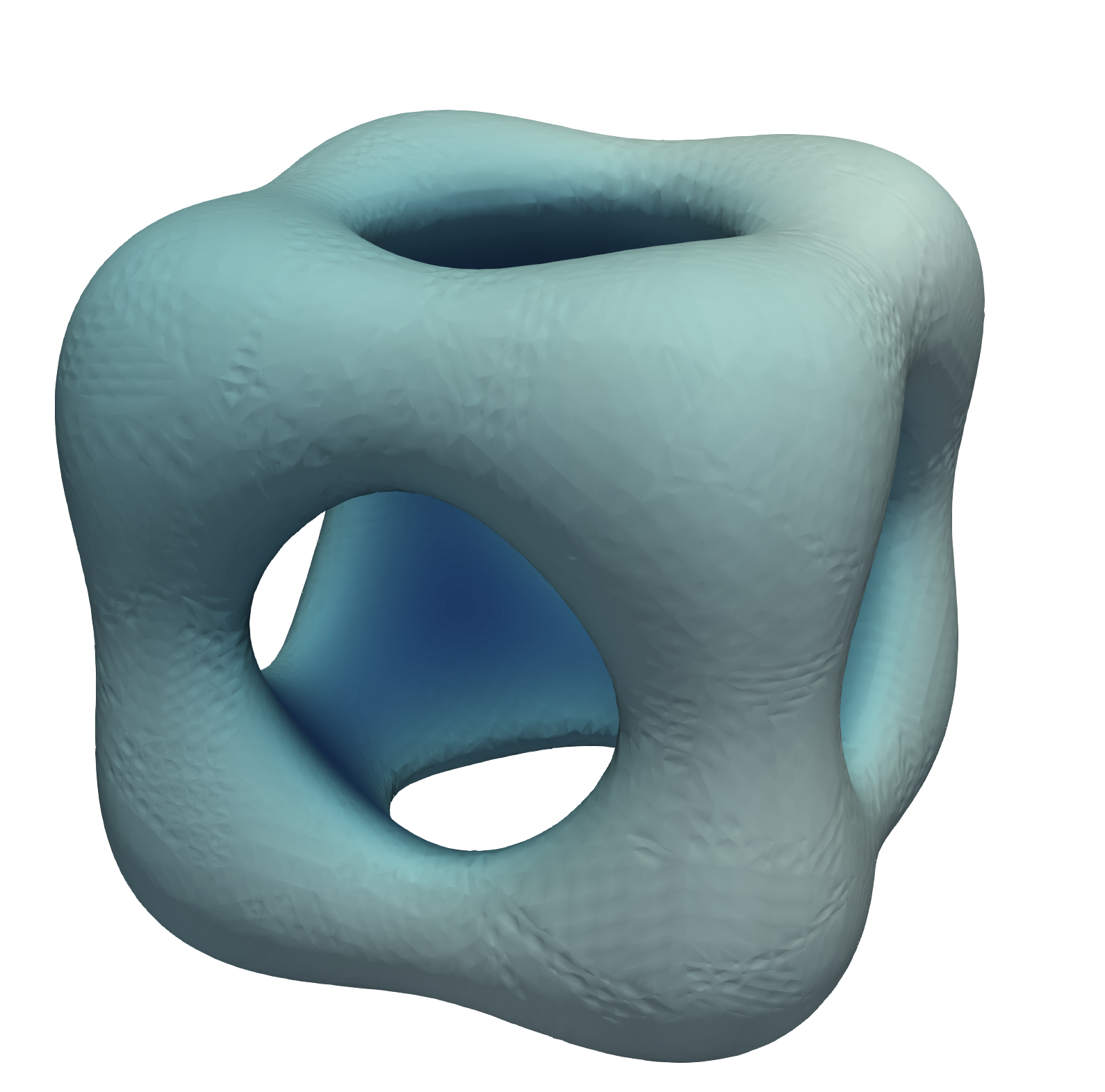}     
	\end{minipage}\hfill
	\begin{minipage}{0.33\textwidth}
		\includegraphics[width=\textwidth]{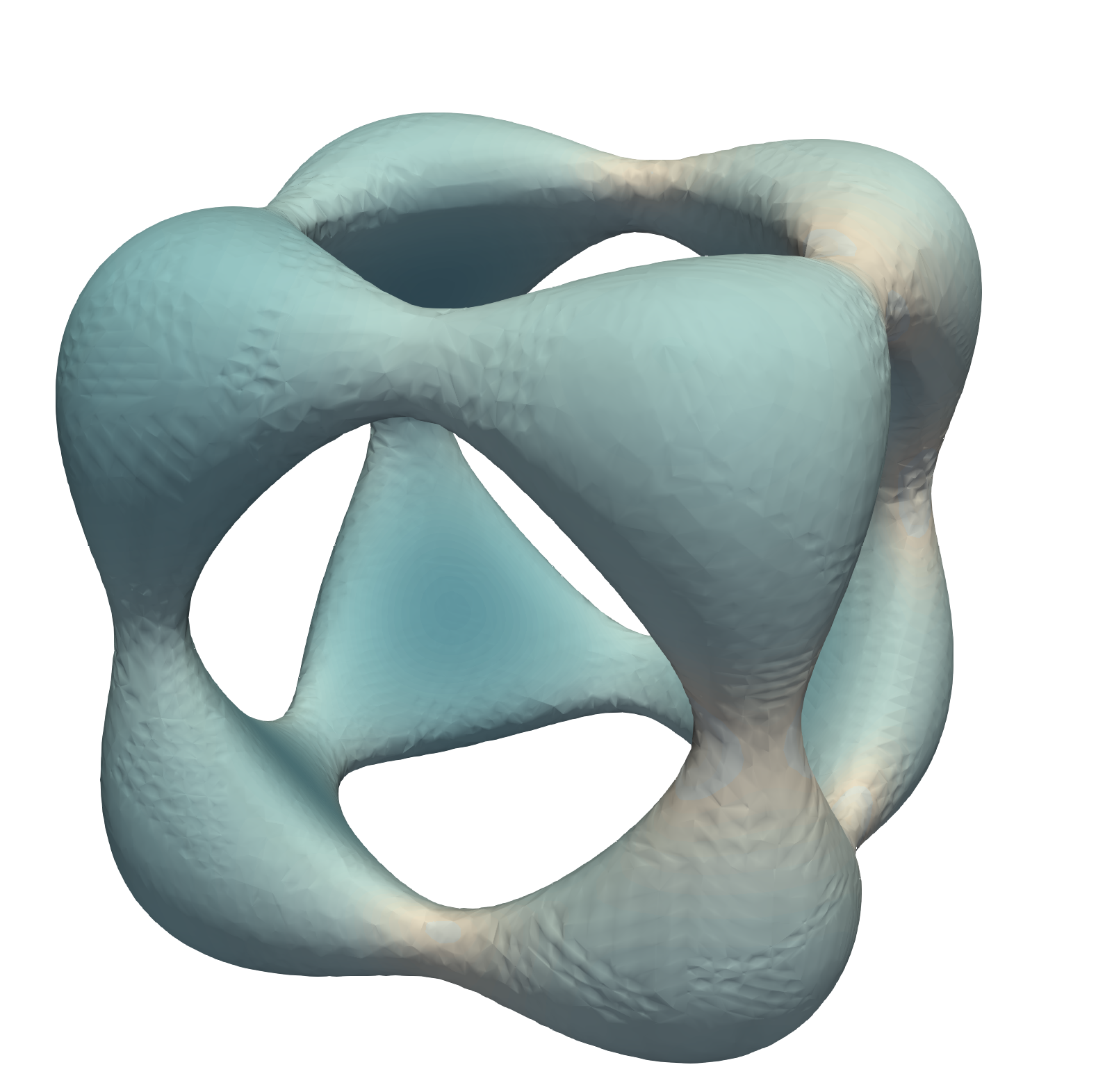}     
	\end{minipage}
	\begin{minipage}{0.33\textwidth}
		\includegraphics[width=\textwidth]{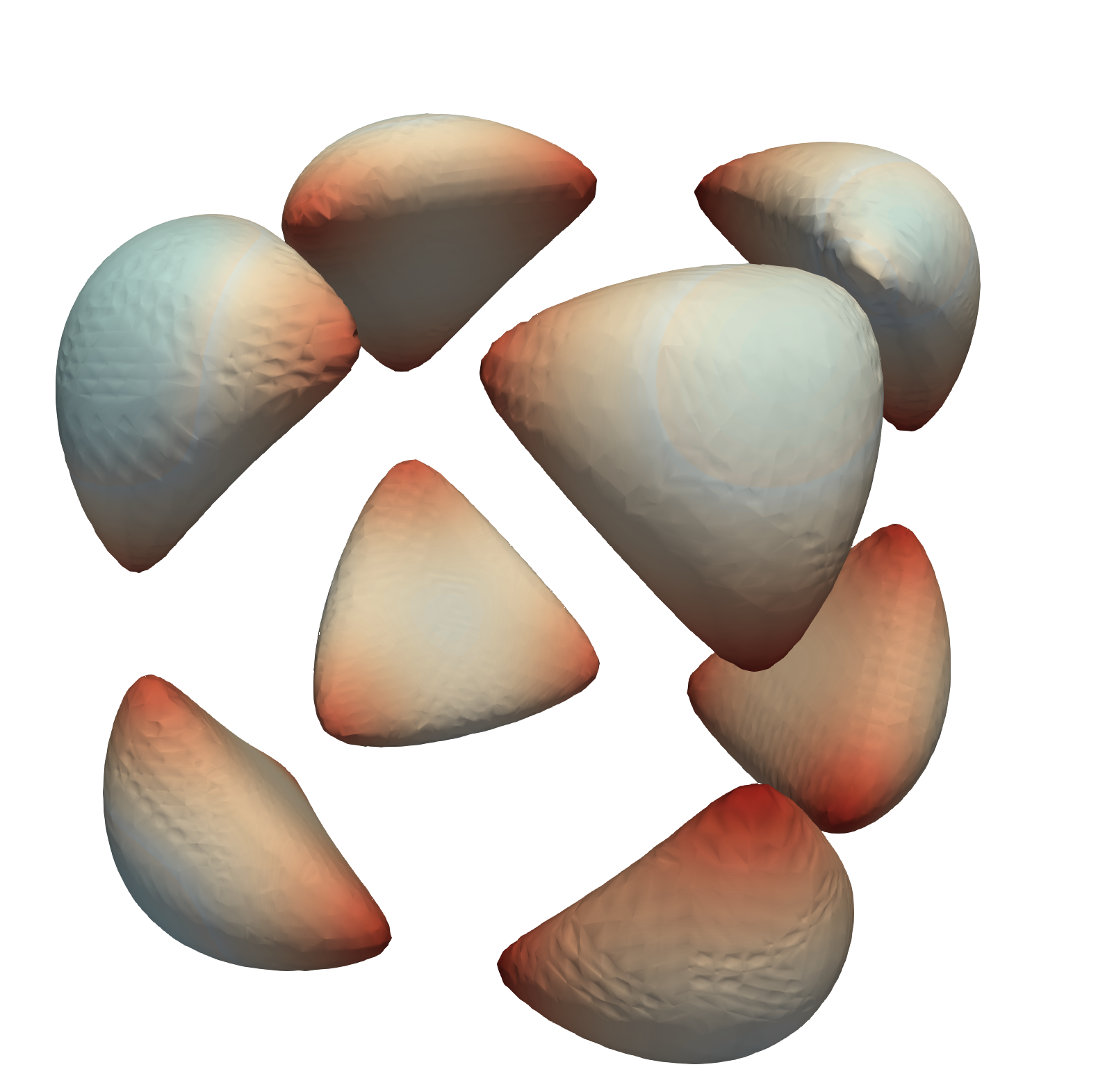}     
	\end{minipage}\hfill
	\begin{minipage}{0.33\textwidth}
		\includegraphics[width=\textwidth]{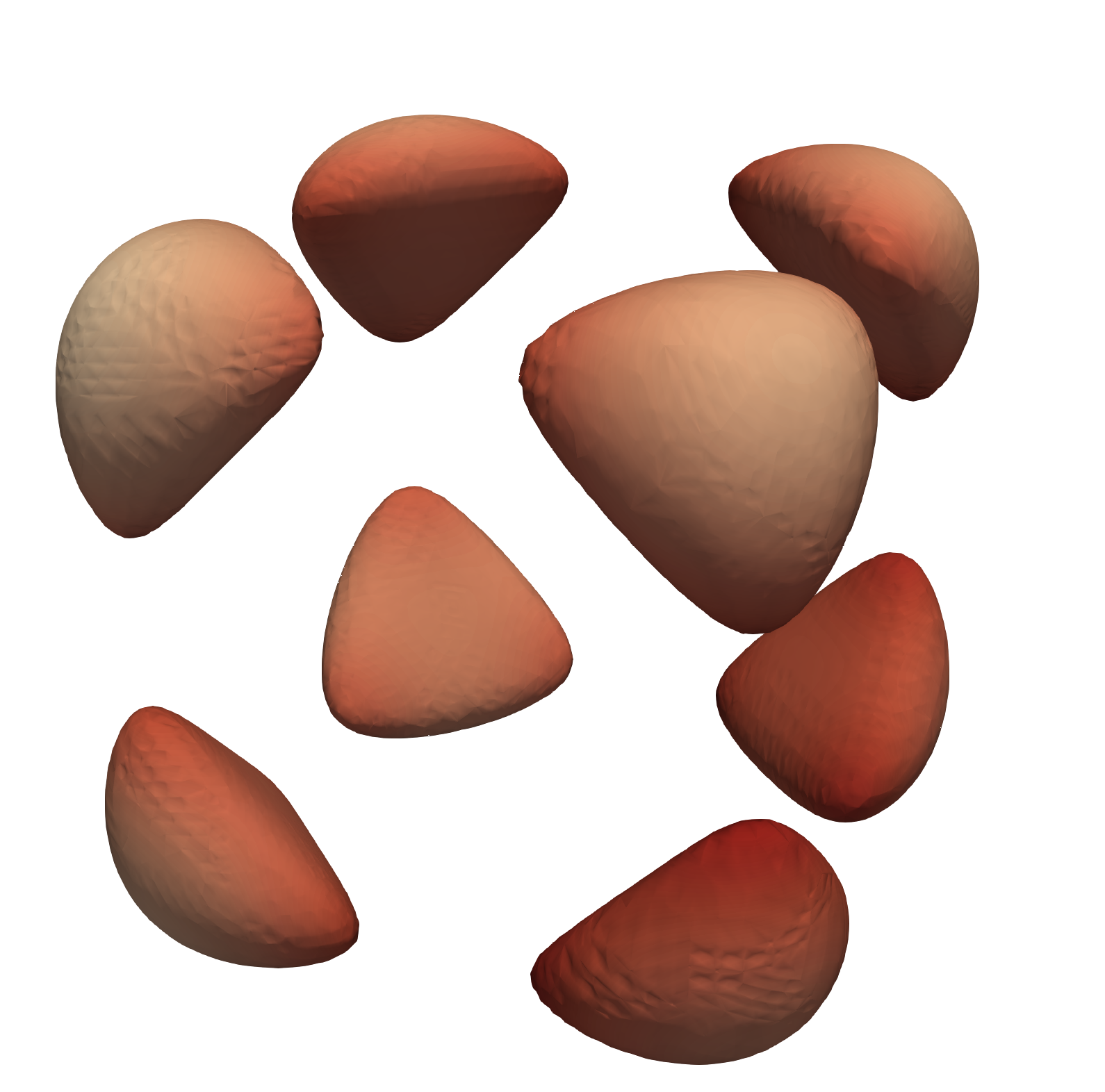}     
	\end{minipage}\hfill
	\begin{minipage}{0.33\textwidth}
		\includegraphics[width=\textwidth]{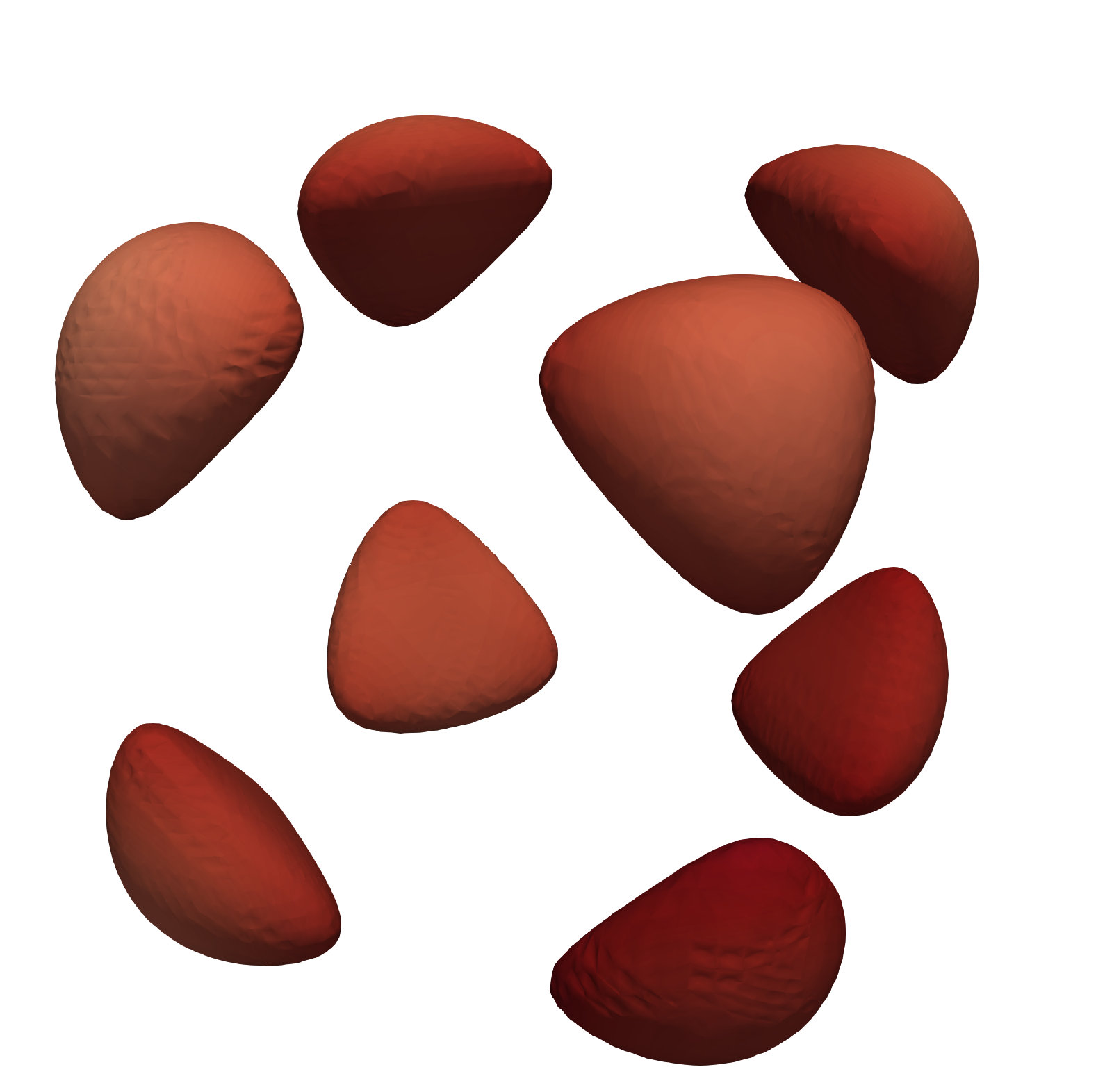}     
	\end{minipage}
	\caption{A sketch of the evolution of the disintegrating sphere. The discrete solution $u_h$ at times $t= \frac{iT}{8}$,  $i\in \{0,\dots,8\}$, is depicted in the images from top left to bottom right. Here, we have $\Delta t=T\cdot2^{-6}$ and $h=2^{-4}$. The full animation can be found \href{https://doi.org/10.5281/zenodo.7385373}{\textcolor{blue}{here}}.}
	\label{figure_deaggro}
\end{figure}
\begin{figure}[!htbp]
	\begin{minipage}{\textwidth}
		\vspace{-0.8cm}
		\includegraphics[width=\textwidth]{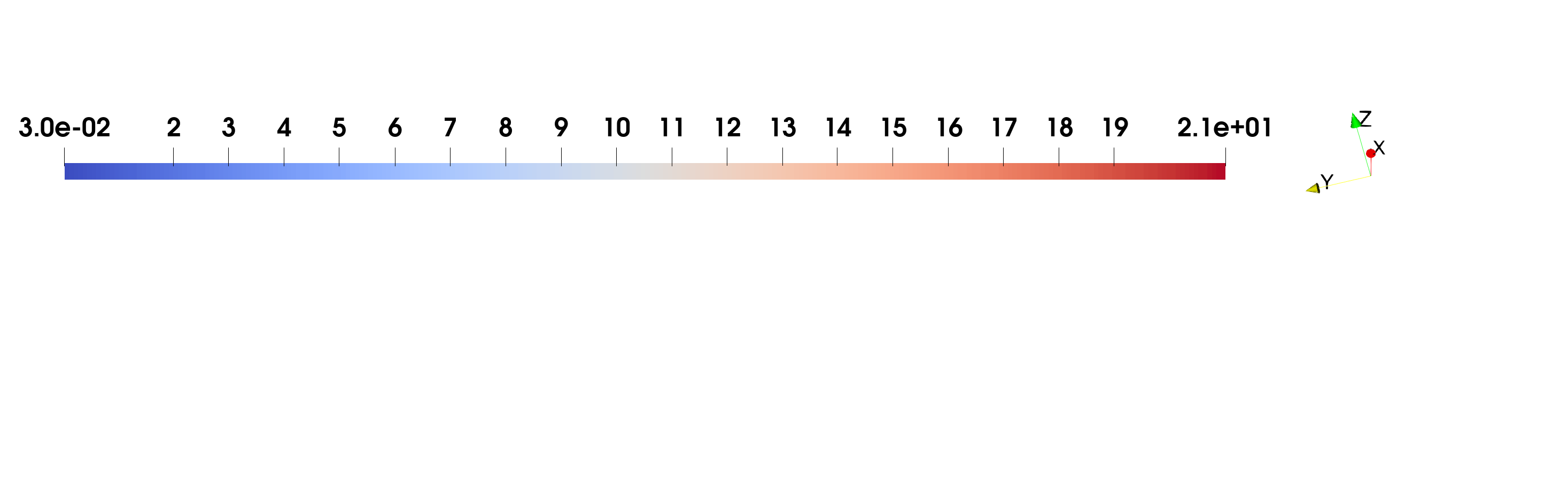}  
		\vspace{-2.7cm}  
	\end{minipage}
	\begin{minipage}{0.495\textwidth}
		\includegraphics[width=\textwidth]{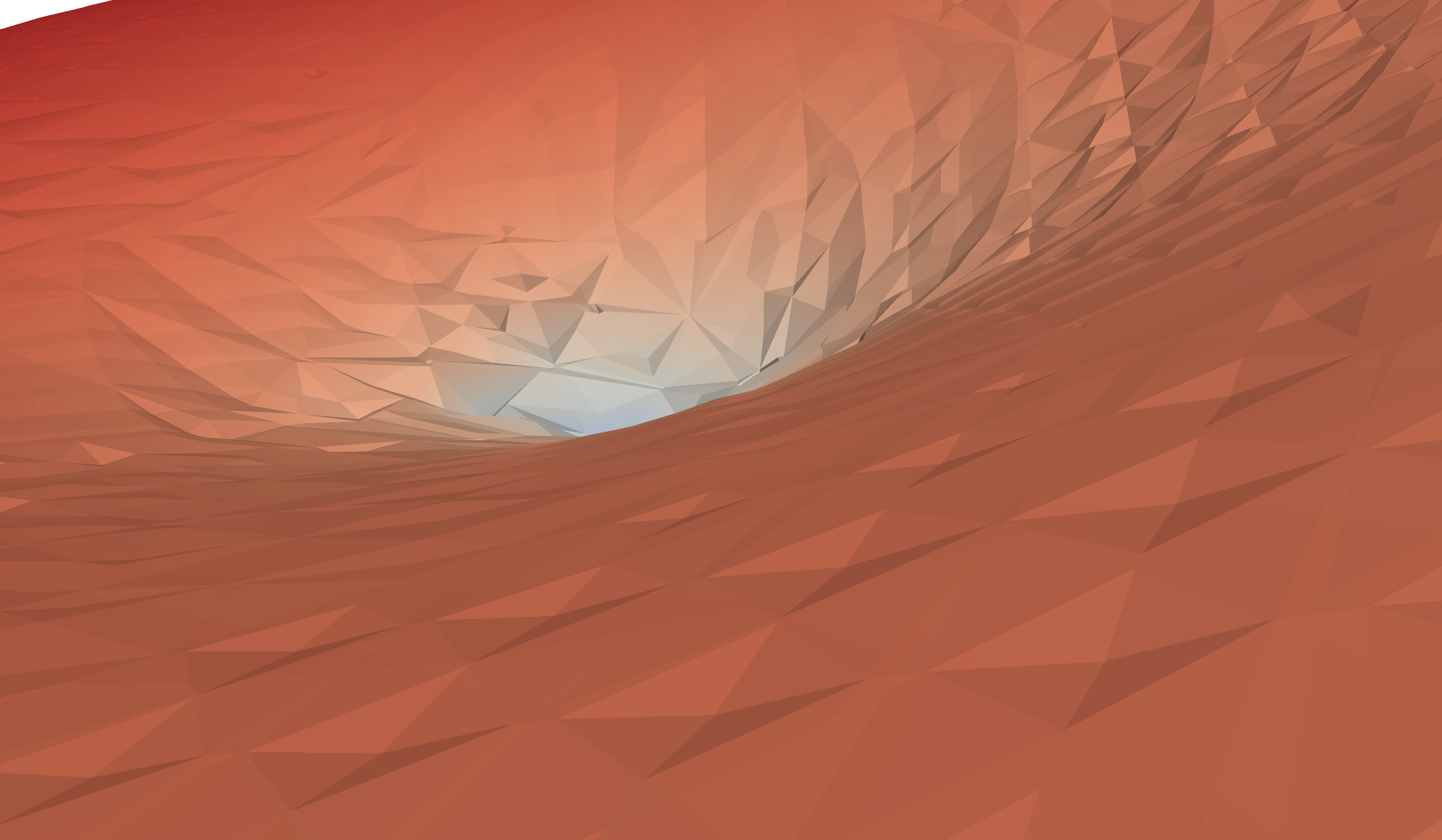}    
	\end{minipage}\hfill
	\begin{minipage}{0.495\textwidth}
		\includegraphics[width=\textwidth]{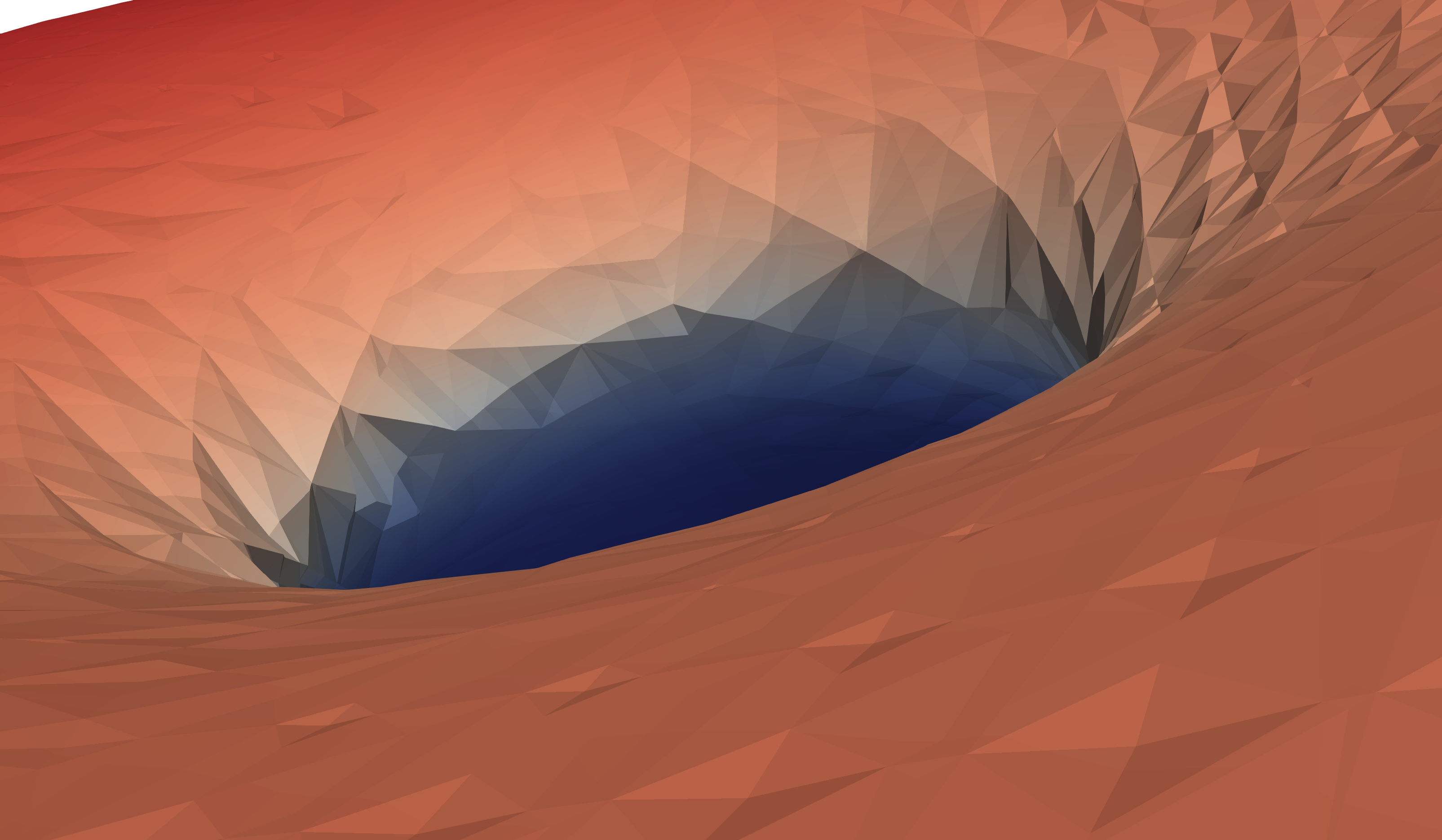}     
	\end{minipage}
	\caption{Disintegrating sphere: on the left side we illustrate the solution on the surface $\Gamma_h(t)$, zoomed in near the spot where a hole emerges later. On the right side we have the same perspective one time step later after the hole arose. The images are colored with the discrete solution $u_h$ at times $t=\frac{24}{64}T$ (left) and $t=\frac{25}{64}T$ (right). Here, we have $\Delta t=T\cdot2^{-6}$ and $h=2^{-4}$.}
	\label{figure_deaggro_whole}
\end{figure}
\begin{figure}[!htbp]
	\begin{minipage}{\textwidth}
		\vspace{-2cm}
		\includegraphics[width=\textwidth]{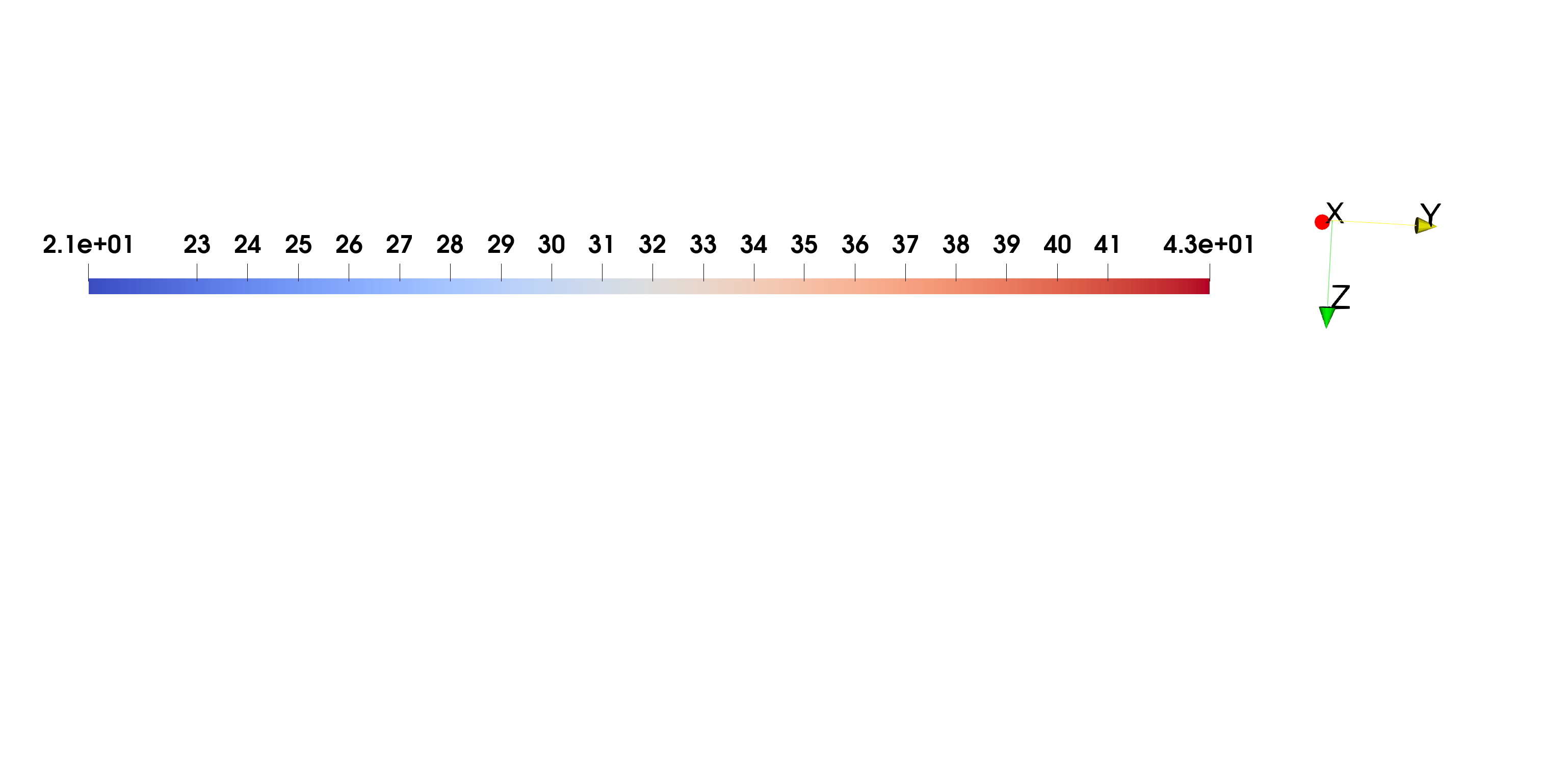}  
		\vspace{-4.2cm}  
	\end{minipage}
	\begin{minipage}{0.495\textwidth}
		\includegraphics[width=\textwidth]{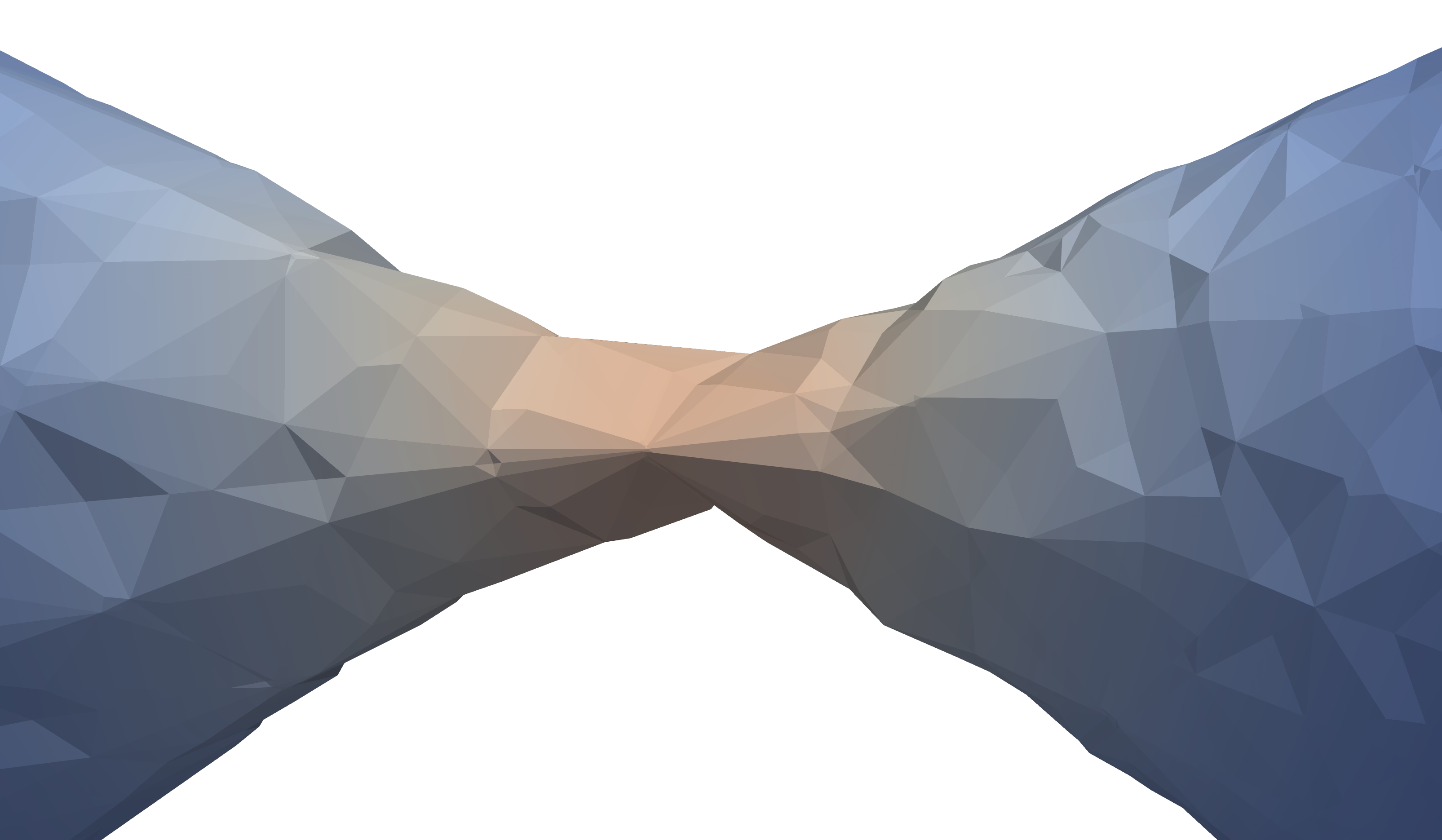}    
	\end{minipage}\hfill
	\begin{minipage}{0.495\textwidth}
		\includegraphics[width=\textwidth]{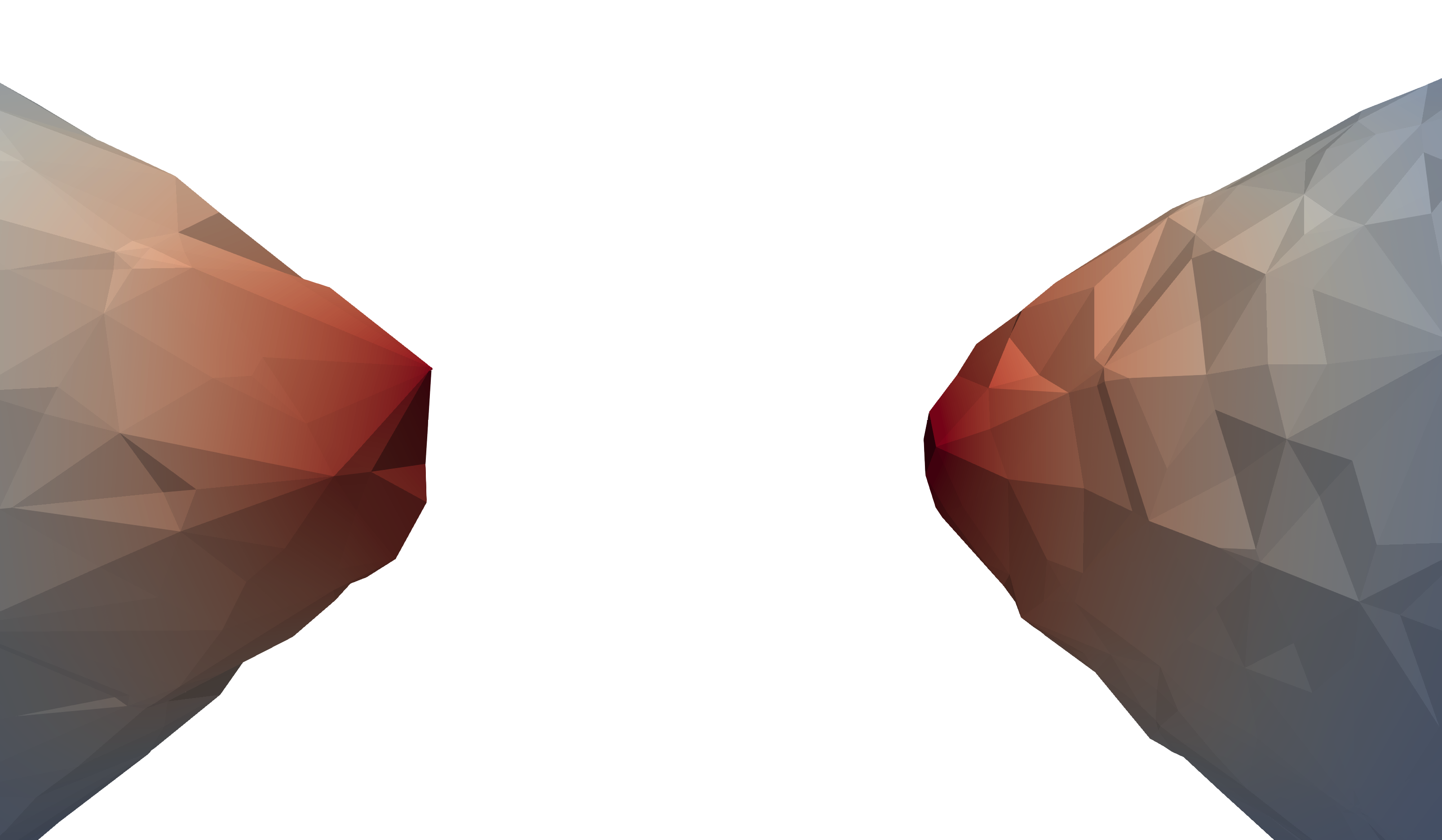}       
	\end{minipage}
	\caption{Disintegrating sphere: on the left side we illustrate the solution on the surface $\Gamma_h(t)$, zoomed in near one of the edges of the decocube, where the surface splits later. On the right side we have the same perspective one time step later after the split happened. The images are colored with the discrete solution $u_h$ at times $t=\frac{44}{64}T$ (left) and $t=\frac{45}{64}T$ (right). Here, we have $\Delta t=T\cdot2^{-6}$ and $h=2^{-4}$.}
	\label{figure_deaggro_split}
\end{figure}

We consider a sphere that gets perforated and splits into eight smaller closed manifolds that are not connected. The level set equation of the perforated manifold is similar to the \textit{decocube} considered in \cite{deckelnick2009h,Grande2017}.
As outer domain we consider $\Omega = (-3,3)^3$ and as time interval we choose $[0,\frac{9}{4}]$. On ${Q=\Omega\times [0,\frac94]}$ the evolution of the surface $\Gamma(t)$ is described by the zero level of the fifth order polynomial
\begin{align*}
	&{\frac { \left( 12{x}^{4}+ \left( 8{y}^{2}+8{z}^{2}-76 \right) {
				x}^{2}+12{y}^{4}+ \left( 8{z}^{2}-76 \right) {y}^{2}+12{z}^{4}-
			76{z}^{2}+205 \right) t}{4}}\\
	&+5{x}^{2}+5{y}^{2}+5{z}^{2}-{\frac
		{125}{4}}\eqqcolon \phi(\bx,t), \quad (\bx,t)=(x,y,z,t)\in Q.
\end{align*}
On $Q$ we construct a normal velocity field that transports $\Gamma(t)$ by
\begin{equation*}
	\bw=-\ddt{\phi}{t}\frac{\nabla\phi}{\norm{\nabla \phi}^2}. 
\end{equation*}
We take the  initial condition 
\begin{equation*}
	u_0(\bx)\coloneqq x+yz+15, \quad \bx\in \Gamma_h(0)
\end{equation*}
and $f=0$ for the right hand side. We choose $k=1$, $h_{\mathrm{init}}=2^{-1}$ and $\Delta t_{\mathrm{init}}=2^{-3}$.
As in \Cref{subsection_smooth_merge} we do not measure the discretization error $u^e-u_h$ because $u$ is not known. In \Cref{figure_deaggro} we show results of numerical simulations.  At $t=0$ we have a sphere which at first deforms to a rounded cube with small bumps. At $t\approx 0.61$ an additional sphere originates from a point singularity in the center of the outer surface. This  sphere gets larger over time and merges with the outer surface at $t\approx 0.87$. The two disjoint surfaces  with genus $0$ merge to a connected surface  with genus $6$ which we call \textit{decocube}. At $t\approx 1.55$ this decocube splits into eight smaller surfaces with genus $0$. Hence, in this example there are three points in time at which a topological singularity occurs. In \Cref{figure_deaggro_whole,figure_deaggro_split} close-ups of the surface near two different topological singularities are shown.

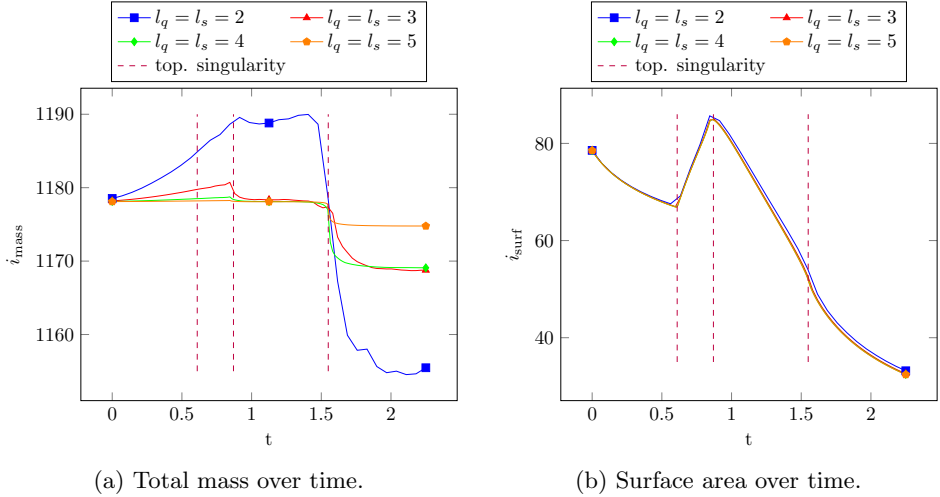
\begin{figure}[!htbp]
	\begin{subfigure}{0.5\textwidth}
		\centering
		\begin{tikzpicture}[scale=0.73]
		\def\varb{1}
		\def\vara{1}
		\begin{axis}[xlabel={t}, ylabel={$i_{\mass}$},legend style={ at={(0.5,1.02)}, anchor=south, legend columns =2}, legend cell align=left,/pgf/number format/1000 sep={}]
		\addplot +[mark repeat = 16] table [x=mass_list_x, y=-, col sep=comma] {Experiments/Mass_Vector_lt_5_lx_3_Diss4_alpha_h_beta0_stabiwith_h_SpatialNormalVolume_userh_1_tmax_5_tstart_5_xmax_3_xstart_3_ks_1_kt_1_kgs_1_kgt_1_gp_0_tempquadorder_4_tempquadorderDiffu_4.csv};
		\addplot +[mark repeat = 32] table [x=mass_list_x, y=-, col sep=comma] {Experiments/Mass_Vector_lt_6_lx_4_Diss4_alpha_h_beta0_stabiwith_h_SpatialNormalVolume_userh_1_tmax_6_tstart_6_xmax_4_xstart_4_ks_1_kt_1_kgs_1_kgt_1_gp_0_tempquadorder_4_tempquadorderDiffu_4.csv};
		\addplot +[mark repeat = 64] table [x=mass_list_x, y=-, col sep=comma] {Experiments/Mass_Vector_lt_7_lx_5_Diss4_alpha_h_beta0_stabiwith_h_SpatialNormalVolume_userh_1_tmax_7_tstart_7_xmax_5_xstart_5_ks_1_kt_1_kgs_1_kgt_1_gp_0_tempquadorder_4_tempquadorderDiffu_4.csv};
		\addplot +[mark repeat = 128] table [x=mass_list_x, y=-, col sep=comma] {Experiments/Mass_Vector_lt_8_lx_6_Diss4_alpha_h_beta0_stabiwith_h_SpatialNormalVolume_userh_1_tmax_8_tstart_8_xmax_6_xstart_6_ks_1_kt_1_kgs_1_kgt_1_gp_0_tempquadorder_4_tempquadorderDiffu_4.csv};
\addplot[dashed, purple,line width=0.5pt] coordinates { 
	(0.87,1155) (0.87,1190)};
\addplot[dashed, purple,line width=0.5pt] coordinates { 
	(1.55,1155) (1.55,1190)};
\addplot[dashed, purple,line width=0.5pt] coordinates { 
	(0.61,1155) (0.61,1190)};
				\legend{{$l_\ti=l_{\spa}=2$},{$l_\ti=l_{\spa}=3$},{$l_\ti=l_{\spa}=4$},{$l_\ti=l_{\spa}=5$}, top. singularity}
		\end{axis}
		\end{tikzpicture}
		\caption{Total mass over time.}
		\label{Dis_sphere_mass} 
	\end{subfigure}\hfill
	\begin{subfigure}{0.5\textwidth}
		\centering
		\begin{tikzpicture}[scale=0.73]
		\def\varb{1}
		\def\vara{1}
		\begin{axis}[xlabel={t}, ylabel={$i_{\surf}$},legend style={ at={(0.5,1.02)}, anchor=south, legend columns =2}, legend cell align=left]
		\addplot +[mark repeat = 32] table [x=surface_list_x, y=-, col sep=comma] {Experiments/Surface_Vector_lt_5_lx_3_Diss4_alpha_h_beta0_stabiwith_h_SpatialNormalVolume_userh_1_tmax_5_tstart_5_xmax_3_xstart_3_ks_1_kt_1_kgs_1_kgt_1_gp_0_tempquadorder_4_tempquadorderDiffu_4.csv};
		\addplot +[mark repeat = 64] table [x=surface_list_x, y=-, col sep=comma] {Experiments/Surface_Vector_lt_6_lx_4_Diss4_alpha_h_beta0_stabiwith_h_SpatialNormalVolume_userh_1_tmax_6_tstart_6_xmax_4_xstart_4_ks_1_kt_1_kgs_1_kgt_1_gp_0_tempquadorder_4_tempquadorderDiffu_4.csv};
		\addplot +[mark repeat = 128] table [x=surface_list_x, y=-, col sep=comma] {Experiments/Surface_Vector_lt_7_lx_5_Diss4_alpha_h_beta0_stabiwith_h_SpatialNormalVolume_userh_1_tmax_7_tstart_7_xmax_5_xstart_5_ks_1_kt_1_kgs_1_kgt_1_gp_0_tempquadorder_4_tempquadorderDiffu_4.csv};
		\addplot +[mark repeat = 256] table [x=surface_list_x, y=-, col sep=comma] {Experiments/Surface_Vector_lt_8_lx_6_Diss4_alpha_h_beta0_stabiwith_h_SpatialNormalVolume_userh_1_tmax_8_tstart_8_xmax_6_xstart_6_ks_1_kt_1_kgs_1_kgt_1_gp_0_tempquadorder_4_tempquadorderDiffu_4.csv};
\addplot[dashed, purple,line width=0.5pt] coordinates { 
	(0.87,35) (0.87,86)
};
\addplot[dashed, purple,line width=0.5pt] coordinates { 
	(1.55,35) (1.55,86)
};
\addplot[dashed, purple,line width=0.5pt] coordinates { 
	(0.61,35) (0.61,86)
};		
				\legend{{$l_\ti=l_{\spa}=2$},{$l_\ti=l_{\spa}=3$},{$l_\ti=l_{\spa}=4$},{$l_\ti=l_{\spa}=5$}, top. singularity}
		\end{axis}
		\end{tikzpicture}
		\caption{Surface area over time.}
		\label{Dis_sphere_surface} 
	\end{subfigure}
	\caption{Disintegrating sphere: the change of mass and surface area for $k=1$.}
	\label{Dis_sphere_over_time}
\end{figure}

As in \Cref{section_smooth_geometry} we run numerical experiments to assess if and how accurate the discrete solution satisfies a mass conservation property. In \Cref{Dis_sphere_over_time} we illustrate the evolution of mass and surface area over time.
In \Cref{Dis_sphere_surface} we see that after the first singularity at $t\approx0.61$ there is a strong increase of the surface area. The merging of the inner and outer surfaces, i.e. the formation of the decocube at $t\approx 0.87$, leads to a decrease of the surface area and the splitting at $t\approx 1.55$ does not change this decreasing  surface area behavior significantly. The mass on the inner sphere, which arises from a point singularity at $t\approx 0.61$, remains  zero until the  merging with the outer surface at $t\approx 0.87$. 
When the two surfaces merge  we observe  an abrupt change of mass, which gets significantly smaller for finer refinement levels. The third topological singularity, the splitting into eight smaller connected surfaces at $t\approx 1.55$, yields a (large) mass loss, even though the surface area does not change significantly at that time. This mass loss becomes smaller on finer levels. The mass function $i_{\rm mass}(t)$ seems to converge to a constant function, which is the expected behavior.

Finally we note one further robustness aspect. The results above show that both in the example of the merging spheres and of the disintegrating sphere one obtains stable and reasonably accurate results on very coarse meshes, e.g., level $l=2$. This indicates that with this discretization method one preserves a main favourable property of the level set technique, namely that it can handle topological singularities in a very stable way.

\section{Conclusions}
In this paper we studied a Eulerian finite element method for the full discretization of scalar partial differential equations on evolving surfaces. The method uses standard space-time finite element spaces on a bulk mesh in combination with  the trace technique known from the literature. The method is applied to the  model problem
 \cref{surfactant1} which describes diffusive transport on an evolving surface. For higher order accuracy we use a space-time parametric mapping. This leads to the discrete problems \eqref{discreteproblem}, parameterized with $\beta \in [0,1]$.  Key properties and implementation issues of this discretization method have been addressed. The performance of the method is illustrated in examples with one and two-dimensional surfaces. The results show that the method has (optimal) higher order convergence (in space and time) in cases with smoothly evolving surfaces and can handle topological singularities in a robust way. 
\\[3ex]
{\bf Acknowledgement.} The authors thank the German Research Foundation
(DFG) for financial support within the Research Unit ”Vector- and tensor valued
surface PDEs” (FOR 3013) with project no. RE 1461/11-2.

\bibliographystyle{siam}
\bibliography{literatur}

\appendix
\section{Proof of \Cref{PI_Theorem}}\label{proof_PI}
\begin{proof}
	By e.g. \cite{GT} we know the partial integration rule for a smooth function $g$ defined on an arbitrary $K_S\in \mathcal{T}_{S_h^n}$
	\begin{equation*}
		\int_{K_S}\nabla_{S_h,i} g \dif \sigma_h=\int_{K_S} g \kappa_{K_S} (\bn_{S_h})_i \dif \sigma_h+\int_{\partial K_S} g (\bnu_{h})_i \dif F,\quad i\in \{1,\dots,4\}.
	\end{equation*}
	Here, $\nabla_{S_h,i}g$ denotes the $i$-th Cartesian component of $\nabla_{S_h}g$. The function $\kappa_{K_S}$ denotes the mean curvature on $K_S$. Replacing $g$ by $\frac{uv}{\alpha_h}\left(\bP_{S_h}\bw_S\right)_i$, summing over $i$ and using the product rule we get
	\begin{align}
		\begin{aligned}
			&\int_{K_S} u \Div_{S_h}\left(\frac{v}{\alpha_h}\bP_{S_h}\bw_S\right)+\left(\frac{v}{\alpha_h}\bP_{S_h}\bw_S\right)\cdot \nabla_{S_h}u\dif \sigma_h\\
			&\qquad=\int_{K_S}\kappa_{K_S} \frac{uv}{\alpha_h} \left(\bP_{S_h}\bw_S\right)\cdot \bn_{S_h}\dif \sigma_h+\int_{\partial K_S}\frac{hv}{\alpha_h}\left(\bP_{S_h}\bw_S\right)\cdot \bnu_{h} \dif F.\label{General_PI}
		\end{aligned}
	\end{align}
	The first term on the right hand side vanishes, since the vectors in the dot product are orthogonal. 
	We calculate using \cref{General_PI} and the definition \eqref{weakmatderivative} of the discrete material derivative
	\begin{align*}
		&\int_{S_h^n}\frac{1}{\alpha_h}\mathring{u}v\dif \sigma_h=\sum_{K_S\in  \mathcal{T}_{S_h^n}}\int_{K_S}\frac{v}{\alpha_h}\left(\bP_{S_h}\bw_S\right)\cdot \nabla_{S_h}u\dif \sigma_h\\
		&\quad=\sum_{K_S\in  \mathcal{T}_{S_h^n}}\left[-\int_{K_S}  u\Div_{S_h}\left(\frac{v}{\alpha_h} \bP_{S_h}\bw_S\right)\dif \sigma_h+\int_{\partial K_S}\frac{uv}{\alpha_h} \left(\bP_{S_h}\bw_S\right)\cdot \bnu_{h}\dif F \right].
	\end{align*}
	We continue using the product rule and the tangency of $\bnu_{h}\restrict{K_S}$ to $K_S$
	\begin{align*}
		\int_{S_h^n}\frac{1}{\alpha_h}\mathring{u}v\dif \sigma_h  &  =\sum_{K_S\in  \mathcal{T}_{S_h^n}}\left[-\int_{K_S} uv\Div_{S_h}\left( \frac{1}{\alpha_h}\bP_{S_h}\bw_S\right) \right. \\ & \quad -\left. \frac{u}{\alpha_h} \left(\bP_{S_h}\bw_S\right)\cdot \nabla_{S_h}v\dif \sigma_h
		+\int_{\partial K_S}\frac{uv}{\alpha_h}\bw_S\cdot \bnu_{h}\dif F\right].
	\end{align*}
	For the second term on the right hand side we use the tangency of $\bP_{S_h}\bw_S$ to $S_h$ and  the definition of the weak material derivative \cref{weakmatderivative}. We split the sum in the third term on the right-hand side. In case $\partial K_S \subset \partial S_h^n$ we know $\bnu_{h}=\bnu_{\partial }$ at the top boundary of the time slab and $\bnu_{h}=-\bnu_{\partial }$ at the bottom boundary of the time slab. We obtain
	\begin{align*}
		&\int_{S_h^n}\frac{1}{\alpha_h}\mathring{u}v\dif \sigma_h=-\int_{S_h^n}\frac{1}{\alpha_h}u\mathring{v}\dif \sigma_h+ \sum_{F\in \mathcal{F}_T^n}\int_F \frac{u_{-}^{n}v_{-}^{n}}{\alpha_{h,-}^n}\bw_S\cdot(\bnu_{\partial })_-^n\dif F\\
		&\qquad\ -\sum_{F\in \mathcal{F}_B^n}\int_F \frac{u_{+}^{n-1}v_{+}^{n-1}}{\alpha_{h,+}^{n-1}}\bw_S\cdot(\bnu_{\partial })_+^{n-1}\dif F+ \sum_{F\in \mathcal{F}_I^n}\int_{F} u v \bw_S\cdot \left[\frac{1}{\alpha_h}\right]_{\bnu}\dif F\\
		& \qquad-\sum_{K_S\in  \mathcal{T}_{S_h^n}}\int_{K_S} uv\Div_{S_h}\left( \frac{1}{\alpha_h} \bP_{S_h}\bw_S\right)\dif \sigma_h.
	\end{align*}
	Due to the definition of $\mathcal{F}_T^n$ and $\mathcal{F}_B^n$ we get 
	\begin{align*}
		&\int_{S_h^n}\frac{1}{\alpha_h}\mathring{u}v\dif \sigma_h=-\int_{S_h^n}\frac{1}{\alpha_h}u\mathring{v}\dif \sigma_h +\int_{\Gamma^n_h(t_n)}\frac{u_{-}^{n}v_{-}^{n}}{\alpha_{h,-}^n}\bw_S\cdot (\bnu_{\partial })_-^n\dif s_h\\
		&\qquad-\int_{\Gamma^n_h(t_{n-1})}\frac{u_{+}^{n-1}v_{+}^{n-1}}{\alpha_{h,+}^{n-1}}\bw_S\cdot (\bnu_{\partial })_+^{n-1}\dif s_h  \\
		& \qquad+ \sum_{F\in \mathcal{F}_I^n}\int_{F} uv\bw_S\cdot \left[\frac{1}{\alpha_h}\right]_{\bnu} \dif F -\sum_{K_S\in  \mathcal{T}_{S_h^n}}\int_{K_S} uv\Div_{S_h}\left( \frac{1}{\alpha_h} \bP_{S_h}\bw_S\right)\dif \sigma_h.
	\end{align*}
	Finally, using the definition $R=\frac{1}{\alpha_h}\bw_S\cdot \bnu_{\partial}$ completes the proof.
\end{proof}
\end{document}